\documentclass[draft]{article}
 
\usepackage{amsmath,amsfonts,amsthm,amssymb,amscd}

\theoremstyle{plain} \newtheorem{theorem}{Theorem}[section]
\newtheorem{lemma}[theorem]{Lemma}

\newtheorem{corollary}[theorem]{Corollary} \theoremstyle{definition}
\newtheorem{definition}[theorem]{Definition} \theoremstyle{remark}
\newtheorem{remark}[theorem]{Remark}

\newcommand{\R}{{\mathbb R}} \newcommand{\U}{{\mathcal U}}

\newcommand{\Z}{{\mathbb Z}}

\newcommand{\N}{{\mathbb N}}

\newcommand{\Ph}{{\mathcal P}}

\newcommand{\resto}{{\mathcal R}}
\def\im{{\rm i}}
\def\rmi{{\rm i}}

\newcommand{\C}{\mathbb{C}}

\def\uno{{\kern+.3em {\rm 1} \kern -.22em {\rm l}}}

\numberwithin{equation}{section}

\begin{document}

  \title{On the Darboux and Birkhoff   steps in the asymptotic  stability of    solitons}

 \author {Scipio Cuccagna}

 \maketitle
\begin{abstract} We   give a unified proof of the step to find   Darboux coordinates and of the ensuing Birkhoff normal forms  procedure,
developed in the course of the proof of  asymptotic stability
of  solitary  waves in \cite{boussaidcuccagna,Cu1,Cu2}.
\end{abstract}

\section{Introduction}

The aim of this paper is to extend in a slightly  more general
and unified
set up two important steps of  the proof of the asymptotic stability
of solitary waves for the Nonlinear Schr\"odinger equation   \cite{Cu2,Cu1,bambusi}
and the particular case of  Nonlinear Dirac system treated  in \cite{boussaidcuccagna}. In both cases
there is a localization at the solitary wave
and a representation of   the system in terms of coordinates arising
from the linearization    at a  solitary wave.
The operators $\mathcal{H}_p$  introduced later play this role.
In general $\mathcal{H}_p$ has both continuous spectrum and non zero eigenvalues. The latter give rise to discrete modes which in  the nonlinear problem could produce chaotic Lissaius like motions. It turns out that in
\cite{bambusicuccagna,Cu2,Cu4,boussaidcuccagna,Cu1,bambusi}
discrete modes relax to  0
  because of a mechanism
of slow leaking of energy away from the discrete modes into the continuous modes, where energy disperses by linear dispersion. The idea was initiated in special situations in \cite{sigal,BP2,SW3}. We refer to
\cite{Cu2} for more comments and references.

The aim of this paper consists in simplifying  two key steps
in the proofs in \cite{Cu1,Cu2,boussaidcuccagna}. The  first step
consists in searching Darboux coordinates. This  allows to decrease the
number of coordinates in the system and to reduce to the study of
the system at an equilibrium point.

The second step consists in the
implementation of the Birkhoff normal forms, to produce a simple {\textit{effective}} Hamiltonian.
After this,  \cite{Cu1,Cu2,boussaidcuccagna}    prove  the  energy
leaking     away from the discrete modes.  In particular the key step is the proof  that certain coefficients
of the discrete modes equations are second powers, the  { \it Nonlinear Fermi Golden Rule} (FGR),
which generically are positive and  yield  discrete mode energy dissipation.

 We do not discuss the  FGR  in this paper    limiting ourselves to the search of Darboux coordinates
and to the Birkhoff normal forms argument.

In this paper we avail ourselves with some ideas and notation  drawn from  early versions of \cite{bambusi} to improve the
presentation in  \cite{Cu1}.

 \cite{bambusi,Cu1}  represent  two   attempts to
 extend the result proved in \cite{Cu2} for standing ground states of the NLS, to the case of moving ground states. A  further   goal  in \cite{bambusi}  is to develop the theory in a more abstract set up.    Early versions of    \cite{bambusi}  did not encompass  a  Birkhoff  step extendable to
 \cite{boussaidcuccagna}.
 \cite{bambusi}   is confined (like us here)  to systems with
	Abelian group of symmetries.

Our present proof was written  before the 3rd version of  \cite{bambusi} was posted on the Arxiv site.
The 2rd version of \cite{bambusi} contained an incorrect effective Hamiltonian,
see Remark \ref{rem:formal2} later.     In the 3rd version of \cite{bambusi} this has been corrected. Still, the discussion in \cite{bambusi}   is  at times sketchy,
for example in Theorems    3.21  and  5.2 
\cite{bambusi}, see Remarks \ref{rem:formal1} and \ref{rem:bamb1}  and the discussion below and at the beginning of Sect.  \ref{subsec:darboux}
about  \eqref{eq:fdarboux}.

 We draw from  \cite{bambusi}
	a better choice of initial  coordinates and set up than \cite{Cu1}.  Some of it  existed also in previous literature, cfr.
	the discussion in Sect.6 \cite{RSS}.
	We also borrow some notation, i.e. symbols $\resto ^{k,m}$ and $\textbf{S} ^{k,m}$.
 Inspired by \cite{bambusi} we simplify the proof of the Birkhoff step in \cite{Cu1}.

      Both
here and in \cite{Cu1} we
consider  initial data in subsets of $\Sigma _n$ for $n\gg 1$ which are
unbounded  in  $\Sigma _n$   and invariant for the system.
We require this  substantial amount of  regularity and spacial decay to 0  for the classes
of solutions of the system,  in order to give a rigorous
treatment of the flows and of the pullbacks.   \cite{bambusi} suggests that
 \cite{Cu1}  should prove decay rates in time. 
We do not know what is the basis for this statement in \cite{bambusi} since,     by the  time invariance of the subsets
$\Sigma _n$ considered, the problem considered in  
 \cite{Cu1}  is very similar in this respect to the one with $\Sigma _n$ replaced by $H^1$.   
Indeed time decay corresponds to bounds on   norms containing time dependent weights.
But if the   problem is invariant by translation in time, the only information that can be 
derived must be invariant by  translation in time, and  bounds on  time weighted norms do not have this property.
We therefore emphasize  that \cite{Cu1}  and the present paper are very different from, say,  \cite{BP2,SW3}, which consider
initial data in subsets of   $H^{k,s}$ which are not invariant by the time evolution.     We also point out that
the 1st version of   \cite{bambusi} contains a false statement on rate of time decay.  The 
2nd version of   \cite{bambusi} in the acknowledgments,  credits us for pointing out this error, although these credits are not any more in  the 3rd version.

To find an effective Hamiltonian, we use in a crucial way the regularity
properties of the flows, which in turn depend  on the fact that we work in  $\Sigma _n$ for $n\gg 1$.
See Theorem \ref{theorem-1.1} where the regularity of the flows is  used to prove that  the coordinate
changes preserve the  system.
To prove for the NLS the same result    in $H^1$,  where the
coordinate changes are continuous only, one needs to explain how they
 preserve the
structure needed to make sense of the NLS.   \cite{bambusi}  claims the result in $H^1$,
	 \cite{bambusi}  claims the result in $H^1$,
	but the  proof  is not spelled out,
see Remark \ref{rem:bamb1}.

			We discuss in some detail
			a key formula on the differentiation   of the pullback of
			a differential form
			along   a flow, see  \eqref{eq:fdarboux}, which is the   basis of
			Moser's method to find Darboux coordinates.
This formula is simple in
			classical set ups, but  in our case  and in \cite{bambusi}  its interpretation  and proof are not obvious.
				In   \cite{bambusi} the formula is stated and used without discussion.
				We treat the issue  rigorously    in
					Sect. \ref{subsec:darboux},    regularizing the flow, using
 \eqref{eq:fdarboux} for the regularized flow,
 and   recovering the  desired equality between differential forms, by
 a  limiting argument. Notice that we do not prove  formula \eqref{eq:fdarboux} for the non regularized flow.

We end with few remarks on the proofs.

 The proof of the Darboux Theorem is   a simplification of that in   \cite{Cu1} in the part discussing
the vector field.  We give in Sect. \ref{subsec:darboux} a detailed proof on the fact that the resulting flow transforms the symplectic
form as desired. See also the introductory remarks in Sect. \ref{sect:symplectic}. Notice that parts of this discussion were skipped in \cite{Cu1}.

 The portion of our paper  on  the Birkhoff normal forms covers from Sect. \ref{sec:pullback}
on and
is quite different from  \cite{boussaidcuccagna,Cu1,Cu2}
mainly because the pullback of the terms of the expansion of the Hamiltonian
cannot be treated on a term by term basis,
see Remark \ref{rem:differences}.  What is important is to get a general structure
of the pullbacks of the Hamiltonian. This is discussed in Sect. \ref{sec:pullback}.
It is likely
that most of the analysis in Lemmas \ref{lem:back},  \ref{lem:back11} and \ref{lem:ExpH11}, is not necessary to the derivation of the effective
Hamiltonian, which is represented by $H'_2$ and the null terms in
$\mathbf{R}_0$ and $\mathbf{R}_1$ of the expansion in Lemma  \ref{lem:ExpH11}, in the final Hamiltonian.
On the other hand, writing the Hamiltonian explicitly should make the arguments   transparent and more clearly
applicable to the part on dispersion and Fermi Golden rule.

In Sect. \ref{sec:speccoo}  we finally distinguish between discrete and continuous modes.

 The present paper treats only equations whose symmetry group is Abelian.  This  limitation will have to be
overcome to extend the theory   to more general systems such for example
  the Dirac system without the symmetry  constraints of \cite{boussaidcuccagna}.

\section{Set up}
\label{sec:setup}

\begin{itemize}
\item
Given two vectors  $ {u},v\in \R ^{2N}$ we denote by $u\cdot v =\sum u_j v_j$
their inner product.

\item We will consider also another quadratic form $|u|_1^2= u\cdot _1 u$
in  $  \R ^{2N}$.

	\item For any   $n \ge 1$ we consider the space $\Sigma _n=\Sigma _n(\R ^3, \R ^{2N})$ defined by
 \begin{equation}\label{eq:sigma}
\begin{aligned} &
      \| u \| _{\Sigma _n} ^2:=\sum _{|\alpha |\le n}  (\| x^\alpha  u \| _{L^2(\R ^3, \R ^{2N})} ^2  + \| \partial _x^\alpha  u \| _{L^2(\R ^3, \R ^{2N})} ^2 )  <\infty .      \end{aligned}\nonumber
\end{equation}
	 We set $\Sigma _0= L^2(\R ^3, \R ^{2N})$. Equivalently we can define $\Sigma _{r}$ for $r\in \R$  by the norm
 \begin{equation}
\begin{aligned} &
      \| u \| _{\Sigma _r}  :=  \|  ( 1-\Delta +|x|^2)   ^{\frac{r}{2}} u \| _{L^2}   <\infty    .   \end{aligned}\nonumber
\end{equation}
 For $r\in \N$ the two definitions are equivalent, see    \cite{Cu2}.	 We will not use another
quite natural class of spaces denoted by $H^{k,s}$ and defined by

\begin{equation}
\begin{aligned} &
      \| u \| _{H^{k,s}}  :=  \|   (1+|x|^2) ^{\frac{s}{2}}    ( 1-\Delta  )   ^{\frac{k}{2}} u \| _{L^2}   <\infty    .   \end{aligned}\nonumber
\end{equation}

\item ${\mathcal S}(\R ^3, \R ^{2N})=\cap _{n\in \N}\Sigma _n(\R ^3, \R ^{2N})$
is  the space of Schwartz functions  and the space of tempered distributions is
${\mathcal S}'(\R ^3, \R ^{2N})=\cup _{n\in \N}\Sigma _{-n}(\R ^3, \R ^{2N})$.

\item  For $X$ and $Y$ two Banach space, we will
denote by $B(X,Y)$ the Banach space of bounded linear operators from
$X$ to $Y$ and by $B^{\ell}(X,Y)= B (  \prod _{j=1}^\ell X ,Y)$.

\item    We denote by $\langle \ , \ \rangle $
the natural  inner product in $L^{2}(\R ^3, \R ^{2N})$.

 \item     $J$ is an invertible antisymmetric matrix in
$\R ^{2N}$. We have also $|J y|_1=|  y|_1$ for any $y\in \R ^{2N}$.
 In $L^{2}(\R ^3, \R ^{2N})$ we consider   the symplectic form $\Omega =\langle  J ^{-1}\ , \ \rangle
$.

 \item We consider in $ L^2(\R ^3, \R ^{2N})$ a linear    selfadjoint elliptic differential operator $\mathcal{D}$
such that $\mathcal{D}\in B(\Sigma _{r},\Sigma _{r-\text{ord}\mathcal{D}} ) $ and  $\mathcal{D} \in B(H^{r},H^{r-\text{ord}\mathcal{D}} ) $ for all $r$ and for a  fixed integer $\text{ord}\mathcal{D}\ge 1$.

\item
We   consider a Hamiltonian  of the form

 \begin{equation} \label{eq:energyfunctional}\begin{aligned}&
 E(U)=E_K(U)+E_P(U)\\&
E_K(U):=\frac{1}{2} \langle  {\mathcal D} U, U \rangle  \, , \quad
E_P(U):=
 \int _{\R ^3}B( |U|_1^2) dx .\end{aligned}
\end{equation}
Here   $B\in C^\infty (\R , \R)$, $B( 0)=B'( 0)=0$ and there exists a $p\in(2,6]$ such that for every
$k\ge 0$ there is a fixed $C_k$ with
\begin{equation}\label{eq:growthB}  \left|  \nabla _\zeta  ^k (B ( |\zeta |_1^2) )\right|\le C_k
|\zeta |^{p-k-1} \quad\text{if $|\zeta |\ge 1$  in $\R ^{2N}$.}
\end{equation}

\end{itemize}
Notice that    $E_P\in C^5(H ^{1}(\R ^3, \R ^{2N})), \R)$.
Consistently with \cite{boussaidcuccagna,Cu1,Cu2}, we focus only on \textit{semilinear} Hamiltonians.  We consider the system
\begin{equation}\label{eq:NLSvectorial}   \dot U = J \nabla E (U)  \quad ,   \quad U(0)=U_0\quad
\end{equation}
where for a Frech\'et differentiable function $F$  the gradient $\nabla F (U)$
is defined by
$\langle \nabla F (U), X\rangle = dF(U) (X)$, with $dF(U)$ the exterior differential calculated at $U$.
We assume that

\begin{itemize}
\item[(A1)]  there exists $d_0$ such that for $d>d_0$ system \eqref{eq:NLSvectorial} is locally well posed in
$H^d$. Furthermore, the space $\Sigma _d$  is invariant by this motion.

\end{itemize}

\noindent We recall the following   definition.

\begin{definition}\label{def:HamField} Given a Frech\'et differentiable function $F$,
the Hamiltonian vectorfield of  $F$ with respect to a \textit{strong} symplectic form $\omega $,  see \cite{amr} Ch. 9, is  the field $X_F $ such that $ \omega (X_F ,Y)=  dF (Y)$ for any given
tangent vector  $Y.$ For $\omega =\Omega$ we have $X_F=J\nabla F$.

For    $F,G  $  two scalar Frech\'et differentiable functions,
we consider the   Poisson bracket
$
  \{ F,G \}  :=  dF (X_G ).
$

If $\mathcal{G}  $   has values in a given
Banach space $\mathbb{E}$ and $G$ is a scalar valued function, then we set  $
  \{ \mathcal{G} ,G \}  :=  \mathcal{G} '(X_G ),
$ for $\mathcal{G}'$ the Frech\'et derivative of $\mathcal{G} $.
\end{definition}
We assume some symmetries in system \eqref{eq:NLSvectorial}. Specifically
we assume what follows.

\begin{itemize}
\item[(A2)] There are    selfadjoint  differential  operators  $\Diamond  _\ell$
for  $\ell =1,..., n_0$ in $L^2$ such that
$ \Diamond  _\ell :
\Sigma _{n} \to \Sigma _{n-d_\ell}   $ for  $\ell =1,..., n_0$.
  We set   $\textbf{d} =\sup _\ell d_\ell$.

\item[(A3)]We assume
$[\Diamond  _\ell  , J] =0$ and $[\Diamond  _\ell  , \Diamond  _k] =0$.

\item[(A4)] We assume   $\{\Pi _\ell , E_K \} = \{\Pi _\ell , E_P \} =0$  for  all $\ell $, where   $\Pi _\ell := \frac{1}{2} \langle \Diamond  _\ell  \  ,
\ \rangle $.

\item[(A5)] Set $
\langle \epsilon \Diamond\rangle ^{ 2}:=1+\sum _j \epsilon ^2 \Diamond ^2_j.$  Then
$ \langle \epsilon \Diamond\rangle ^{ -2}  \in B(\Sigma  _n ,\Sigma  _ n )$ with
\begin{equation} \label{eq:a71}  \begin{aligned}
 &
\quad \text{ $\|     \langle \epsilon \Diamond\rangle ^{ -2}  \| _{B(\Sigma  _n ,\Sigma  _ n ) } \le C_n<\infty$ for any $|\epsilon |\le 1 $ and $n\in \N$.}
 \end{aligned}
\end{equation}
Furthermore, for any $n\in \Z$ we have
\begin{equation} \label{eq:a72}   \begin{aligned}
 &    strong-\lim _{\epsilon\to 0}  \langle \epsilon \Diamond\rangle ^{ -2} =1 \text{  in $B(\Sigma  _n ,\Sigma  _ n )$ } \\&  \lim _{\epsilon\to 0}  \| \langle  \epsilon \Diamond\rangle ^{ -2} -1 \|  _{B(\Sigma  _n ,\Sigma  _ {n'} ) }=0
\quad \text{  for any   $n'\in \Z$ with $n'<n$. }
 \end{aligned}
\end{equation}

 \item[(A6)]  Consider the
  groups
$e^{J \langle  \epsilon \Diamond\rangle ^{ -2}\Diamond \cdot \tau  } $   defined in $L^2$.
We assume that  for any $n\in \N$ these groups leave  $ \Sigma _n$ invariant and  that for any  $n\in \N$ and $c>0$ there a $C$ s.t.
$\| e^{J\langle  \epsilon \Diamond\rangle ^{ -2}\Diamond \cdot \tau  } \| _{B(\Sigma _n,\Sigma _n)}\le C$ for any $|\tau |\le c$  and  any $|\epsilon |\le 1 $ .\end{itemize}

We introduce now our  \textit{solitary} waves.

\begin{itemize}
\item[(B1)]  We assume   that for  $  {\mathcal O}$ an open subset of $\R ^{n_0}$
we have a function  $p\to  \Phi _p  \in {\mathcal S}(\R ^3, \R ^{2N}) $ which is   in
$C^\infty ( {\mathcal O} ,  {\mathcal S} )$, with $\Pi _\ell (\Phi _p )=p_\ell$, where the
$\Phi _p   $ are constrained critical points of $E$ with associated Lagrange multipliers $\lambda _\ell (p)$
so that

 \begin{equation} \label{eq:lagr mult}\begin{aligned}&
 \nabla E(\Phi _p  )= \lambda   (p) \cdot \Diamond \Phi _p \end{aligned}
\end{equation}

\item[(B2)]   We will assume that the map $p\to \lambda   (p)$ is a diffeomorphism.
In particular this means that the following matrix has rank $n_0 $

\begin{equation} \label{eq:nondegen} \begin{aligned} & \text{rank}  \left [ \frac{\partial \lambda _i}{\partial p _j}   \right ]      _{ \substack{ i\downarrow  \  , \  j \rightarrow}}= n_0 .
\end{aligned}
\end{equation}

\end{itemize}

A function    $U(t):=  e^{J( t   \lambda (p) +\tau _0)\cdot \Diamond}\Phi _p  $
		 is a  {solitary wave} solution of \eqref{eq:NLSvectorial}
 for any fixed  vector $\tau _0$.

\subsection{The linearization}
\label{subsec:linearization}
		
		Set ${\mathcal H}_p:= J(\nabla ^2 E(\Phi _p  )- \lambda  (p) \cdot \Diamond ) $.
		Notice that   $E (e^{J \tau   \cdot \Diamond}U) \equiv E ( U)  $ for any $U$ yields
		$\nabla E (e^{J \tau   \cdot \Diamond}U) = e^{J \tau   \cdot \Diamond} \nabla E ( U)  $ and
		$\nabla ^2 E (e^{J \tau   \cdot \Diamond}U) = e^{J \tau   \cdot \Diamond} \nabla ^2E ( U) e^{-J \tau   \cdot \Diamond} .$
		Then
		\eqref{eq:lagr mult}    implies $\nabla E( e^{J \tau   \cdot \Diamond}\Phi _p  )= e^{J \tau   \cdot \Diamond}  \lambda  (p) \cdot \Diamond \Phi _p $. So applying $\partial _{ \tau _j}$ we obtain
		$ (\nabla ^2 E(  \Phi _p  )  -\lambda   (p) \cdot \Diamond )   J \Diamond _j    \Phi _p  =0$ and so
 \begin{equation} \label{eq:kernel}\begin{aligned}&
 {\mathcal H}_p   J \Diamond _j    \Phi _p  =0\end{aligned}
\end{equation}
\begin{itemize}
\item [(C1)]  We will assume

\begin{equation} \label{eq:kernel1}\begin{aligned}&\ker
 {\mathcal H}_p  =\text{Span}\{ J \Diamond _j    \Phi _p:j=1,..., n_0   \} .\end{aligned}
\end{equation}
\end{itemize}

	\noindent Applying    $\partial _{ \lambda _j}$    to \eqref{eq:lagr mult}  yields
		$ (\nabla ^2 E(  \Phi _p  )  -\lambda  (p) \cdot \Diamond )     \partial _{ \lambda _j} \Phi _p  =   \Diamond _j\Phi _p  $. This yields
		
		\begin{equation} \label{eq:gen kernel}\begin{aligned}&
 {\mathcal H}_p   \partial _{ \lambda _j}     \Phi _p  =J \Diamond _j    \Phi _p\end{aligned}
\end{equation}
We have
\begin{equation} \label{eq:nondegen1} \begin{aligned} &  \langle \partial _{ \lambda _j}     \Phi _p  , \Diamond _k    \Phi _p\rangle =\frac 12 \partial _{ \lambda _j}  \langle    \Phi _p  , \Diamond _k    \Phi _p\rangle =\partial _{ \lambda _j}p_k.
\end{aligned}
\end{equation}
 Necessarily, by (B2) there exists $j$ such that $\partial _{ \lambda _j}p_k\neq 0.$
 This implies that the \textit{generalized kernel} is
\begin{equation} \label{eq:gen kernel1}\begin{aligned}&N_g (
 {\mathcal H}_p ) =\text{Span}\{ J \Diamond _j    \Phi _p, \partial _{ \lambda _j}     \Phi _p :j=1,..., n_0   \} .\end{aligned}
\end{equation}

		\noindent The map $(p , \tau ) \to e^{J \tau _0 \cdot \Diamond}\Phi _p $  is in
		$C^\infty ( {\mathcal O} \times  \R ^{n_0 },  {\mathcal S} )$.
\begin{itemize}
\item [(C2)]
We assume this map is a local embedding
		and that the image is a manifold $\mathcal{G}$.
\end{itemize}
At any given point  $e^{J \tau   \cdot \Diamond}\Phi _{p }$
		the tangent space of $\mathcal{G}$ is  given by
		
\begin{equation*} \label{ }\begin{aligned}&  T_{e^{J \tau   \cdot \Diamond}\Phi _{p }} \mathcal{G} = \text{Span}\{ e^{J \tau   \cdot \Diamond}\partial _{p_j}\Phi _{p },
e^{J \tau   \cdot \Diamond}\Diamond  _{j }\Phi _{p }  :j=1,..., n_0\}  . \end{aligned}
\end{equation*}
We  have $\Omega (   e^{J \tau   \cdot \Diamond}\partial _{p_j}\Phi _{p },  e^{J \tau   \cdot \Diamond}\partial _{p_k}\Phi _{p } ) = \Omega (    \partial _{p_j}\Phi _{p },   \partial _{p_k}\Phi _{p } )$.

\begin{itemize}
\item [(C3)] We
assume  that
 \begin{align} \label{eq:gen kernel2} &  \Omega (    \partial _{p_j}\Phi _{p },   \partial _{p_k}\Phi _{p } ) = 0
\text{ for all $j$ and $k$}  \\&   \label{eq:gen kernel20} \Omega (    \partial _{p_j}\Phi _{p },    \Phi _{p } ) = 0
\text{ for all $j$} . \end{align}
\end{itemize}
Notice that \eqref{eq:gen kernel20}  is not required in \cite{bambusi} but in any case is true for the applications
in \cite{bambusi,boussaidcuccagna,Cu1,Cu2}. Here we use it   in Lemma \ref{lem:1forms}.

We have  the following beginning of Jordan block decomposition of $\mathcal{H}_p$.
\begin{lemma}
  \label{lem:begspectdec}  Consider the operator $\mathcal{H}_p$. We
  have
  \begin{equation} \label{eq:begspectdec1} \begin{aligned} &
 J^{-1}\mathcal{H}_p=-\mathcal{H}_p^{\ast}J^{-1} \, , \quad  \mathcal{H}_pJ=-J\mathcal{H}_p^{\ast} .
\end{aligned}\end{equation}
Assume (B1)--(B2) and (C1). Then we have

\begin{align}
	 \label{eq:begspectdec2}& L^2= N_g(\mathcal{H}_p)\oplus N_g^\perp (\mathcal{H}_p^{\ast}) \  ,
   \\& N_g (\mathcal{H}_p^{\ast})  =\text{Span}\{   \Diamond _j    \Phi _p, J^{-1}\partial _{ \lambda _j}     \Phi _p :j=1,..., n_0   \}   .  \label{eq:begspectdec3}
\end{align}

\end{lemma}
\proof We have $\mathcal{H}_p=JA$ for a selfadjoint operator $A$ and with
$J$ a bounded antisymmetric operator. Then  $\mathcal{H}_p^{\ast}=-AJ$
and \eqref{eq:begspectdec1} follows by direct inspection.
Recall that (B1)--(B2) and (C1)  imply \eqref{eq:gen kernel1}.
Then
\eqref{eq:begspectdec1}  implies
 \eqref{eq:begspectdec3}.

\noindent   The map $\psi \to \langle \ , \psi \rangle $ establishes a map $N_g  (\mathcal{H}_p^{\ast})
\to B (N_g(\mathcal{H}_p), \R )$. By \eqref{eq:nondegen1}, formulas \eqref{eq:gen kernel1} and \eqref{eq:begspectdec3}
imply that this map is an isomorphism.   For any  $u \in L^2$ there is exactly one $v\in N_g(\mathcal{H}_p)$ such that $\langle u , \ \rangle$ and $\langle v , \ \rangle$ coincide as elements in $B (N_g  (\mathcal{H}_p^{\ast}) , \R )$.  Then  $ u-v \in  N_g^\perp (\mathcal{H}_p^{\ast}) $ and we get
   \eqref{eq:begspectdec2}.

\qed

Obviously Lemma \ref{lem:begspectdec} holds true only because  our $J$ is very special. For the KdV, where $J=\frac{\partial}{\partial x}$,
\eqref{eq:begspectdec2}--\eqref{eq:begspectdec3} are not true.

\noindent   
Denote by $P_{N_g }(p) =P_{N_g(\mathcal{H}_p)}$ the projection onto
$N_g(\mathcal{H}_p) $ associated to \eqref{eq:begspectdec2}  and  by  $P(p):=1-P_{N_g }(p)$   the projection on $N_g^\perp (\mathcal{H}_p^{\ast})$. We have, summing
on repeated indexes,

\begin{equation} \label{eq:projNg} \begin{aligned} & P_{N_g }(p)X= -J\Diamond
_j \Phi _p\  \langle X ,J^{-1} \partial _{p_j}\Phi _p\rangle +\partial
_{p_j}\Phi _p \  \langle X ,\Diamond _j \Phi _p \rangle. \end{aligned}
\end{equation}

\begin{lemma} 
  \label{lem:projections} Assume   (B1)--B(2) and (C1). 
	Then:
	
	\begin{itemize}
\item[(1)]   $P_{N_g }(p)\in  B (\mathcal{S}',\mathcal{S}) $ for any $p\in  \mathcal{O}$ and 
$P_{N_g }(p)\in C^{\infty}
(\mathcal{O}, B (\Sigma _{-k},\Sigma _{ k}))$ for any $k\in \N $.
	\item[(2)]   For any $p_0\in  \mathcal{O}$ and $k $ there exists an $\varepsilon _k>0$
	such that for $|p-p_0|< \varepsilon _k $

	\begin{equation}\label{eq:projections}
		P (p )P ( p_0) :N_g^\perp (\mathcal{H} _{p_0 }^{\ast})\cap \Sigma _{k}\to N_g^\perp (\mathcal{H}_{p }^{\ast})\cap \Sigma _{k}
	\end{equation}
	is  an isomorphism.
	
	\item[(3)]  For $h >k$ we have $\varepsilon _h \ge \varepsilon _k $.

	\end{itemize} 
\end{lemma}
\proof  Claim (1)  is elementary and we skip the proof.   
 
\noindent Consider the map  $P (p )P ( p_0) P(p)= 1 +P (p )(P_{N_g }(p)  -P_{N_g }(p_0) ) P(p) $  from   
$N_g^\perp (\mathcal{H}_{p }^{\ast})\cap \Sigma _{k }$ into itself.  By   Claim (1) and by the Fredholm alternative,  this is an isomorphism for  $|p-p_0|< \varepsilon _k $  with  $\varepsilon _k>0$ sufficiently small.
This implies that the  $P (p )P ( p_0)$ in \eqref{eq:projections} is onto. For the same reasons also 
$P (p_0 )P ( p ) P(p_0) $   is an isomorphism from
$N_g^\perp (\mathcal{H}_{p_0 }^{\ast})\cap \Sigma _{k }$ into itself.  Then  $P (p )P ( p_0)$ in \eqref{eq:projections} is one to one. This yields  Claim (2).

\noindent For  $h>k$ we have the commutative diagram

\begin{equation}
\begin{aligned}  
N_g^\perp (\mathcal{H} _{p_0 }^{\ast})&\cap \Sigma _{h}
\stackrel{P (p )P ( p_0)}{\rightarrow} &N_g^\perp (\mathcal{H}_{p }^{\ast})\cap \Sigma _{h} \\& \downarrow  & \downarrow 
\\  N_g^\perp (\mathcal{H} _{p_0 }^{\ast})&\cap \Sigma _{k}\stackrel{P (p )P ( p_0)}{\rightarrow} &N_g^\perp (\mathcal{H}_{p }^{\ast}) \cap \Sigma _{k}
\end{aligned}\nonumber
\end{equation}
with the vertical maps two embedding. This implies that for  $|p-p_0|< \varepsilon _k $
we have $\ker  P (p )P ( p_0)=0$ in  $N_g^\perp (\mathcal{H} _{p_0 }^{\ast}) \cap \Sigma _{h}$.
To complete the proof of Claim (3), we need to show that given $u\in  N_g^\perp (\mathcal{H} _{p  }^{\ast}) \cap \Sigma _{h}$ and the resulting  $v\in  N_g^\perp (\mathcal{H} _{p _0 }^{\ast}) \cap \Sigma _{k}$
with $u=P (p )P ( p_0)v$, we have   $v\in    \Sigma _{h}$.
But this follows immediately from 
 \begin{equation*}   v= u  + (P_{N_g }(p)  -P_{N_g }(p_0) )v   \text{ where $u\in  \Sigma _{h}$ and $(P_{N_g }(p)  -P_{N_g }(p_0) )v\in \mathcal{S}$.}
\end{equation*}

\qed

We will denote   the inverse of  \eqref{eq:projections} by 
	\begin{equation}\label{eq:projections1}
		(P (p )P ( p_0))^{-1} :N_g^\perp (\mathcal{H} _{p }^{\ast})\cap \Sigma _{k}\to N_g^\perp (\mathcal{H}_{p_0 }^{\ast})\cap \Sigma _{k}.
	\end{equation}

We have the  following \textit{Modulation} type lemma.
	
\begin{lemma}[Modulation]
  \label{lem:modulation} Assume (A2), (B.1), (B.2),  (C.1) and (C.3). Fix $n \in \Z$, $n\ge 0$ 
	and
  fix $\Psi _0=e^{J \tau _0 \cdot \Diamond}\Phi _{p_0}$. Then $\exists$ a neighborhood $\U  $     in $ \Sigma _{-n}(\R ^3, \R ^{2N}) $   of $U_0$
    and   functions $p  \in C^\infty (\U  , \mathcal{O})$
	and $\tau \in C^\infty (\U  , \R ^{n_0   })$ s.t.  $p(\Psi _0)=p_0 $   and $\tau (\Psi _0) =\tau _0$
  and s.t.  $\forall U\in \U  $\begin{equation}\label{eq:anzatz}\begin{aligned} &
 U =   e^{J \tau \cdot \Diamond} (  \Phi _{p } +R)
 \text{  and $R\in N^{\perp}_g (\mathcal{H}_p ^*)$.}
\end{aligned}\end{equation}
\end{lemma}
\proof   Consider the  following  $2n_0$ functions:

\begin{equation} \label{eq:nals}\begin{aligned} & \mathcal{F}_j (U,p, \tau ):= \Omega (
U  -e^{ J \tau \cdot \Diamond}  \Phi _ p, e^{ J \tau \cdot \Diamond} \partial _{p_j}\Phi _p )
   \\& \mathcal{G}_j (U,p, \tau ) :=\Omega (
U  -e^{ J \tau \cdot \Diamond}  \Phi _ p, Je^{ J \tau \cdot \Diamond} \Diamond _j \Phi _p )    . \end{aligned}
\end{equation}
These functions  belong to 
$C^\infty ( \Sigma _{-n} \times \mathcal{O} \times \R ^{n_0}, \R )$.
We introduce the notation $R = e^{ -J \tau \cdot \Diamond}U-\Phi _ p $.
Notice that $R=0$ for $U=\Phi _ p$. Then

\begin{equation*} \label{ }\begin{aligned} &  \partial _{\tau _k}\mathcal{F}_j (U,p, \tau ) = \Omega (
e^{ J \tau \cdot \Diamond} R  ,e^{ J \tau \cdot \Diamond} J   \Diamond _k \partial _{p_j}\Phi _p ) - \Omega (  J   \Diamond _k   e^{ J \tau \cdot \Diamond}   \Phi _p,
  e^{ J \tau \cdot \Diamond}  \partial _{p_j}\Phi _p )
   \\& =  -  \langle  R  ,    \Diamond _k \partial _{p_j}\Phi _p   \rangle  -  \langle  \Diamond _k      \Phi _p ,
     \partial _{p_j}\Phi _p  \rangle   = -  \langle  R  ,    \Diamond _k \partial _{p_j}\Phi _p   \rangle  -  \frac{1}{2}\partial _{p_j} \langle  \Diamond _k      \Phi _p  ,
     \Phi _p  \rangle  \\&  = -  \langle  R  ,    \Diamond _k \partial _{p_j}\Phi _p   \rangle  -  \delta _{jk}  . \end{aligned}
\end{equation*}
By \eqref{eq:gen kernel2} we have
\begin{equation*} \label{ }\begin{aligned} &  \partial _{p _k}\mathcal{F}_j (U,p, \tau ) = \Omega (
e^{ J \tau \cdot \Diamond} R  ,e^{ J \tau \cdot \Diamond}  \partial _{p_k}  \partial _{p_j}\Phi _p ) - \Omega (  J       e^{ J \tau \cdot \Diamond}   \partial _{p_k} \Phi _p,
  e^{ J \tau \cdot \Diamond}  \partial _{p_j}\Phi _p )
   \\& =  \Omega (
 R  ,  \partial _{p_k}  \partial _{p_j}\Phi _p )  . \end{aligned}
\end{equation*}
By (A3) we have

\begin{equation*} \label{ }\begin{aligned} &  \partial _{\tau _k}\mathcal{G}_j   = \Omega (
e^{ J \tau \cdot \Diamond} R  ,e^{ J \tau \cdot \Diamond} J ^2  \Diamond _k \Diamond _j\Phi _p ) - \Omega (  J   \Diamond _k   e^{ J \tau \cdot \Diamond}   \Phi _p
  e^{ J \tau \cdot \Diamond} J   \Diamond _j\Phi _p )
   \\& =   - \langle  R  ,    J\Diamond _k \Diamond _j\Phi _p    \rangle  -  \langle J \Diamond _k      \Phi _p ,
     \Diamond _j\Phi _p   \rangle   =    - \langle  R  ,    J\Diamond _k \Diamond _j\Phi _p    \rangle   , \end{aligned}
\end{equation*}
We have
\begin{equation*} \label{ }\begin{aligned} &  \partial _{p _k}\mathcal{G}_j  = \Omega (
e^{ J \tau \cdot \Diamond} R  ,e^{ J \tau \cdot \Diamond} J \Diamond _j  \partial _{p_k}  \Phi _p ) - \Omega (         e^{ J \tau \cdot \Diamond}   \partial _{p_k} \Phi _p
  e^{ J \tau \cdot \Diamond}  J \Diamond _j\Phi _p )
   \\& =  -\langle  R  ,  \Diamond _j  \partial _{p_k}  \Phi _p \rangle  +\langle  \partial _{p_k}  \Phi _p   ,  \Diamond _j     \Phi _p \rangle  = -\langle  R  ,  \Diamond _j  \partial _{p_k}  \Phi _p \rangle  + \delta _{jk}
	  . \end{aligned}
\end{equation*}
At $U=\Psi _0$, $\tau  =\tau _0$ and $p=p_0$ we have $\mathcal{F}_j  =\mathcal{G}_j  =0$. Since in this case $R=0$
we get the desired result by the  Implicit Function Theorem. \qed

\subsection{Spectral coordinates}
\label{subsec:coordinates}
Lemmas \ref{lem:modulation}--\ref{lem:begspectdec} lead to a  natural
  decomposition of \eqref{eq:NLSvectorial}.  To write it we need
further notation.

\noindent We are ready for the  natural
coordinates decomposition.
Let $\Pi (U_0)=p_0$.
We  consider
for $R\in N_g^\perp (\mathcal{H} _{ p_0 }^{\ast})$ the map
\begin{equation}\label{eq:coordinate}\begin{aligned} &
(\tau , p, R) \to U= e^{ J \tau \cdot \Diamond} (  \Phi _{p } +P(p)R).
  \end{aligned}  \end{equation}
We have the following formulas,
\begin{equation}\label{eq:vectorfields} \begin{aligned} &
\frac \partial {\partial  {\tau _j}}   = J    \Diamond  _jU
\, ,\quad  \frac \partial {\partial  {p _j}}= e^{ J \tau \cdot \Diamond}  (  \partial _ {p _j}\Phi _{p } + \partial _ {p _j}P(p)R) ,\end{aligned}
\end{equation}
with $\frac \partial {\partial  {p _j}} \in C^{\infty} (\U \cap \Sigma _{k }, \Sigma _{k' })$
for any pair  $(k ,k' )\in \N ^2$, with $\U  \subset \Sigma _{-n}$ the neighborhood
of $ e^{J\tau  _0\cdot \Diamond }\Phi _{p _0} $ in
Lemma \ref{lem:modulation}.  Similarly, $\frac \partial {\partial  {\tau _j}} \in C^{0} (\U \cap \Sigma _{k }, \Sigma _{k-d_j })$.  We have what follows.

\begin{lemma} \label{lem:gradient R}  Consider the $n\ge 0$ and $\U$ 
in Lemma  \ref{lem:modulation} and fix  an integer $k\ge -n$.
Then the map
$U\to R(U)=R$  is $ C^{0} (\U \cap \Sigma _{k }, \Sigma _k)$.  For $k\ge -n +\textbf{d}$
  we have
$R\in C^1(\U \cap \Sigma _{k }, \Sigma _{k -\textbf{d}})$. 
  For $\U$ sufficiently small in $\Sigma _{-n }$ 
the
Frech\'et derivative
$R'(U) $   of $R(U)$ is defined by the following formula, summing on the repeated index $j$,
\begin{equation*}  \begin{aligned}
R'( U)&=   (P (p )P ( p_0))^{-1}  P (p )\big [ e^{ - J \tau \cdot \Diamond}    \uno-
 J \Diamond _j P (p )  R \, d\tau _j-      \partial _{p_j} P (p )  R \, dp _j  \big  ]  ,
\end{aligned}
\end{equation*} 
where $(P (p )P ( p_0))^{-1} :N_g^\perp (\mathcal{H} _{p }^{\ast})\cap \Sigma _{k -\textbf{d}}\to N_g^\perp (\mathcal{H}_{p_0 }^{\ast})\cap \Sigma _{k -\textbf{d}} $   is well defined by  Lemma \ref{lem:projections}.

\end{lemma}
\proof The continuity of $ R(U)$ follows   from $R=e^{-J\tau \cdot
\Diamond}U-\Phi _p$ and
 \begin{equation*}  \begin{aligned} &
 R-R'= e^{-J \tau \cdot
\Diamond}U    -e^{-J\tau '\cdot
\Diamond}U '+\Phi  _{p'} -\Phi _p =\\&\Phi  _{p'} -\Phi _p  +   (e^{-J \tau \cdot
\Diamond}  -e^{-J\tau '\cdot
\Diamond}  )U      +e^{-J\tau '\cdot
\Diamond}(U-U ' ) .
\end{aligned}
\end{equation*}
 Then use $p\to \Phi _p\in C^\infty ( {\mathcal O} ,  {\mathcal S} )$,
 the fact that $e^{ J \tau \cdot
\Diamond}$ is strongly continuous in $\Sigma _k$ and locally uniformly bounded therein. The fact that  $ R(U)$  has  Frech\'et derivative
follows by the chain rule. To get the formula for $R'(U)$ notice that the equalities $R'\frac \partial {\partial  {p _j}}=R'\frac \partial {\partial  {\tau _j}}=0$ and $R'e^{   J \tau \cdot \Diamond}   P (p )P ( p_0)= \uno _{|N_g^\perp (\mathcal{H}_{p_0}^{\ast})} $ characterize $R'$.
We claim we have
\begin{equation}\label{eq:frechr}
 R'=\mathbf{a} _j d\tau _j + \mathbf{b} _j dp _j +(P (p )P ( p_0))^{-1}  P (p )  e^{ - J \tau \cdot \Diamond}
\end{equation}
for some $\mathbf{a} _j $ and $\mathbf{b} _j$. First of all, by the independence of coordinates $(\tau ,p)$ from
$R\in N_g^\perp (\mathcal{H}_{p_0}^{\ast})$,
\begin{equation*}
 d\tau_j \circ e^{   J \tau \cdot \Diamond}   P (p )P ( p_0)=dp_j \circ e^{   J \tau \cdot \Diamond}   P (p )P ( p_0)=0.
\end{equation*}
Indeed for $g \in  N_g^\perp (\mathcal{H}_{p_0}^{\ast}) $ we have for instance
\begin{equation*}\begin{aligned}
0= \frac {d}{dt} \tau _j (u (\tau , p , R+tg)) _{|t=0} &= \frac {d}{dt} \tau _j ( e^{   J \tau \cdot \Diamond}   ( \Phi
_p + P(p) P(p_0) (R+tg))) _{|t=0} \\& =
d\tau_j \circ e^{   J \tau \cdot \Diamond}   P (p )P ( p_0)g . \end{aligned}
\end{equation*}
Secondarily, by the definition of  $(P (p )P ( p_0))^{-1}$,
\begin{equation*}
 (P (p )P ( p_0))^{-1}  P (p )  e^{ - J \tau \cdot \Diamond}  \circ  e^{   J \tau \cdot \Diamond}   P (p )P ( p_0) =\uno _{N_g^\perp (\mathcal{H}_{p_0}^{\ast})}.
\end{equation*}
Hence we get the claimed equality  \eqref{eq:frechr}.

\noindent To get    $\mathbf{a} _j $ and $\mathbf{b} _j$ notice that
by $R' \frac \partial {\partial  {\tau _j}}=0$ and $P(p)J\Diamond _j  \Phi _p
=0$
\begin{equation*}  \begin{aligned} &
\mathbf{a} _j= -(P (p )P ( p_0))^{-1}  P (p )  e^{ - J \tau \cdot \Diamond}
\frac \partial {\partial  {\tau _j}}=\\& -(P (p )P ( p_0))^{-1}  P (p )  e^{ - J \tau \cdot \Diamond} e^{   J \tau \cdot \Diamond}  J\Diamond _j  (\Phi _p+P(p)R)=\\&   -(P (p )P ( p_0))^{-1}  P (p )     J\Diamond _j   P(p)R  .
  \end{aligned}
\end{equation*}
Similarly by $R'\frac \partial {\partial  {p _j}} =0$ and $P(p)\partial _{p_j}\Phi _p
=0$

\begin{equation*}  \begin{aligned} &
\mathbf{b} _j= -(P (p )P ( p_0))^{-1}  P (p )  e^{ - J \tau \cdot \Diamond}
\frac \partial {\partial  {p _j}}=\\& -(P (p )P ( p_0))^{-1}  P (p )     (
\partial _{p_j}\Phi _p+\partial _{p_j}P(p)R)=\\&   -(P (p )P ( p_0))^{-1}  P (p )    \partial _{p_j}P(p)R   .
  \end{aligned}
\end{equation*}

\qed

\bigskip
A crucial point in the  stability proofs in  \cite{bambusicuccagna,boussaidcuccagna,Cu1,Cu2}, first realized
and used in \cite{Cu3},
is the importance not to loose track of the Hamiltonian nature of
\eqref{eq:NLSvectorial}, in whichever coordinates the system is written.
Thus we have what follows.

\begin{lemma}\label{lem:HamFor} In the coordinate system \eqref{eq:coordinate}, system \eqref{eq:NLSvectorial} can be written  as
\begin{equation} \label{eq:SystPoiss} \begin{aligned} &
  \dot p  = \{ p , E \} \, ,  \,    \dot \tau   = \{ \tau   , E \} \, ,  \,    \dot R   = \{ R   , E \}
	    . \end{aligned}
\end{equation}
\end{lemma}
\proof
The statement is not standard only for   $\dot R   = \{ R   , E \}$.  Notice that \eqref{eq:NLSvectorial} can be written as
\begin{equation}\label{eq:sys1}\begin{aligned} &
  \dot  U=  J \dot \tau \cdot  \Diamond U+
	e^{ J \tau \cdot \Diamond} \dot p\cdot  \nabla _p(  \Phi _{p } +P(p)R) +e^{ J \tau \cdot \Diamond}   P(p)
	\dot R \\& =\sum _j  \dot \tau _j   \frac{\partial }{\partial \tau _j}+
  \dot p_j   \frac{\partial }{\partial p _j} +e^{ J \tau \cdot \Diamond}   P(p)
	\dot R =   J \nabla E  (  U) .\end{aligned}
\end{equation}
When we apply the  derivative $R'(U)$ to \eqref{eq:sys1},  all the terms in the lhs of the last
line cancel except for
\begin{equation}
\begin{aligned} &  R'(U)
  e^{ J \tau \cdot \Diamond}   P(p)
	\dot R  = R'(U) J
  \nabla E  (  U) =R'(U)X_E(U)=\{ R,E\} ,
\end{aligned}\nonumber
\end{equation} from
the definition of hamiltonian field  and of Poisson bracket.
 Finally we use
\begin{equation}
\begin{aligned} &  R'(U)
  e^{ J \tau \cdot \Diamond}   P(p)
	\dot R  = \frac{d}{ds}_{|_{s=0}}R( U (\tau , p , R+s\dot R) )   =\frac{d}{ds}_{|_{s=0}} (R+s\dot R) = \dot R   .
\end{aligned}\nonumber
\end{equation}

\qed

\subsection{Reduction of order of system \eqref{eq:SystPoiss}}
\label{subsec:reduction}

The following Poisson bracket identities are useful.
\begin{lemma}
  \label{lem:Involutions} Consider the  functions $\Pi _j$.  Then
$ X  _{\Pi _j }=
\frac{\partial}{\partial \tau _j} .$
In particular
  \begin{equation} \label{eq:Ham.VecField1}  \{ \Pi _j,\tau _k   \}  =-\delta _{jk} \ , \quad   \{ \Pi _j,p _k   \} \equiv 0  \ , \quad   \{ R, \Pi _j   \}  = 0.
\end{equation}

\end{lemma}
\proof   \eqref{eq:Ham.VecField1}   follows from the first claim, which is a consequence of      \eqref{eq:vectorfields}:

  \begin{equation}
\begin{aligned} &  X  _{\Pi _j } (U) = J\nabla \Pi _j (U)  = J \Diamond _jU=    \frac{\partial}{\partial \tau _j}  . \end{aligned}\nonumber
\end{equation}

\qed

\noindent We introduce now a new Hamiltonian:
 \begin{equation}     \label{eq:K}  \begin{aligned} &
 {K}(U):=E(U)- E\left (   \Phi _{p _0}\right ) -  \lambda _j(p)    \left ( \Pi _j (U)- \Pi _j (U_0)\right ) .
\end{aligned}
\end{equation}
Notice that $K(e^{J\tau \cdot \Diamond}U)\equiv K(U)$. Equivalently,  $\partial _{\tau _j} {K}  \equiv 0$.
We know that  for solutions of \eqref{eq:NLSvectorial} we have $\Pi _j (U (t))= \Pi _j (U_0) $ and

\begin{equation*}       \begin{aligned} &   \{ p _j, K   \} =  \{ p _j, E   \} \, , \,  \{ R , K   \} =  \{ R, E   \} \, , \,
  \{ \tau _j  , K   \} =  \{ \tau _j , E   \}+   \lambda _j(p) .
\end{aligned}
\end{equation*}
By   $\partial _{\tau _j} {K}  \equiv 0$,   the evolution    of the variables $p ,R$ is unchanged if we consider the following new
Hamiltonian system:
\begin{equation} \label{eq:SystK} \begin{aligned} &
\dot p_j  =   \{ p _j, K   \}  \, , \quad   \dot \tau _j   = \{  \tau _j , K  \}  \, , \quad    \dot R= \{ R,  {K}
\}     .  \end{aligned}
\end{equation}
It is elementary that      the \textit{momenta} $\Pi _j(U)$ are invariants of motion of \eqref{eq:SystK}.

Before exploiting the invariance of $ \Pi _j(U)$  to reduce the order of
the system, we introduce appropriate notation.
First of all we  set \begin{equation}\label{eq:PhaseSpace}\begin{aligned} &
  {\Ph}^{r }:=\R ^{ n_0}\times   (\Sigma _{r } \cap N_g^\perp ({\mathcal H}_{p_0}))=\{(\tau   , R)\}
\, , \\&
\widetilde {\Ph}^{r }:= \R ^{ n_0}\times {\Ph}^{r }=\{(\Pi ,\tau   , R)\}  .  \end{aligned}
\end{equation}
 We set  ${\Ph} ={\Ph}^{0 }$ and $\widetilde {\Ph}=
\widetilde {\Ph}^{0 }$.

\begin{definition}\label{def:scalSymb}
We will say that   $F(t,\varrho , R)\in C^{M}(I\times \mathcal{A},\R)$ with
$I$ a neighborhood of 0 in $\R$ and
 $\mathcal{A}$   a neighborhood of 0 in
$  {\Ph}^{-K }$
 is $\mathcal{R}^{i,j}_{K, M}$   and we will write $F=\mathcal{R}^{i,j}_{ K,M}$, or more specifically $F=\mathcal{R}^{i, j}_{ K,M} (t,\varrho , R)$,
 if    there exists    a $C>0$   and a smaller neighborhood  $\mathcal{A}'$ of 0   s.t.
 \begin{equation}\label{eq:scalSymb}
  |F(t,\varrho , R)|\le C \|  R\| _{\Sigma   _{-K}}^j (\|  R\| _{\Sigma   _{-K}}+|\varrho |)^{i} \text{  in $I\times \mathcal{A}'$}.
\end{equation}
We say  $F=\mathcal{R}^{i, j} _{K, \infty}$  if $F=\mathcal{R}^{i,j}_{K, m}$ for all $m\ge M$.
We say    $F=\mathcal{R}^{i, j}_{\infty, M} $       if   for all   $k\ge K$    the above   $F$ is the restriction  of an 
$F(t,\varrho , R)\in C^{M}(I\times \mathcal{A}_{k },\R)$ with  $\mathcal{A}_k$   a neighborhood of 0 in
$  {\Ph}^{-k }$ and 
 which is
$F=\mathcal{R}^{i,j}_{k, M}$. Finally we say 
$F=\mathcal{R}^{i, j}  $   if $F=\mathcal{R}^{i, j} _{K, \infty}$ and  $F=\mathcal{R}^{i, j} _{ \infty , M}$.

\end{definition}

\begin{definition}\label{def:opSymb}  We will say that an   $T(t,\varrho , R)\in C^{M}(I\times \mathcal{A},\Sigma   _{K}  (\R^3, \R ^{2N}))$,  with $I\times \mathcal{A}$    like above,
 is $ \mathbf{{S}}^{i,j}_{K,M} $   and we will write $T= \mathbf{{S}}^{i,j}_{K,M}$  or more specifically $T= \mathbf{{S}}^{i,j}_{K,M} (t,\varrho , R)$,
 if     there exists  a $C>0$   and a smaller neighborhood  $\mathcal{A}'$ of 0   s.t.
 \begin{equation}\label{eq:opSymb}
  \|T(t,\varrho , R)\| _{\Sigma   _{K}}\le C \|  R\| _{\Sigma   _{-K}}^j (\|  R\| _{\Sigma   _{-K}}+|\varrho |)^{i}  \text{  in $I\times \mathcal{A}'$}.
\end{equation}
We  use notation  $T=\mathbf{{S}}^{i, j} $, $T=\mathbf{{S}}^{i,j}_{K,\infty }$    or  $T=\mathbf{{S}}^{i,j}_{\infty,M}$  as above.

\end{definition}
These notions will be often used also for functions $F=\mathcal{R}^{i, j}_{ K,M} ( \varrho , R)$ and  $T=\mathbf{{S}}^{i, j}_{ K,M} ( \varrho , R)$
  independent of $t$.

\begin{remark}
\label{rem:formal1}
We will see later that the coefficients  of the vector fields
whose flows are used to change coordinates are symbols as of
Definitions \ref{def:scalSymb} and \ref{def:opSymb}.     
 The definitions of  the   symbols    
$ \mathcal{R}^{i, j} $ and  $ \mathbf{{S}}^{i, j}   $ 
in   Def. 3.9 and 3.10  \cite{bambusi}  are very restrictive, since they 
require for the symbols  to be   defined in the whole 
$I\times \mathcal{S}'$.   The proofs in  \cite{bambusi}  at most prove
that the   coefficients  of the vector fields  in  fact  are symbols 
of the form  $ \mathcal{R}^{i, j} _{K,M}  $ and  $ \mathbf{{S}}^{i, j}  _{K,M}    $ in our sense.
As an example we refer to  Lemmas 3.26   and 5.5  in \cite{bambusi}.  In  Lemma  3.26  \cite{bambusi}
   the fact that  the $b_i$ and the  $\langle W^l;Y\rangle $ are symbols of the form  $ \mathcal{R}^{j, k}   $  
	for some $(j,k)$  
			in the sense of  Def.  3.10   in \cite{bambusi},   requires preliminarily to show  
			at least that they are functions of $(\varrho ,R)$  for $(\varrho ,R)$ in some neighborhood  $\U$ of $(0,0)$ in
			$\R ^{n_0}\times \mathcal{S}'$. 
			This  is not addressed in  \cite{bambusi} and      is far from trivial,   since the coefficients of the   linear system 
			right  above
			formula (3.60) are    are    unbounded in any such $\U$.   The justification that the coefficients
			$\Phi _{\mu \nu}(M)$ of $\chi $ in Sect. 5  in  \cite{bambusi} are in $\mathcal{S}$ is  similarly  inconclusive.
				The key step should be that the homological equation in 
				Lemma 5.5  can be solved  for all parameters $k$  uniformly in the variable   $M\in \R ^n$, at least
				for $|M|< a $  for a fixed $a$.   But the  homological equations involve the perturbation of an operator
				and in \cite{bambusi} the perturbation is not fully analyzed. For example there is no discussion of
				the norm  $\| V_M-V_0\| _{\mathcal{W}^k\to \mathcal{W}^k}$  as $k$ grows
				and $|M|< a $.  This norm should be expected to grow  
				and become large, possibly breaking  down  the proof of  $\Phi _{\mu \nu}(M)\in  \mathcal{S}$.
				In fact it is plausible that  $\Phi _{\mu \nu}(M)\in  \mathcal{S}$ only for $M=0$.

		From the above remarks we can see that no  coordinate  change   in the Birkhoff or in the Darboux steps   in \cite{bambusi} is  shown to be an \textit{almost smooth} transformation  in the sense of Definition 3.15  in \cite{bambusi}.  So  for instance the proof  of the Birkhoff normal forms,
that  is   Theor. 5.2   \cite{bambusi},  is  inconclusive.   The proof of  the Darboux step, that is  Theor. 3.21  \cite{bambusi}, 
is even sketchier and is similarly inconclusive.
\end{remark}

\bigskip
We   proceed  now to  a  reduction of order in \eqref{eq:SystK}. Write
	\begin{equation}\label{eq:variables} \begin{aligned} &
 \Pi _j(U) =\Pi _j(e^{ J \tau \cdot \Diamond} (  \Phi _{p } +P(p)R))=\Pi _j(  \Phi _{p } +P(p)R) \\& = \frac{1}{2}
\langle \Diamond _j(  \Phi _{p } +P(p)R) ,   \Phi _{p } +P(p)R = p_j+ \Pi _j(P(p)R)\\& = p_j+ \Pi _j( R)  + \Pi _j( (P(p)-P(p_0)) R) +\langle R, \Diamond _j(P(p)-P(p_0)) R\rangle    .
\end{aligned}
\end{equation}
We well move from variables $(\tau , p , R)$   to variables $(\tau , \Pi , R)$. Setting  $ \varrho _j=\Pi _j( R)$,
we have

	\begin{equation}\label{eq:variables1} \begin{aligned} &
 p_j=\Pi _j-\varrho _j+ \widetilde{\Psi} _j (p-p_0, R)
\end{aligned}
\end{equation}
 with
$ \widetilde{\Psi} _j  =\mathcal{R}^{0, 2}(p-p_0, R)$.
The implicit function theorem yields:
\begin{lemma}
\label{lem:var} There   are functions $p_j =p_j ( \Pi ,\Pi  ( R) , R)$
  defined implicitly  by  \eqref{eq:variables},or \eqref{eq:variables1},
  such that $ p_j=\Pi _j-\varrho _j+{\Psi} _j(\Pi , \varrho ,R)
  $ with
  ${\Psi}  (p_0 , \varrho ,R)=
   \mathcal{R}^{0, 2} (\varrho , R)
   $.
\end{lemma}
We consider now $(\tau , \Pi , R)$ as a new coordinate system.
By   $
    \frac{\partial}{\partial \tau _k } \Pi _j(U)\equiv 0 $ it follows that the vectorfields $ \frac{\partial}{\partial \tau _k }$ are the same for the two systems of coordinates.
 In the new variables, system \eqref{eq:SystK} reduces to the pair
 of systems

 \begin{align} &  \label{eq:SystK1}   \dot \tau _j   = \{ \tau _j   ,K \} \, , \quad    {\dot
 \Pi _j }=  0   \, ,  \\&  \label{eq:SystK21}
	   \dot R= \{ R, K
\}     .
 \end{align}
 System  \eqref{eq:SystK21}   is closed  because of $\partial _{\tau _j}K=0$.

\section{Darboux Theorem}
\label{sect:symplectic}

In this section we present one of the two main results of this paper.
We seek to reproduce Moser's proof of the Darboux theorem. Specifically
we look for a vector field ${\mathcal X}^t$ that will produce a flow
as in \eqref{eq:fdarboux} below.  The proof of the existence and properties of
${\mathcal X}^t$   is similar to \cite{Cu2}, but     influenced
by the  choice of coordinates in \cite{bambusi}. We also add material
to justify, once ${\mathcal X}^t$ has been found,  the formal formula    \eqref{eq:fdarboux}. Notice that for \cite{boussaidcuccagna,Cu2} formula
\eqref{eq:fdarboux} does not require justification because ${\mathcal X}^t$ is a smooth
vectorfield on a given manifold.  But the situation in  \cite{Cu1,bambusi} is different since now ${\mathcal X}^t$ is not
a standard vectorfield on a manifold  and  $\Omega $ is not a regular
differential form on the same manifold, so
Lie derivative,
pullbacks, push forwards  and the related differentiation formulas, require justification.

Notice that, to be useful in the asymptotic stability theory,  the change of variables has to be such that the new Hamiltonian equations is semilinear.
This is why even in  \cite{boussaidcuccagna,Cu2},   where
we could apply the standard Darboux theorem for strong symplectic forms
on Banach manifolds, see \cite{amr} Ch. 9,  it is important to select ${\mathcal X}^t$ with an {\it ad hoc} process.

\subsection{Search of a vectorfield}
\label{subsec:vector}
Recall that  $\Omega =\langle  J ^{-1}\ , \ \rangle
$ and  consider
\begin{equation} \label{eq:Omega0} \Omega _0 :=  d  \tau _j \wedge
d\Pi  _j  + \langle  J ^{-1}R' , R' \rangle
 .
\end{equation}

\begin{lemma}
  \label{lem:1forms}
  At the points  $e^{ J \tau \cdot \Diamond}   \Phi _{p_0 }$ for all $\tau  \in \R  ^{n_0 }$  we have $\Omega _0 =\Omega   .$

  \noindent Consider the following  forms:
\begin{equation} \label{eq:1forms1}\begin{aligned} &
 \mathrm{ B} _0 :=\tau _j
d\Pi  _j  +\frac{1}{2} \langle  J ^{-1}R  , R' \rangle ;  \quad \mathrm{B}:=\mathrm{B_0}+\alpha  \text{ for }
\end{aligned}
\end{equation}
   \begin{equation} \label{eq:1forms2}\begin{aligned} &
\alpha  :=-\beta _j(p,R)  d\Pi _j  +   \left \langle  \Gamma (p)  R +
\beta  _j(p,R) P^*(p) \Diamond _j P(p)R,R'
\right \rangle \ ,
     \\&
 \Gamma (p):=\frac{1}{2}    J ^{-1}  \left ( P(p) -P(p_0)  \right ) \ ,
\\& \beta  _j(p,R) :=\frac{1}{2}\
\frac{  \langle   P^*(p) J ^{-1} R, \partial  _{p_{ j}}P(p)R   \rangle  }{ 1+
\langle  \Diamond _j P(p)R, \partial  _{p_{ j}}P(p)R   \rangle  } \
   .
\end{aligned}
\end{equation}

   Then
  $ d\mathrm{B_0}=\Omega _0$ and $
  d\mathrm{B } =\Omega   .$

\end{lemma}
\proof     $ d\mathrm{B_0}=\Omega _0$ follows from the definition of exterior differential.
Set $\widetilde{B} :=\frac{1}{2}\langle  J ^{-1} U , \ \rangle $. Notice that $d \widetilde{B}=\Omega $.
By \eqref{eq:coordinate} we get:

\begin{equation} \label{eq:beta1} \begin{aligned} &
\widetilde{B} (X)= \frac{1}{2}\langle  J ^{-1}e^{  J \tau \cdot \Diamond}     \Phi _{ {p} }   ,   X \rangle
+ \frac{1}{2}\langle  J ^{-1}     P(p)R , e^{ -J \tau \cdot \Diamond} X \rangle .
\end{aligned}
\end{equation}
Set    $\psi (U):=\frac{1}{2}\langle J ^{-1}  e^{  J \tau \cdot \Diamond}     \Phi _{ {p} } ,U
   \rangle .$  Then  we claim
\begin{equation*} \label{eq:dPsi0} \begin{aligned} & d\psi =\frac{1}{2}\langle J ^{-1}   e^{  J \tau \cdot \Diamond}    \Phi _{ {p} } ,\
   \rangle  + p_jd\tau _j,\end{aligned}
\end{equation*}
where in this proof we will sum on repeated indexes.
The last formula implies
\begin{equation} \label{eq:beta11} \begin{aligned} &
\widetilde{B}  =  d\psi - p_jd\tau _j
+ \frac{1}{2}\langle  J ^{-1}     P(p)R , e^{ -J \tau \cdot \Diamond} \ \rangle .
\end{aligned}
\end{equation}
The desired formula on $d\psi$  follows by

\begin{equation*}  \begin{aligned} & d\psi =\frac{1}{2}\langle J ^{-1}   e^{  J \tau \cdot \Diamond}    \Phi _{ {p} } ,\
   \rangle  + \frac{1}{2}  \langle     e^{  J \tau \cdot \Diamond}   \Diamond _j \Phi _{ {p} } ,U
   \rangle  d\tau _j \\&  + \frac{1}{2}  \langle     e^{  J \tau \cdot \Diamond}  J ^{-1}  \partial  _{ {p}_j }\Phi _{ {p} } ,U
   \rangle  dp_j   = \frac{1}{2}\langle J ^{-1}   e^{  J \tau \cdot \Diamond}    \Phi _{ {p} } ,\
   \rangle  +\\&  \frac{1}{2}  \langle         \Diamond _j \Phi _{ {p} } ,\Phi _{ {p} }+P(p)R
   \rangle  d\tau _j    + \frac{1}{2}  \langle        J ^{-1}  \partial  _{ {p}_j }\Phi _{ {p} } ,\Phi _{ {p} } +P(p)R
   \rangle  dp_j \stackrel{\text{by \eqref{eq:begspectdec3}}}{=}
	   \\&  \frac{1}{2}\langle J ^{-1}   e^{  J \tau \cdot \Diamond}    \Phi _{ {p} } ,\
   \rangle+ \underbrace{\frac{1}{2}  \langle         \Diamond _j \Phi _{ {p} } ,\Phi _{ {p} }
   \rangle}_{p_j}  d\tau _j    + \frac{1}{2} \underbrace{ \langle        J ^{-1}  \partial  _{ {p}_j }\Phi _{ {p} } ,\Phi _{ {p} }
   \rangle}_{\qquad \quad 0 \text{ by \eqref{eq:gen kernel20}}}  dp_j .
\end{aligned}
\end{equation*}
By Lemma \ref{lem:gradient R} and using $ P(p)^*J ^{-1} =J ^{-1}P(p)$ we have
\begin{equation*}   \begin{aligned} & \frac{1}{2}\langle  J ^{-1}     P(p)R , e^{ -J \tau \cdot \Diamond} \ \rangle =\frac{1}{2}\langle  J ^{-1}      R ,P(p) R'  \ \rangle \\& +\frac{1}{2}\langle  J ^{-1}      R ,P(p)J \Diamond _jP( {p})R   \rangle d\tau _j+\frac{1}{2}\langle  J ^{-1}      R ,P(p)  \partial _{p_j}P( {p})R   \rangle dp _j \\& =\frac{1}{2}\langle  J ^{-1}      R , R'  \ \rangle   + \frac{1}{2}\langle  J ^{-1}      R ,(P(p) -P( {p}_0) ) R'  \ \rangle \\& -\Pi _j(P( {p})R)  d\tau _j+\frac{1}{2}\langle  J ^{-1}      R ,P(p)   \partial _{p_j}P( {p})R   \rangle dp _j .
\end{aligned}
\end{equation*}
So by \eqref{eq:beta11}  and using $P(p)J =JP^*(p)$  we get
\begin{equation*}  \begin{aligned} &
\widetilde{B}  - d\psi = - (\overbrace{ p_j+ \Pi _j(P( {p})R)}^{\Pi_j} ) d\tau _j
+ \frac{1}{2}\langle  J ^{-1}      R , R'  \ \rangle   \\&+ \frac{1}{2}\langle  J ^{-1}      R ,(P(p) -P(p_0) ) R'  \ \rangle   -\frac{1}{2}\langle     P^*( {p}) J ^{-1}      R ,   \partial _{p_j}P( {p})R   \rangle dp _j  .
\end{aligned}
\end{equation*}
Then  $d\alpha =\Omega -\Omega _0$ for
\begin{equation*} \label{eq:beta12} \begin{aligned} & \alpha :=
\widetilde{B}  - d\psi -B_0 +d(\Pi_j \tau _j) =    \\&  \frac{1}{2}\langle  J ^{-1}      R ,(P(p) -P( {p}_0) ) R'  \ \rangle   -\frac{1}{2}\langle       P^*( {p}) J ^{-1}  R ,   \partial _{p_j}P( {p})R   \rangle dp _j  .
\end{aligned}
\end{equation*}
By $p_j=\Pi _j - \Pi _j(P(p)R)$  we get

 \begin{equation*} \label{eq:beta13} \begin{aligned} & dp_j=d\Pi _j  -\langle \Diamond _j P(p) R, P(p)R'\rangle - \langle \Diamond _j P(p) R, \partial _{p_j}P(p)R \rangle  dp_j.
\end{aligned}
\end{equation*}
	 Then inserting  the next formula  in  the formula for $\alpha$, we obtain \eqref{eq:1forms2}:

 \begin{equation} \label{eq:beta14} \begin{aligned} & dp_j=    \frac{ d\Pi _j -\langle \Diamond _j P(p) R, P(p)R'\rangle}{1+\langle \Diamond _j P(p) R, \partial _{p_j}P(p)R \rangle}   .
\end{aligned}
\end{equation}
\qed

In the Lemmas \ref{lem:dalpha1}--\ref{lem:vectorfield} we will initially
consider the regularity of the functions in terms of the coordinates
$(\tau , p, R)$.

\begin{lemma}
  \label{lem:dalpha1}   We have $\beta _j  \in C^\infty ( {\mathcal O} \times \Sigma _{-n }, \R )$
	for any   $n$.  For any pair    $(n,n' )$ we have $\Gamma \in C^\infty ( {\mathcal O} , B( \Sigma _{-n' }, \Sigma _{ n })  )$.
	  Summing on repeated indexes, we have

  \begin{equation} \label{eq:dalpha1}
\begin{aligned} & d\alpha  =-\partial _{p_k}\beta _j dp_k \wedge  d\Pi _j - \langle \nabla _R\beta _j,R'
\rangle \wedge  d\Pi _j \\& + dp_k \wedge  \langle  \partial _{p_k}[\Gamma (p)  R +
\beta  _j(p,R) P^*(p) \Diamond _j P(p)R],R' \rangle
\\& +       \langle \nabla _R \beta  _j , R'\rangle  \wedge \langle   P^*(p) \Diamond _j P(p)R,R'\rangle +2\langle \Gamma R',R'\rangle
   . \end{aligned}
\end{equation}

\end{lemma}
\proof    Follows from a simple computation.  In particular, for  a $\mathbf{L}\in B(\Sigma _1, L^2)$ fixed, we use the formula
\begin{equation*}  \begin{aligned} &
 d\langle \mathbf{L}R, R'\rangle (X,Y):=X\langle \mathbf{L}R, R'Y\rangle-Y\langle \mathbf{L}R, R'X\rangle
  - \langle \mathbf{L}R, R'[X,Y]\rangle \\&
 =
 \langle \mathbf{L}R'X, R'Y\rangle -\langle \mathbf{L}R'Y, R'X\rangle . \end{aligned}
\end{equation*}

\qed

\begin{lemma}
  \label{lem:dalpha2}     Summing on repeated indexes, we have

  \begin{equation*}
\begin{aligned}   d\alpha  &= \widehat{\delta } _   k    \partial _{p_k}{\beta } _ {j } d\Pi _j \wedge  d\Pi _k + \langle \widehat{\Gamma}_j+
(  \widehat{\delta } _   k    \partial _{p_k}{\beta } _ {j } -\widehat{\delta } _   j    \partial _{p_j}{\beta } _ {k } )    \Diamond _k P(p)R,R'
\rangle \wedge  d\Pi _j \\& +   2 \langle  \Gamma  (p)  R ', R'\rangle
+
\langle \widetilde{\beta}  _j   , R'\rangle \wedge \langle  P^*(p) \Diamond _j P(p)R,R' \rangle
 \  ,\end{aligned}
\end{equation*}
	where we have  (this time not summing on repeated indexes)

  \begin{equation*}
\begin{aligned}
\widehat{\delta} _k &:= \frac 1{1+\langle \Diamond _k P(p) R, \partial _{p_k}P(p)R \rangle}
\  ,   \\  \widehat{\Gamma}_j&:=-\nabla _R \beta _j-\widehat{\delta} _j  [  \partial _{p_j}\Gamma    R +\sum _{i=1}^{n_0}\beta  _i\partial _{p_j}
\left (   P^*(p) \Diamond _i P(p)  \right )R  ]\\&
+  \sum _{k=1}^{n_0}(  \widehat{\delta } _   k    \partial _{p_k}{\beta } _ {j } -\widehat{\delta } _   j    \partial _{p_j}{\beta } _ {k } )   ( P^*(p)-1)\Diamond _k P(p)   R
\\
\widetilde{\beta}  _j&:= \nabla _R \beta  _j+   \widehat{\delta} _j  \partial _{p_j}  (  \Gamma      +\sum _{k=1}^{n_0}
  \beta  _k  P^*(p) \Diamond _k P(p)    )R
   . \end{aligned}
\end{equation*}

\end{lemma}
\proof  Follows by an elementary computation substituting \eqref{eq:beta14} in \eqref{eq:dalpha1}
\qed

\begin{lemma}
  \label{lem:vectorfield0}   For any fixed large $n$ and for $\varepsilon _0>0$,
	   consider the set $\U _{\mathbf{d}} \subset \widetilde{\Ph}^{\textbf{d}}=\{  (p,R) \}$ defined by   $  \|R\| _{\Sigma  _{-n}}\le   \varepsilon _0$  and
   $| p -p_0|\le   \varepsilon _0$.
  Then for $ \varepsilon _0 $  small enough
     there exists a unique vectorfield  $ \mathcal{X}^t  :\U  _{\mathbf{d}}  \to  \widetilde{\Ph}  $       which solves
$
 i_{\mathcal{X}^t} \Omega _t=-  \alpha
   $, where $\Omega _t:=\Omega _0+t(\Omega -\Omega _0)$.

\end{lemma}

\proof
First of all we consider $Y$ such that $
 i_{Y} \Omega _0=-  \alpha
   $, that is to say

\begin{equation*}    \begin{aligned} &
   (Y ) _{\tau _j} d\Pi _j -  (Y ) _{\Pi  _j} d\tau  _j  +\langle J^{-1} (Y ) _{R} ,R' \rangle \\ &= \beta _j(p,R)  d\Pi _j  -   \left \langle  \Gamma (p)  R +
\beta  _j(p,R) P^*(p) \Diamond _j P(p)R,R'
\right \rangle   .
   \end{aligned}
\end{equation*}
This yields

\begin{equation}   \label{eq:Yvec}   \begin{aligned} &
   (Y ) _{\tau _j} =\beta _j(p,R)=\mathcal{R}^{0,2}(p,R)   \ , \quad  (Y ) _{\Pi  _j}=0\  , \\ &
	(Y ) _{R}= - P(p_0)J \Gamma (p)  R  -
\beta  _j(p,R)P(p_0)J P^*(p) \Diamond _j P(p)R \\& = \mathbf{S}^{1,1}  (p-p_0,R) +\mathcal{R}^{0,2}(p,R) P(p_0)P (p)J   \Diamond _j P(p)R  .
   \end{aligned}
\end{equation}
Equation  $
 i_{\mathcal{X}^t} \Omega _t=-  \alpha
   $   is equivalent to
\begin{equation} \label{eq:fred1}  \begin{aligned}
 & (1+ t\mathcal{K})\mathcal{X}^t=Y
\end{aligned}
\end{equation}
where the operator $\mathcal{K}$ is defined by $ i_{X}d  \alpha =
 i_{\mathcal{K} X} \Omega _0$. In coordinates,  \eqref{eq:fred1} becomes
   $
 (\mathcal{X}^t)_{\Pi _j} =0$ and, for $P=P(p)$,
\begin{align} &   (\mathcal{X}^t)_{\tau _j} + t \langle  \widehat{\Gamma} _j+
(  \widehat{\delta } _   k    \partial _{p_k}{\beta } _ {j } -\widehat{\delta } _   j    \partial _{p_j}{\beta } _ {k } )    \Diamond _k P  R, (\mathcal{X}^t)_{R}
\rangle  =-\beta _j  , \label{eq:fred20b} \\&    (\mathcal{X}^t)_{R} +t{\mathcal L}(\mathcal{X}^t)_{R}  =(Y)_{R}
 \  ,\text{where for $X\in N_g^\perp (\mathcal{H}_{p_0}^*)$} \label{eq:fred20t}\\&  {\mathcal L}X:=  P ( {p}_0)J\left [   2     \Gamma  X +
\langle \widetilde{\beta}  _j   , X\rangle P^*  \Diamond _j P R  -   \langle  P^* \Diamond _j P R ,X\rangle    \widetilde{\beta}  _j \right ]  \label{eq:fred20} .\end{align}

\eqref{eq:fred20} implies the following lemma.

\begin{lemma}\label{lem:fred11}
  We have, summing on repeated indexes,  with $i$ varying in some finite set,
\begin{equation} \label{eq:fred11}  \begin{aligned}
 &{\mathcal L}X = {\mathcal A}_j(X)  J\Diamond _j  R +{\mathcal B}_i(X) \Psi _i
\end{aligned}
\end{equation}
where:   $\Psi _i = \mathbf{S}^{0,0}(p-p_0,R)$;
for $L= {\mathcal A}_j, {\mathcal B}_i$, we have
  $ L \in C^\infty (\mathcal{U}_{\mathbf{d}} , B(L^2 ,\R ))$  with
\begin{equation}   \label{eq:fred12} \begin{aligned}
 &    L(X)  = L  _{j} \left \langle \Diamond _j  R,  X  \right \rangle + \langle \widetilde{L} ,X \rangle ,
\end{aligned}
\end{equation}
where we have $ \widetilde{L} =  \mathbf{S}^{1,0}(p-p_0,R)$  and $ L  _{j}\in \resto ^{0,0}(p-p_0,R)$.
 \end{lemma}
\proof Schematically, for $\widetilde{L} _i =\mathbf{S}^{0,0}(p-p_0,R)$ and  $\Psi _i  =\mathbf{S}^{0,0}(p-p_0,R)$ we have

\begin{equation*}  \begin{aligned} & P(p)  R =R -P_{N_g}(p ) R
=R+\sum _i
\langle \widetilde{L} _i,R \rangle \Psi _i  \, ,   \\&  P^*(p) \Diamond _k  R =\Diamond _k R-P_{N_g}^*(p) \Diamond _k  R = \Diamond _kR+\sum _i
\langle \widetilde{L} _i,R \rangle \Psi _i  .\end{aligned}\end{equation*}
Then $(P^*(p) \Diamond _k P(p) -\Diamond _k  )R=\mathbf{S}^{0,1}(p-p_0,R)$.

\noindent By the definition of $\widetilde{\beta}  _j$ we have
\begin{equation}  \begin{aligned} &   \widetilde{\beta}  _j    =\sum _k \widehat{\delta} _j  (\partial _{p_j}\beta  _k)    \Diamond _k R    +   \widehat{L}     \\& \widehat{L}  :=   \nabla _R \beta  _j+ \frac 12 J^{-1}  \widehat{\delta} _j  \partial _{p_j}     P(p)  R     +\sum _{k } \beta  _k  \partial _{p_j} ( P^*(p) \Diamond _k P(p)    )R\\& -\sum _{k }\widehat{\delta} _j \partial _{p_j}\beta  _k  \left [  P_{N_g}^*(p)\Diamond _k P (p)R  + \Diamond _k P_{N_g}(p) R \right ]  ,
\end{aligned}  \nonumber \end{equation}
where $\widehat{L}= \mathbf{S}^{0,1}(p-p_0,R)$.

\noindent We also have  $ \Gamma  X =\frac 12 J^{-1} (P_{N_g} (p_0)- P_{N_g} (p ) ) X =\sum _i
\langle \widetilde{L} _i, X \rangle \Psi _i $ with $\widetilde{L} _i= \mathbf{S}^{1,0}(p-p_0,R)$
and $\Psi _i = \mathbf{S}^{0,0}(p-p_0,R)$.
   This yields the result.

\qed

\begin{lemma}
  \label{lem:vectorfield} System \eqref{eq:fred20b}--\eqref{eq:fred20}
   admits exactly one solution $\mathcal{X}^t$.
   For $\mathcal{A }_j=\mathcal{R}^{0,2}_{n,\infty}(t,p-p_0,R)$,  $\mathcal{D } =\mathbf{S}^{1,1}_{n,\infty}(t,p-p_0,R)$ with $|t|<3$, we have
\begin{equation}\label{eq:quasilin1}
\begin{aligned} &
    (\mathcal{X}^t)_R = \mathcal{A }_j J\Diamond _j R  + \mathcal{D}  .   \end{aligned}
\end{equation}

\end{lemma}
\proof   Recall   $Y$   defined by $i_Y\Omega _0 =-\alpha $. By \eqref{eq:Yvec} with $\widetilde{\mathcal{A }}_j=\mathcal{R}^{0,2}_{n,\infty}(p-p_0,R)$ and  $\widetilde{\mathcal{D }} =\mathbf{S}^{1,1}_{n,\infty}(p-p_0,R)$
we have  $  ({Y})_{R} = \widetilde{\mathcal{A }}_jJ\Diamond _j R  + \widetilde{\mathcal{D} } $. By
$ (\mathcal{X}^t)_{R} +t{\mathcal L}(\mathcal{X}^t)_{R}  =(Y)_{R}$ and
Lemma \ref{lem:fred11} this implies for $X=(\mathcal{X}^t)_{R}$
\begin{equation}  \begin{aligned} &  \langle \Diamond _k R, X\rangle  +t\mathcal{B}_i(X)\langle \Diamond _k R, \Psi _i \rangle   =  \langle \Diamond _k R, (Y)_{R} \rangle  \\&  \langle  \widetilde{L}, X\rangle +t\mathcal{A }_j(X) \langle  \widetilde{L}, J\Diamond _j R\rangle   +t\mathcal{B}_i(X)\langle  \widetilde{L}, \Psi _i \rangle   =  \langle  \widetilde{L}, (Y)_{R} \rangle ,\end{aligned}\nonumber
\end{equation}
as $L$ runs through all the $L= {\mathcal A}_j, {\mathcal B}_i$.  Taking appropriate linear combinations of these
equations with the coefficients $L_j$  of  $L= {\mathcal A}_j, {\mathcal B}_i$, see Lemma \ref{lem:fred11},
for a   matrix $\textbf{R}^{0,1}(p-p_0, R)$ whose coefficients are $\resto^{0,1}(p-p_0, R)$,
we get  $$(1+t\textbf{R}^{0,1}(p-p_0, R)) \begin{pmatrix} \mathcal{A }_j((\mathcal{X}^t)_{R}) \\   \mathcal{B}_i((\mathcal{X}^t)_{R})\end{pmatrix} =   \begin{pmatrix}      \mathcal{A }_j((Y)_{R}) \\   \mathcal{B}_i((Y)_{R})\end{pmatrix}.$$
Then we get

\begin{equation} \label{eq:fredcorr1} \begin{aligned} &  \begin{pmatrix} \mathcal{A }_j((\mathcal{X}^t)_{R}) \\   \mathcal{B}_i((\mathcal{X}^t)_{R})\end{pmatrix} =(1+t\textbf{R}^{0,1}(p-p_0, R))^{-1}  \begin{pmatrix}      \mathcal{A }_j((Y)_{R}) \\   \mathcal{B}_i((Y)_{R})\end{pmatrix}. \end{aligned}
\end{equation}
Using the left hand side of
\eqref{eq:fredcorr1}  set

\begin{equation}\label{eq:fredcorr2}
 {\mathcal L}(\mathcal{X}^t)_{R} := {\mathcal A}_j((\mathcal{X}^t)_{R})  J\Diamond _j  R +{\mathcal B}_i((\mathcal{X}^t)_{R}) \Psi _i.
\end{equation}
The rhs of \eqref{eq:fredcorr2} satisfies the properties stated for the rhs of
\eqref{eq:quasilin1}. Finally  set
$ (\mathcal{X}^t)_{R}   :=(Y)_{R}- t{\mathcal L}(\mathcal{X}^t)_{R}$.
This is  a solution of \eqref{eq:fred20t}. It is elementary to see from the argument that such solution is unique and that it satisfies
the properties of the statement.  \qed

Turning to coordinates $(\tau , \Pi , R)$ and by
 Lemma \ref{lem:var} we conclude what follows.
\begin{lemma}\label{lem:fred12}
  Consider the coordinate system $(\tau ,\Pi ,R)$.
For  $G$ any of the   $\mathcal{A}  _{j}$,  $ \mathcal{D} $  in Lemma
\ref{lem:vectorfield},
 we have $G=G ( \Pi ,   \Pi (R),R )$, with  $G ( \Pi ,   \varrho ,R )$
      smooth w.r.t.  $ (\Pi,\varrho ,R)\in   \mathcal{U}_{ \textbf{d}}
$, with    $\mathcal{U}_{ \textbf{d}}$  formed by the $ (\Pi,\varrho ,R)\in \R ^{2n_0}
\times (\Sigma _{\textbf{d}}\cap N^\perp _g (\mathcal{H}_{p_0}))$
  defined by the inequalities $  \|R\| _{\Sigma  _{-n}}\le   \varepsilon $,
   $|\varrho | \le   \varepsilon  $ and
   $| \Pi  -p_0|\le   \varepsilon  $  for $\varepsilon >0$ small enough.
\end{lemma}

\subsection{Flows}
\label{subsec:ode}

The following lemma is repeatedly used in the sequel, see Lemma 3.24
\cite{bambusi}.

\begin{lemma} \label{lem:ODE}  Below we pick $r,M,M_0,s,s',k,l\in \N\cup \{ 0 \}$ with $1\le l\le M$.
Consider  a system

\begin{equation} \label{eq:ODE}\begin{aligned} &
  \dot \tau   _j = T_j (t,\Pi , \Pi (R), R  ) \  , \quad \dot \Pi  _j    =0 \  ,  \\&  \dot R
	= \mathcal{A }_j(t,\Pi , \Pi (R), R   ) J\Diamond _j R  + \mathcal{D}(t,\Pi , \Pi (R), R   )
,
\end{aligned}   \end{equation}
  where  we assume what follows.

\begin{itemize}
\item $P_{N  _g ( p_0)} (\mathcal{A }_j J\Diamond _j R  + \mathcal{D}) \equiv 0$.

\item  At $\Pi =p_0$, dropping the dependence  on $\Pi$ and for $\U _{-r}$ a neighborhood
of  0 in $\Ph  ^{-r}$,  we have
$\mathcal{A }(t, \varrho ,R ) \in C^M ((-3,3)\times \U _{-r},  \R ^{n_0}  )$ and $ \mathcal{D}(t, \varrho ,R ) \in C^M ( (-3,3)\times\U _{-r},  \Sigma _{r}  )$

\item In $  (-3,3)\times \U _{-r}$  for a fixed $i$ in $ \{ 0,1\}$, and a fixed $C_r$, we have:

 \begin{equation} \label{eq:symbol}   \begin{aligned} &   | \mathcal{A }(t, \varrho ,R ) |\le C \| R\| _{\Sigma _{-r}}^{M_0+1}    ,  \\&   \| \mathcal{D }(t, \varrho ,R )\| _{\Sigma _{r}} \le C    (|\varrho| +\| R\| _{\Sigma _{-r}})^i  \| R\| _{\Sigma _{-r}}^{M_0 } .
\end{aligned}    \end{equation}

 \end{itemize}
Let  $k\in \Z\cap [0,r-(l+1)\textbf{d}]$ and set for  $s^{\prime \prime}\ge \textbf{d}$ (or $s^{\prime \prime}\ge \textbf{d}/2$  if $\textbf{d}/2\in \N $)
\begin{equation} \label{eq:domain0}   \begin{aligned}   &   \U _{\varepsilon _1,k}^{s^{\prime \prime}} :=
\{ (\tau, \Pi  , R) \in \widetilde{ {\Ph}}^{s^{\prime \prime}} \  : \  \Pi =p_0 \  ,  \  \|  R \| _{\Sigma _{-k }} + |\Pi (R)| \le \varepsilon _1\}.
\end{aligned}    \end{equation}
Then for $\varepsilon _1>0$ small enough,
	the initial value
	problem associated to \eqref{eq:ODE}  for $\Pi =p_0$
	defines a      flow $\mathfrak{F} ^t =(  \mathfrak{F}  ^t  _{\tau    } ,  \mathfrak{F}  ^t   _{R} )$
	  for $t\in [-2, 2]$   in   $\U _{\varepsilon _1,k}^{\textbf{d}}$.  In particular
for $\Pi =p_0$, for  $R$ in
 a  neighborhood $  B_{\Sigma _{-k}}$   of
0 in $\Sigma _{-k}$ and $\Pi (R)   $  in a  neighborhood $ B_{\R ^{ n_0}} $ of
0 in $\R ^{ n_0}$,   we have
	\begin{equation} \label{eq:ODE1}\begin{aligned} &
\mathfrak{F}  ^t   _{R} (  \Pi (R), R) = e^{Jq(t,  \Pi (R), R)\cdot \Diamond } ( R+ \textbf{S} (t,  \Pi (R), R    ))
,
\end{aligned}   \end{equation}
\begin{equation} \label{eq:ODEpr210}\begin{aligned} \text{with   }  &
     \textbf{S}  \in C^l((-2,2)
\times B_{\R ^{ n_0}}\times   B_{\Sigma _{-k}} , \Sigma _{r- (l+1)\textbf{d}}
 )    \\&  {q}  \in C^l((-2,2)
\times B_{\R ^{ n_0}} \times   B_{\Sigma _{-k}}, \R ^{ n_0}
 ).
\end{aligned}   \end{equation}
For fixed $C>0$ we have
 \begin{equation} \label{eq:symbol1}   \begin{aligned}   &   | q (t,\varrho ,R ) |\le C \| R\| _{\Sigma _{(l+1)\textbf{d}-r}}^{M_0+1}  \, , \\&     \| \textbf{S} (t,\varrho ,R ) \| _{\Sigma _{r- (l+1)\textbf{d}}} \le C   (|\varrho| +\| R\| _{\Sigma _{(l+1)\textbf{d}-r}})^i \| R\| _{\Sigma _{(l+1)\textbf{d}-r}}^{M_0 }  .
\end{aligned}    \end{equation}
Furthermore we have  $\textbf{S}= \textbf{S}_1+ \textbf{S} _2$  with
\begin{equation} \label{eq:duhamel}   \begin{aligned}   &  \textbf{S}_1 (t,  \Pi (R), R    ) =
\int _0^t \mathcal{D }(t',\Pi (R (t') ) ,R (t')    ) dt'\\&   \|   \textbf{S}_2 (t,  \varrho, R   )    \| _{\Sigma _{s}}
\le C\| R  \| _{\Sigma _{(l+1)\textbf{d}-r} } ^{2M_0+1} (|\varrho  | +\| R \| _{\Sigma _{(l+1)\textbf{d}-r}} )^i .
\end{aligned}    \end{equation}
For $r-(l+1)\textbf{d}\ge  s' \ge  s+l\mathbf{d}\ge l\mathbf{d}$ and    $k\in \Z\cap [0,r-(l+1)\textbf{d}]$ and
     for  $ \varepsilon _1>0$ sufficiently small, we have
 \begin{equation} \label{eq:reg1}\begin{aligned} &\mathfrak{F}  ^t \in C^l((-2,2)\times \U _{\varepsilon _1,k}^{s'}
  , \widetilde{\Ph}  ^{s }
 )
.
\end{aligned}   \end{equation}
Furthermore, there exists $ \varepsilon _2>0$
such that
\begin{equation} \label{eq:main1}\begin{aligned} &  \mathfrak{F}  ^t( \U _{\varepsilon _2,k}^{s'})\subset  \U _{\varepsilon _1,k}^{s'}  \text{ for all $|t|\le 2$
.}
\end{aligned}   \end{equation}

\noindent  We have      \begin{equation} \label{eq:ODE11}\begin{aligned} &
\mathfrak{F}  ^t   (e^{J\tau \cdot \Diamond } U)  \equiv e^{J\tau \cdot \Diamond }\mathfrak{F}  ^t   (  U)
.
\end{aligned}   \end{equation}
	\end{lemma}
\proof
It is enough to focus on the equation for $R$.   Set $S= e^{-Jq \cdot \Diamond }R$ for $q\in \R ^{n_0}$.  Then  consider the following system:

\begin{equation} \label{eq:ODEpr}\begin{aligned} &
    \dot S
	= e^{-Jq \cdot \Diamond } \mathcal{D}(  t,\varrho ,e^{ Jq \cdot \Diamond }S   ) \  ,
\\& \dot q = \mathcal{A } (t,\varrho ,e^{ Jq \cdot \Diamond }S    ) \quad , \quad q(0)=0  ,  \\& \dot \varrho _j= \langle S,  e^{-Jq \cdot \Diamond } \Diamond _j\mathcal{D}( t, \varrho , e^{ Jq \cdot \Diamond }S   ) \rangle  \ .
\end{aligned}   \end{equation}
For $ l\le M$ and $    k, s^{\prime \prime}\in [0, r-(l+1)\textbf{d}]$ the field in \eqref{eq:ODEpr} is $C^l( (-3,3)\times \U _{-k},  \Sigma _{ s^{\prime \prime}}
\times \R ^{2n_0})$  with
$\U _{-k}\subset
\Sigma _{-k}
\times \R ^{2n_0} $ a neighborhood of the equilibrium 0.
 This follows from the fact that $(q,X)\to  e^{ Jq \cdot \Diamond }X$
 is in $C^{l} ( \R ^{ n_0} \times \Sigma _{ \ell} ,\Sigma _{ \ell -l\mathbf{d}} )$ for all $\ell \in \Z$ and from the hypotheses on $\mathcal{A}$ and $ \mathcal{D}$. For example
 \begin{equation*}\begin{aligned} &
  (t,q,\varrho ,S) {\rightarrow}  e^{-Jq \cdot \Diamond } \Diamond _j \mathcal{D}( t, \varrho , e^{ Jq \cdot \Diamond }S   ) \in C^l ( (-3,3)\times \R ^{2n_0} \times \Sigma _{  l\mathbf{d}-r}  , \Sigma _{   r-(l+1)\mathbf{d}}
   ), \end{aligned}
 \end{equation*}
 (more precisely for $(q,\varrho ,S)$ in a neighborhood of the origin). So
 \begin{equation*}\begin{aligned} &
  (t,q,\varrho ,S) {\rightarrow} \langle S, e^{-Jq \cdot \Diamond } \Diamond _j \mathcal{D}( t, \varrho , e^{ Jq \cdot \Diamond }S   )\rangle  , \end{aligned}
 \end{equation*}
 is in
 $ C^l ( (-3,3)\times \R ^{2n_0} \times \Sigma _{  -k}  , \R
   )$ for $k\le     r -  (l+1)\mathbf{d} $ (for $(q,\varrho ,S)$ near origin).

 \noindent
 For $l\ge 1$ we can apply  to \eqref{eq:ODEpr} standard
theory of ODE's to conclude that  there are neighborhoods of the origin $B_{\R ^{2n_0}}\subset \R ^{2n_0}$ and $B_{\Sigma _{-k}}\subset  \Sigma _{-k} $
such that   the flow  is of the form

\begin{equation} \label{eq:ODEpr1}\begin{aligned} &
      S(t)
	=  R+ \textbf{S} (t, \varrho , R )  \ , \quad  \textbf{S} (0, \varrho ,R) =0\ ,  \\& q(t)=  {q} (t,\varrho , R ) \ , \quad  {q}(0, \varrho , R) =0\ ,  \\& \varrho (t)=\varrho + \overline{\varrho} (t, \varrho , R)  \   ,   \quad   \overline{\varrho} (0,\varrho , R ) =0\ ,
\end{aligned}   \end{equation}
\begin{equation} \label{eq:ODEpr2}\begin{aligned} \text{with   }  &
     \textbf{S}  \in C^l((-2,2)
\times B_{\R ^{ n_0}}\times   B_{\Sigma _{-k}} , \Sigma _{r- (l+1)\textbf{d}}
 )    \\&   \overline{\varrho},{q} (t,\varrho , R )  \in C^l((-2,2)
\times B_{\R ^{ n_0}} \times   B_{\Sigma _{-k}}, \R ^{ n_0}
 ).
\end{aligned}   \end{equation}

 \noindent For $S\in \Sigma  _{ {\textbf{d}}  }\cap B_{\Sigma _{-k}}$ and $S(0)=S$,
choosing $s^{\prime \prime}\ge \textbf{d}$  we have $S(t)\in\Sigma  _ {\textbf{d}}  $
with $\Pi (S(t))=\varrho (t)$  for  $\varrho (0)=\varrho =\Pi (S)$.  Then \eqref{eq:ODEpr2} yields
\eqref{eq:ODEpr210}  (we can replace $\Sigma  _{ {\textbf{d}}  }$ with $\Sigma  _{ \frac{{\textbf{d}} }{2} }$ if $ \frac{\textbf{d} }{2}\in \N$). \eqref{eq:ODE1} and \eqref{eq:ODEpr210} yield \eqref{eq:reg1}.

\noindent We have  for $R(0)=R$
\begin{equation} \label{eq:duh2}\begin{aligned} & R(t)= e^{ Jq (t) \cdot \Diamond } (R+\int _0^t
e^{- Jq (t') \cdot \Diamond  }  \mathcal{D}( t',   \varrho (t'),R (t')  ) dt' ).
\end{aligned}   \end{equation}
By (A6),  for $\epsilon =0$,  and by \eqref{eq:symbol}, for  $   |s ^{\prime\prime}|\le    r-(l+1)\textbf{d} $ we have
\begin{equation} \label{eq:duh21}\begin{aligned} & \| R(t) \| _{\Sigma _{ s ^{\prime\prime}} }\le   C\| R  \| _{\Sigma_{ s ^{\prime\prime}} }+C\int _0^t  \| \mathcal{D}(  t',  \varrho (t),R (t')  ) \| _{\Sigma _{ r} } dt'    \\& \le  C\| R  \| _{\Sigma _{ s ^{\prime\prime}} }
+C\int _0^t   \| R (t')\| _{\Sigma _{-r}}^{ M_0 } (|\varrho (t')| +\| R (t')\| _{\Sigma _{-r}}) ^i dt'     \\& \le  C\| R  \| _{\Sigma _{ s ^{\prime\prime}} }
+C\int _0^t   \| R (t')\| _{\Sigma _{ s ^{\prime\prime}} }^{ M_0 } (|\varrho (t')| +\| R (t')\|_{\Sigma _{ s ^{\prime\prime}} }) ^i dt',
\end{aligned}   \end{equation}
with the caveat that the second line is purely formal and is used to
get the third line, where the integrand is continuous.
Proceeding similarly, for $\varrho (0)  = \varrho$
\begin{equation} \label{eq:duh22}\begin{aligned} &  | \varrho (t) - \varrho    |
\le  \int _0^t | \langle R(t'), \Diamond \mathcal{D}( t',  R (t'), \varrho (t')  ) \rangle |  dt'\\&  \le  \int _0^t  \| R (t')\| _{\Sigma _{(l+1)\textbf{d} -r}}  \| \mathcal{D}( t',  \varrho (t),R (t')  ) \| _{\Sigma _{r-l\mathbf{d}}}
  dt'\\&
 \le
C\int _0^t  \| R (t')\| _{\Sigma _{(l+1)\textbf{d} -r}}   ^{ M_0+1 }
(|\varrho (t')| +\| R (t')\| _{\Sigma _{(l+1)\textbf{d} -r}}) ^i     dt'.
\end{aligned}   \end{equation}
So for  $| s ^{\prime\prime}|\le r-(l+1)\textbf{d}$, using the continuity
in $t'$ of the integrals in the last lines of \eqref{eq:duh21} and \eqref{eq:duh22},
by the Gronwall inequality  there is  a fixed $C$ such that for all $|t|\le 2$  we have
\begin{align} \label{eq:gronwall0} & \| R(t) \| _{\Sigma _{ s ^{\prime\prime}} }\le   C\| R  \| _{\Sigma _{ s ^{\prime\prime}} } \, , \\&    | \varrho (t) - \varrho    |    \le  C  \| R \| _{\Sigma _{(l+1)\textbf{d} -r}}^{ M_0+1 }
(|\varrho | +\| R \| _{\Sigma _{(l+1)\textbf{d} -r}}) ^i . \label{eq:gronwall1}
\end{align}
 By
\eqref{eq:gronwall0}  for $s ^{\prime\prime}=s'$ and   $s ^{\prime\prime}=-k$ and  by
$    | \varrho (t) - \varrho    |    \le  C  \| R \| _{\Sigma _{-k}}^{ M_0+1 }
(|\varrho | +\| R \| _{\Sigma _{-k}}) ^i$, we get   $\mathfrak{F}  ^t( \U _ {\varepsilon _2,k}^{s^{\prime  }})\subset  \U _ {\varepsilon _1,k}^{s^{\prime  }}$ for all $|t|\le 2$ for   $\varepsilon _1 \gg \varepsilon _2,$ that is \eqref{eq:main1}.

 \noindent We have \begin{equation*}  \begin{aligned} & S(t,\varrho ,R)=  \int _0^t
e^{- Jq (t') \cdot \Diamond  }  \mathcal{D}( t',   \varrho (t'),R (t')  ) dt' ),
\end{aligned}   \end{equation*}
 Proceeding as for \eqref{eq:duh21} and using  \eqref{eq:gronwall0}--\eqref{eq:gronwall1}
 we get the estimate for $\textbf{S}$  in
\eqref{eq:symbol1}. The estimate on $q$ is obtained similarly
integrating the second equation in \eqref{eq:ODEpr1}.

\noindent We have
\begin{equation} \label{eq:duh3}\begin{aligned} & \textbf{S}_2 (t,  R, \varrho    )  =  \int _0^1 dt^{\prime \prime }\int _0^t
  e^{- t^{\prime \prime }q (t') \cdot \Diamond  }  q (t') \cdot \Diamond \mathcal{D}( t',   \varrho (t),R (t')  ) dt'
\end{aligned}   \end{equation}
Then by \eqref{eq:gronwall0}--\eqref{eq:gronwall1} we get
 \begin{equation} \label{eq:duh4}\begin{aligned} & \| \textbf{S}_2 (t,  R, \varrho    )  \| _{\Sigma _{r-\textbf{d}}}\le C ^{\prime \prime}    \int _0^t
    |q (t')  |      \| \mathcal{D}(t',    \varrho (t),R (t'))  \| _{\Sigma _{r-\textbf{d} }}  dt'    \\&  \le C'      \int _0^t
     \| R (t')\|  _{\Sigma _{(l+1)\textbf{d}-r} }^{ 2M_0+1 } (|\varrho (t')| +\| R (t')\|  _{\Sigma _{(l+1)\textbf{d}-r} })^i dt'  \\&\le  C
		\| R  \| _{\Sigma _{(l+1)\textbf{d}-r} }^{2M_0+1} (|\varrho  | +\| R \| _{\Sigma _{(l+1)\textbf{d}-r}} )^i.
\end{aligned}   \end{equation}
This yields \eqref{eq:duhamel}.  \eqref{eq:reg1} follows by \eqref{eq:ODE1}--\eqref{eq:ODEpr210}.
Finally,
\eqref{eq:ODE11}  follows immediately from \eqref{eq:ODE1}.

\qed

\begin{lemma}
  \label{lem:ODEdomains} Assume hypotheses and conclusions of Lemma  \ref{lem:ODE}.
	Consider the  flow of system    \eqref{eq:ODEpr}   for $\Pi =p_0$ .  Denote  the flow in the space	 with variables $\{  ( \varrho  ,R)\}$ by  $\mathfrak{F} ^t =(     \mathfrak{F}  ^t   _{\varrho } ,  \mathfrak{F}  ^t   _{R} )$.   Then we have
	   \begin{equation} \label{eq:ODEdomains1} \begin{aligned} &
\mathfrak{F}  ^t   _{R} (  \varrho  ,R ) = e^{Jq(t,  \varrho  ,R)\cdot \Diamond } ( R+ \textbf{S} (t,  \varrho  ,R    )) \\&  \mathfrak{F}  ^t   _{\varrho} (  \varrho  ,R ) = \varrho + \overline{\varrho} (t, \varrho , R)
.
\end{aligned}   \end{equation}
Furthermore, the following facts hold.

\begin{itemize}
\item[(1)]
Let $k \in \Z \cap [0,r-(l+1)\textbf{d}]$ and $h\ge \max\{ k +l\textbf{d},(2 l+1)\textbf{d}-r\}$. Then  we have
 $\mathfrak{F}  ^t \in C^l((-2,2)\times \U _{-k  }
  , \Ph ^{-h }
 )
$
for a
neighborhood  of the origin
$\U _{-k }\subset \Ph ^{-k }$.

\item[(2)]   Let $h$ and $k$ be like above with $h\le r-( l+1)\textbf{d}$. Then given a function $\resto ^{a,b}_{h, l}( \varrho  ,R)$,  we have $\resto ^{a,b}_{h, l}\circ \mathfrak{F}  ^t=\resto ^{a,b}_{k, l}( t, \varrho  ,R) $ and given  a function $\textbf{S}^{a,b}_{h, l}( \varrho  ,R)$,  we have $\textbf{S} ^{a,b}_{h, l}\circ \mathfrak{F}  ^t=\textbf{S} ^{a,b}_{k, l} ( t,\varrho  ,R)$.

\end{itemize}
	\end{lemma}
\proof   \eqref{eq:ODEdomains1} follows  by \eqref{eq:ODEpr1}.
By \eqref{eq:ODEpr2} we have\begin{equation*}  \begin{aligned}    &
     \textbf{S}  \in C^l((-2,2)
\times \U _{-k } , \Sigma _{r- (l+1)\textbf{d}}
 )    \ , \quad   q \text{ and }\mathfrak{F}  ^t   _{\varrho} \in C^l((-2,2)
\times \U _{-k }, \R ^{ n_0}
 ).
\end{aligned}   \end{equation*}
By the above formulas we have $\mathfrak{F}  ^t   _{R}\in C^l((-2,2)
\times \U _{-k }, \Sigma _{r- (2l+1)\textbf{d}} \cap \Sigma _{-k- l\textbf{d}}
 ).$  This yields   $\mathfrak{F}  ^t   _{R}\in C^l((-2,2)
\times \U _{-k }, \Sigma _{-h}
 ) $  and yields Claim (1).

\noindent By Claim (1),     $\resto ^{a,b}_{h, l}\circ \mathfrak{F}  ^t \in C^l( (-2,2)
\times  \U _{-k }, \R ^{ n_0}
 )$. Let $(\varrho ^t , R^t)=\mathfrak{F}  ^t(\varrho  , R )$. Then
\begin{equation*}  \begin{aligned}    &
   |\resto ^{a,b}_{h, l}\circ \mathfrak{F}  ^t(\varrho  , R )|=| \resto ^{a,b}_{h, l}(\varrho ^t , R^t )| \le C' \|  R^t\| _{\Sigma   _{-h}}^b (\|  R^t\| _{\Sigma   _{-h}}+|\varrho ^t|)^{a}\\& \le  C \|  R\| _{\Sigma   _{-h}}^b (\|  R\| _{\Sigma   _{-h}}+|\varrho |)^{a}\le  C \|  R\| _{\Sigma   _{-k}}^b (\|  R\| _{\Sigma   _{-k}}+|\varrho |)^{a},
\end{aligned}   \end{equation*}
 where the first inequality uses Definition \eqref{eq:scalSymb}, the second uses \eqref{eq:gronwall0}--\eqref{eq:gronwall1} for $s^{\prime\prime}=-h$ and the last is obvious. Similarly  by  Claim (1),     $\textbf{S} ^{a,b}_{h, l}\circ \mathfrak{F}  ^t \in C^l(  (-2,2)
\times \U _{-k }, \Sigma _h
 )\subset  C^l( (-2,2)
\times  \U _{-k }, \Sigma _k )$ and
\begin{equation*}  \begin{aligned}    &
    \| \textbf{S} ^{a,b}_{h, l}(\varrho ^t , R^t )\| _{ \Sigma _k}\le  \| \textbf{S} ^{a,b}_{h, l}(\varrho ^t , R^t )\| _{ \Sigma _h} \le C' \|  R^t\| _{\Sigma   _{-h}}^b (\|  R^t\| _{\Sigma   _{-h}}+|\varrho ^t|)^{a}\\& \le  C \|  R\| _{\Sigma   _{-h}}^b (\|  R\| _{\Sigma   _{-h}}+|\varrho |)^{a}\le  C \|  R\| _{\Sigma   _{-k}}^b (\|  R\| _{\Sigma   _{-k}}+|\varrho |)^{a}.
\end{aligned}   \end{equation*}

 \qed

\bigskip To prove Theorem \ref{th:main} we will need more information on $(\Pi (R (1)) , R(1)) $. This is provided by the following
lemma.

\begin{lemma} \label{lem:ODEbis}  Consider,  for $\mathcal{D}$ the function in
\eqref{eq:ODE} at $\Pi =p_0$, the system
\begin{equation} \label{eq:ODEbis}\begin{aligned} &   \dot S(t)
	=  \mathcal{D}(  t, \Pi (R_0),S(t)  ) \, , \quad S(0)=R_0 .
\end{aligned}   \end{equation}
Then  for $S'=S(1)$ and  for $R'=R(1)$ with $R(t)$ the solution of  \eqref{eq:ODE}
with $R(0)=R_0$, we have (same indexes of Lemma \ref{lem:ODE})
\begin{equation} \label{eq:ODE1bis}\begin{aligned} &
    \| R'-S'\|  _{\Sigma _{-s' }}   \le  C  \| R _0\|  _{\Sigma _{ -s }}^{M_0+2}  \, , \\&
      \Pi (R')- \Pi (S')=\resto ^{i,2M_0+1}_{s , l}
( \Pi (R_0), R_0).
	 \end{aligned}   \end{equation}
\end{lemma}
\proof  Recall that for $\varrho =\Pi (R)$ we have $\dot \varrho =\langle R , \Diamond  \mathcal{D}( t, \varrho ,R   ) \rangle .$ Similarly, for  $\sigma  =\Pi (S)$ we have  $\dot \sigma = \langle S,\Diamond  \mathcal{D}(t,   \varrho _0  ,S) \rangle $, where $\varrho _0=\Pi (R_0)$. So  we have
\begin{equation*} \label{eq:ODE2bis}\begin{aligned}    \dot \varrho -\dot \sigma  & =
\langle R, \Diamond  \mathcal{D}(  t,\varrho ,R    ) \rangle  -  \langle S,\Diamond   \mathcal{D}(t,    \varrho _0,S   ) \rangle  \\&  =\langle R -S, \Diamond  \mathcal{D}(t,  \varrho ,R   ) \rangle +\langle S,  \Diamond  ( \mathcal{D}(t,   \varrho _0 ,S  ) -\mathcal{D}(t,  \varrho ,R    )) \rangle .
	 \end{aligned}   \end{equation*}
By  \eqref{eq:symbol}   for   fixed constants
and using $s'\le r- \mathbf{d}$,   we have

\begin{equation*} \label{eq:ODE3bis}\begin{aligned} &  |\dot \varrho -\dot \sigma   |\lesssim \| R-S\|  _{\Sigma _{-s ' }}  \| \mathcal{D}(  t, \varrho ,R   ) \|  _{\Sigma _{r }} + \|  S\|  _{\Sigma _{ -s'   }}
\| \mathcal{D}( t,   \varrho _0 ,S  ) -\mathcal{D}(t,    \varrho ,R  ) \|  _{\Sigma _{r }}
\\&   \lesssim \| R-S\|  _{\Sigma _{-s'  }}  \| R \|  _{\Sigma _{ -s   }}^{M_0}  (|\varrho |+ \| R \|  _{\Sigma _{ -s '  }}) ^i+|  \varrho -  \varrho _0    | \  \|  S\|  _{\Sigma _{ -s  ' }}  \|(R,  S)\|  _{\Sigma _{ -s'  }}^{M_0}
\\& + \| R-S\|  _{\Sigma _{-s' }}  \|  S\|  _{\Sigma _{ -s '  }}  \|(R,  S)\|  _{\Sigma _{ -s'}}^{M_0-1}  (|  (\varrho ,  \varrho _0  )   |+ \|(R,  S)\|  _{\Sigma _{ -s'}})^i.
	 \end{aligned}   \end{equation*}
We have
$ \dot R -\dot S  = \mathcal{D}(t,   \varrho ,R  ) -\mathcal{D}(t,   \varrho _0  ,S ) + J  \mathcal A(t,   \varrho ,R  ) (t,   \varrho ,R  ) \cdot \Diamond R$
and hence  for   fixed constants  we have, using $s\le s'-\mathbf{d}$,

\begin{equation*} \label{eq:ODE4bis}\begin{aligned} &  \| R-S\|  _{\Sigma _{-s ' }}   \le \int _0^t [ \|
\mathcal {D}(   \varrho , R) -\mathcal{D}(    \varrho _0 ,S ) \|  _{\Sigma _{-s'  }} + |{\mathcal A} |  \|   R \|  _{\Sigma _{-s  }} ]dt'   \\& \lesssim  \int _0^t \big [
\| R-S\|  _{\Sigma _{-s'  }}   \|(R,  S)\|  _{\Sigma _{ -s' }}^{M_0-1}  (|  (\varrho ,\varrho _0 )   |+ \|(R,  S)\|  _{\Sigma _{ -s'}})^i \\&  + | \varrho - \varrho _0 |   \  \|(R,  S)\|  _{\Sigma _{ -s' }}^{M_0}
 +   \| R \|  _{\Sigma _{ -s }}^{M_0+2}  \big ] dt' .
	 \end{aligned}   \end{equation*}
	Recall that $ | \varrho - \varrho _0 |  \le C  \| R _0 \|  _{\Sigma _{(l+1)\textbf{d}-r   }} ^{M_0+1}(|\varrho _0|+ \| R_0 \|  _{\Sigma _{ (l+1)\textbf{d}-r   }}) ^i $  by \eqref{eq:gronwall1}, that $s<r-(l+1)\textbf{d}$ and that we have \eqref{eq:gronwall0} for $s^{\prime \prime}=-s,-s'$.
	Then by Gronwall inequality,  the above inequalities yield

\begin{equation} \label{eq:ODE6bis}\begin{aligned} &
 \| R(t)-S(t)\|  _{\Sigma _{-s' }}   \le  C  \| R _0\|  _{\Sigma _{ -s }}^{M_0+2}   \\&
|  \varrho (t)-  \sigma   (t)|\le    C     \| R _0\|  _{\Sigma _{-s   }}^{ 2M_0+1 }    (|\varrho _0|+ \| R_0 \|  _{\Sigma _{ -s   }}) ^i .
	 \end{aligned}   \end{equation}
This yields the  bounds implicit in  \eqref{eq:ODE1bis}. The regularity follows from  Lemma \ref{lem:ODE}.

\qed

\subsection{Darboux Theorem: end of the proof}
\label{subsec:darboux}

Formally the proof should follow by $ i_{\mathcal{X}  ^t} \Omega  _t=-\alpha$, where $ \Omega _t =(1-t)\Omega _0+t
\Omega$, and by
\begin{equation}\label{eq:fdarboux}  \begin{aligned} &
\frac{d}{dt}
\left ( \mathfrak{F}_{  t}^*\Omega _t\right )     =   \mathfrak{F}_{  t}^*
\left (L_{\mathcal{ X}_t} \Omega _t+\frac{d}{dt}\Omega _t\right )     =  \mathfrak{F}_{  t}^*
\left ( d i_{\mathcal{X}  ^t} \Omega  _t+d\alpha \right )      =0.
\end{aligned}
\end{equation}
But while for  \cite{boussaidcuccagna,Cu2} the above formal
computation falls within the classical framework of flows, fields and
differential forms,    in  the case of  \cite{bambusi,Cu1}
 this is not rigorous. In order to   justify rigorously this computation, we will consider
first a   regularization of  system \eqref{eq:ODE}.
\begin{lemma} \label{lem:moll}  Consider  the system

\begin{equation} \label{eq:ODEmoll}\begin{aligned} &
  \dot \tau   _j = T_j (t,\Pi , \Pi (R), R  ) \  , \quad \dot \Pi  _j    =0 \  ,  \\&  \dot R
	= \mathcal{A }_j(t,\Pi , \Pi (R), R   ) J\langle \epsilon \Diamond\rangle ^{-2}\Diamond _j R  + \mathcal{D}_ \epsilon    (t,\Pi , \Pi (R), R   )
,
\end{aligned}   \end{equation}
  where   $\mathcal{D}_  \epsilon =\mathcal{D}+\mathcal{A }_j P_{N  _g ( p_0)}   J\Diamond _j ( 1- \langle \epsilon \Diamond\rangle ^{-2})R  . $

\begin{itemize}
\item[(1)]
	For $|\epsilon  |\le 1$  system \eqref{eq:ODEmoll}  satisfies all the conclusions of Lemma \ref{eq:ODEmoll},
	if we replace $\Diamond $ in \eqref{eq:ODE1}  with $\langle \epsilon \Diamond\rangle ^{-2}\Diamond $
	(resp. ${\mathcal D}$ in \eqref{eq:duhamel} with  ${\mathcal D} _\epsilon$),
	with a fixed
	choice of constants  $\varepsilon _1  $, $\varepsilon _2  $, $C$,
			and with a fixed choice of   sets $B_{\R^{n_0}}$, $B_{\Sigma _{-s}}$.
		
		\item[(2)] For ${\mathcal X}^t$ the vector field of  \eqref{eq:ODE}, denote by ${\mathcal X}^t_\epsilon$
		the vector field of  \eqref{eq:ODEmoll}. Let
$ n' >  n+\textbf{d}$  with $n,n'\in \N$.  Then for $k\in \Z \cap [0,r]$
we have
		\begin{equation} \label{eq:reg3}\begin{aligned} & \lim _{\epsilon \to 0}{\mathcal X}^t_\epsilon  = {\mathcal X}^t \text{  in } C^{M}((-3,3)\times     \U _ {\varepsilon _0,k}^{n'}
  , \widetilde{ {\Ph} }^{n }
 )  \text{   uniformly locally}
,
\end{aligned}   \end{equation}
	that is uniformly on subsets of $(-3,3)\times   \U ^{n'}_{\varepsilon _0,k}$ bounded  in $(-3,3)\times   \widetilde{ {\Ph} }^{n '}$.

			\item[(3)]	 Denote  by $\mathfrak{F}  ^t_\epsilon =(\mathfrak{F}  ^t_{\epsilon \tau},\mathfrak{F}  ^t_{\epsilon R}) $ the flow associated to  \eqref{eq:ODEmoll} at $\Pi =p_0$.
Let $s'$,$s$ and $k$ as in the statement of Lemma \ref{lem:ODE}. Then there is  a pair $0<\varepsilon _ 1<\varepsilon _0$ such that
		\begin{equation} \label{eq:reg2}\begin{aligned} & \lim _{\varepsilon \to 0}\mathfrak{F}  ^t_\epsilon  = \mathfrak{F}  ^t \text{  in } C^{l-1}([-1,1]\times     \U _ {\varepsilon _1,k}^{s^{\prime  }}
  ,  \U _ {\varepsilon _0,k}^{s }
 )  \text{   uniformly locally}
.
\end{aligned}   \end{equation}

\end{itemize}
			
	\end{lemma}
\proof  For claim (1), it  is enough to check that ${\mathcal D} _\epsilon$ satisfies
an estimate like the one of  ${\mathcal D} $ in  \eqref{eq:symbol1} for a fixed $C$ for all $|\epsilon |\le 1$.
Indeed, after this has been checked, the proof of  Lemma \ref{eq:ODE} can be repeated verbatim,
exploiting (A6) for $\epsilon \neq 0$ and with $\Diamond$ replaced by $\langle \epsilon \Diamond\rangle ^{-2}\Diamond$.

\noindent  The   estimate on  ${\mathcal D} _\epsilon$ needed for Claim (1)
follows by
the definition of  ${\mathcal D} _\epsilon$ ,  by the estimate on ${\mathcal D} $,   by
$P_{N  _g ( p_0)}=  \textbf{e}_a\langle \textbf{e}^*_a, \ \rangle $ (sum on repeated indexes) for Schwartz functions $\textbf{e}_a$ and $\textbf{e}^*_a$ and, for $n\in \N$ with $n-1\ge s+\textbf{d}$,
 and by
\begin{equation}\label{eq:reg4}  \begin{aligned} &  \|  P_{N  _g ( p_0)}
		J\Diamond  _i ( 1- \langle \epsilon \Diamond\rangle ^{-2})   \| _{{B(\Sigma _{-r},\Sigma _{r}) }} \le
		\|  \textbf{e}_a \langle J\Diamond  _i ( 1- \langle \epsilon \Diamond\rangle ^{-2})\textbf{e}^*_a, \ \rangle
		   \| _{{B(\Sigma _{-r},\Sigma _{r}) }}
		\\& \le  \| \textbf{e}_a\| _{\Sigma _r}
		\|     ( 1- \langle \epsilon \Diamond\rangle ^{-2}) \textbf{e}_a^*\| _{\Sigma _{r+\textbf{d}}}
		\le C (\epsilon) \| \textbf{e}_a\| _{\Sigma _r} \    	\|      \textbf{e}_a^*\| _{\Sigma _{r'}}
\end{aligned}   \end{equation}
   $C (\epsilon)= \|
		 \Diamond   ( 1- \langle \epsilon \Diamond\rangle ^{-2})   \| _{B(\Sigma _{r'},\Sigma _{r+\textbf{d}}) } $ is bounded  by \eqref{eq:a71}  for   $|\epsilon |\le 1$ for any pair $(r',r)$ with $r'>r+\textbf{d}$.

\noindent We consider now Claim (2). We have

\begin{equation*}  \begin{aligned}& {\mathcal X}^t -{\mathcal X}^t_\epsilon=  \mathcal{A }_j(t,\varrho ,R)\left (  J ( 1- \langle \epsilon \Diamond\rangle ^{-2})\Diamond _j R  -  P_{N  _g ( p_0)}   J\Diamond _j ( 1- \langle \epsilon \Diamond\rangle ^{-2})R\right ).
 \end{aligned}
\end{equation*}
We have $P_{N  _g ( p_0)}   J\Diamond _j ( 1- \langle \epsilon \Diamond\rangle ^{-2})R\stackrel{\epsilon\to 0}{\rightarrow} 0$ for $R\in \Sigma_{n'}$ for any $n'\in \Z$
 because in fact
$C (\epsilon)   \stackrel{\epsilon\to 0}{\rightarrow} 0$   by \eqref{eq:a72}, with $C (\epsilon)$ defined like above for any pair  $( r',r)$ with $r'>r+\textbf{d}$.

\noindent Still by \eqref{eq:a72}, for  $  n > n'+\textbf{d}$ and for $R\in \Sigma_{n'}$  we have by (A5)

\begin{equation} \label{eq:reg5} \begin{aligned} &  \|
		J\Diamond   ( 1- \langle \epsilon \Diamond\rangle ^{-2}) R  \| _{ \Sigma _n  } \le
\|
		 \Diamond   ( 1- \langle \epsilon \Diamond\rangle ^{-2})   \| _{B(\Sigma _{n'},\Sigma _n ) } \| R  \| _{\Sigma_{n'}  } \\& \le C   \|
		  ( 1- \langle \epsilon \Diamond\rangle ^{-2})   \| _{B(\Sigma _{n'} ,\Sigma _{n +\textbf{d}}) } \| R  \| _{\Sigma_{n'}  } \stackrel{\epsilon\to 0}{\rightarrow} 0
.
\end{aligned}   \end{equation}
These facts yield   \eqref{eq:reg3}.

 \noindent We turn now to Claim (3) and to \eqref{eq:reg2}.  By the Rellich criterion,   the embedding
$\Sigma _a\hookrightarrow\Sigma _b$ for $a>b$ is compact. Hence also $\Ph ^a\hookrightarrow\Ph ^b$
 is compact.
Then \eqref{eq:reg2} follows by the  Ascoli--Arzela Theorem by a standard argument.

\qed

\begin{corollary}\label{cor:darboux} Consider \eqref{eq:ODE} defined by
the field
${\mathcal X}^t$ and consider indexes and notation
of Lemma \ref{lem:ODE} (in particular we have $M_0=1$ and $i=1$ in
\eqref{eq:symbol} and elsewhere; $r$ and $M$ can be arbitrary).
Consider $s'$,$s$ and $k$ as in \ref{lem:ODE}. Then  for the map $\mathfrak{F}  ^t \in C^l(  \U _ {\varepsilon _1,k}^{s^{\prime  }}
  , \widetilde{\Ph}  ^{s })$
derived
from \eqref{eq:reg1},
 we have $ \mathfrak{F}  ^{1*}\Omega =\Omega _0$.
\end{corollary}
\proof   $\Omega _0$ is constant  in the coordinate system $(\tau , \Pi , R)$ where $R\in N^\perp _g({\mathcal H}_{p_0}^*),$
  with $\Omega _0=d\tau _j\wedge d\Pi _j +\langle J^{-1} \  ,  \ \rangle $,
where we apply $\langle J^{-1} \  ,  \ \rangle$ only to vectors in the $R$ space. Hence $\Omega _0$ is $C^\infty$
 in $R\in L^2$, $\tau $ and $\Pi$,  with values in
   $B^2 (L^2, \R ) . $ From Lemma \ref{lem:dalpha2}  we have that $d\alpha$, so also $\Omega $ by
	$\Omega =\Omega _0+d\alpha $,  belongs to $C^\infty ( \U _ {\varepsilon _0,k}^{s}, B^2 (\widetilde{ \Ph}, \R ))$  for an $\epsilon _0>0$, and so  also to $C^\infty ( \U _ {\varepsilon _0,k}^{s }, B^2 (\widetilde{ \Ph}^s, \R ))$. Let now $r-(l+1)\textbf{d}\ge s'\ge  s+l\mathbf{d}$ and $k\in \Z \cap [0,r-(l+1)\textbf{d}]$. Then  for a fixed $0<\varepsilon  _2\ll \varepsilon _1$
   and for all $|\epsilon |\le 1$  we have
 \begin{equation} \label{eq:darboux2} \begin{aligned} & \mathfrak{{F}}_{\epsilon}  ^t \in C^l((-2,2)\times \U _ {\varepsilon _2,k}^{s^{\prime  }}
  ,    \U _ {\varepsilon _1,k}^{s }
 )
, \quad  \mathfrak{F}_{\epsilon}  ^t( \U _ {\varepsilon _2,k}^{s^{\prime }}  )\subset
 \U _ {\varepsilon _1,k}^{s^{ \prime}}   \text{ for all $|t |\le 2$}
\end{aligned}   \end{equation}
  by Lemma \ref{lem:ODE}, for a fixed $l\ge 2$.  By Lemma \ref{lem:moll}
  we have uniformly locally
   \begin{equation} \label{eq:darboux3} \begin{aligned} & \lim _{\varepsilon \to 0}\mathfrak{F}  ^t_\epsilon  = \mathfrak{F}  ^t \text{  in } C^{l }([-1,1]\times    \U _ {\varepsilon _2,k}^{s^{\prime }}
  ,  \U _ {\varepsilon _1,k}^{s }
 )  .
\end{aligned}   \end{equation}

  \noindent Let us take  $0<\varepsilon  _3\ll \varepsilon _2$ s.t.
   $\mathfrak{F}_{\epsilon}  ^t(\U _ {\varepsilon _3,k}^{s^{\prime }}   )\subset  \U _ {\varepsilon _2,k}^{s^{\prime }}
$ for all $|t |\le  2$ and   $|\epsilon |\le 1$.

\noindent In $\U _ {\varepsilon _3,k}^{s^{\prime }}  $ the following computation
is valid because $\mathcal{X}_\epsilon ^t$ is a standard vector field in
$\U _ {\varepsilon _1,k}^{s^{\prime }}  $ and similarly $\Omega _t$ is a regular
differential form therein:

\begin{equation*}   \label{eq:darboux1}   \begin{aligned} & \mathfrak{F}_{\epsilon  }^{1*}\Omega  -\Omega _0
=\int _0^1
\frac{d}{dt}
\left ( \mathfrak{F}_{\epsilon  }^{t*}\Omega _t\right )    dt = \int _0^1 \mathfrak{F}_{\epsilon }^{t*}
\left (L_{\mathcal{X}_\epsilon ^t} \Omega _t+\frac{d}{dt}\Omega _t\right )   dt\\&= d
\int _0^1 \mathfrak{F}_{\epsilon }^{t*}
\left (   i_{\mathcal{X}_\epsilon ^t} \Omega  _t+ \alpha \right )    dt
 ,
\end{aligned}
\end{equation*}
where we recall $\Omega  _t=\Omega  _0+t(\Omega -\Omega  _0)$.

\noindent If  we consider a ball
 $\mathbf{B}$ in  $\U _ {\varepsilon _3,k}^{s^{\prime }}   $,
   in the notation of Lemma \ref{lem:1forms},
for some function $\psi _{\epsilon}\in C^1(\mathbf{B},\R )$ we can   write
\begin{equation} \label{eq:darboux11}  \begin{aligned} &   \mathfrak{F}_{\epsilon  }^{1*}(B_0+\alpha )  -B_0 +d\psi _{\epsilon}
=
\int _0^1 \mathfrak{F}_{\epsilon }^{t*}
\left (   i_{\mathcal{X}_\epsilon ^t} \Omega  _t+ \alpha \right )    dt
 ,
\end{aligned}
\end{equation}

 \noindent   By  \eqref{eq:darboux2}--\eqref{eq:darboux3}
we have
\begin{equation*}
 \lim _{\epsilon \to 0}(\mathfrak{F}_{\epsilon  }^{1*}(B_0+\alpha )  -B_0)=\mathfrak{F} ^{1*}(B_0+\alpha )  -B_0 \text{ in } C^{l -1}(    \U _ {\varepsilon _3,k}^{s^{\prime }}
  , B  (\widetilde{ \Ph}^{s' }, \R )
 ) .
\end{equation*}
The set   $\Gamma :=\{\mathfrak{F}_{\epsilon  }^{t}(\mathbf{B}): |t|\le 2, |\epsilon |\le 1   \}$ is a bounded subset in  $\U _ {\varepsilon _2,k}^{s^{\prime }}  $ because of \eqref{eq:gronwall0}--\eqref{eq:gronwall1}. Then we have
\begin{equation*}  \begin{aligned}& \lim _{\epsilon \to 0}{\mathcal X}^t_\epsilon ={\mathcal X} ^t  \text{  in } C^{0}((-2,2)\times  \Gamma
  ,  \widetilde{{\Ph}} ^{s }
 )  \text{   uniformly }.
   \end{aligned}
\end{equation*}
Hence by $i_{\mathcal{X}  ^t} \Omega  _t=-\alpha$ we get
\begin{equation*}
  \lim _{\epsilon \to 0}\left (   i_{\mathcal{X}_\epsilon ^t} \Omega  _t+ \alpha \right )=     i_{\mathcal{X}  ^t} \Omega  _t+ \alpha =0 \text{  in } C^{0}((-2,2)\times  \Gamma
  , B( \widetilde{{\Ph}} ^{s },\R )
 )  \text{   uniformly. }
\end{equation*}
This implies  \begin{equation*}  \begin{aligned}& \lim _{\epsilon \to 0}
\|
\int _0^1 \mathfrak{F}_{\epsilon }^{t*}
\left (   i_{\mathcal{X}_\epsilon ^t} \Omega  _t+ \alpha \right )    dt \| _{L^\infty (\mathbf{B} ,B ( \widetilde{{\Ph}} ^{s' },\R ))}\\& \le C \lim _{\epsilon \to 0}\| i_{\mathcal{X}_\epsilon ^t} \Omega  _t+ \alpha \| _{L^\infty ([0,1]\times \Gamma ,B ( \widetilde{{\Ph}} ^{s },\R ))}=0,
 \end{aligned}
\end{equation*}
for $C$ an upper bound to the  norms  $\| (\mathfrak{F}_{\epsilon  }^{t*}) _{|\mathfrak{F}_{\epsilon  }^{t }(\upsilon)} :B ( \widetilde{{\Ph}} ^{s '},\R )\to B ( \widetilde{{\Ph}} ^{s  },\R ) \|  $  as $\upsilon $ varies in  $  \mathbf{B}$. Notice that $C<\infty$
by \eqref{eq:reg2}.

\noindent  By \eqref{eq:darboux11} we conclude that uniformly
\begin{equation*}
 \lim _{\epsilon \to 0} d\psi _{\epsilon} =B_0-\mathfrak{F} ^{1*}(B_0+\alpha )\text{ in } C^{0}(   \mathbf{B}
  , B  (\widetilde{ \Ph}^{s'}, \R )).
\end{equation*}
Normalizing $\psi _{\epsilon} (\upsilon _0)=0$ at some given $\upsilon _0\in \mathbf{B}$, it follows that
also $\psi _{\epsilon}$ converges locally uniformly to a function
$\psi _{0}$ with $d\psi _{0}= B_0-\mathfrak{F} ^{1*}(B_0+\alpha )$.
Taking the exterior differential, we conclude that
$
   \mathfrak{F}  ^{1*}\Omega  =\Omega _0$  { in } $C^{\infty}( \U _ {\varepsilon _3,k}^{s^{\prime }}
  , B^2 (\widetilde{ \Ph}^{s'}, \R )
 ) .$

\qed

\section{Pullback of the Hamiltonian}
\label{sec:pullback}

In  the somewhat abstract set up of this paper it is particularly important to have a general description
of the pullbacks of the Hamiltonian  $K$.  Our main  goal in this section is formula  \eqref{eq:back11}. This formula and its  related expansion
in Lemma \ref{lem:ExpH11} obtained splitting  $R$ in discrete and continuous modes, play a key role  in the   Birkhoff normal forms argument.

The first and quite general result is the following consequence of
  Lemma \ref{lem:ODE}.
  \begin{lemma} \label{lem:ODE1}   Consider
  $\mathfrak{F} =\mathfrak{F}_1 \circ \cdots \circ \mathfrak{F}_L$
  with  $   \mathfrak{F}_j= \mathfrak{F}_j ^{t=1}$ transformations as of Lemma \ref{lem:ODE}. Suppose that for $j$
      we have $M_0=m_j$, with given numbers  $1\le m_1\le ...\le m_L$.
      Suppose also that all the $j$ we have the same pair   $ r$ and $ M $, which we assume sufficiently large. Let $i_j=1$ if $m_j=1$. Fix $0<m'<M$

      \begin{itemize}
\item[(1)]
      Let  $ r >  2L(m'+1)  \mathbf{d} + s'_L > 4L(m'+1) \mathbf{d} + s_1$, $s_1\ge \textbf{d}$.
	Then,  for any $\varepsilon >0$ there exists a $\delta  >0$
such that $\mathfrak{F}\in C^{m'} (\U  ^{s'_L}  _{\delta  ,a}, \U ^{s _1}  _{\varepsilon ,h})$    for $0\le a\le h $ and
$0\le h < r-(m'+1) \mathbf{d}$.

\item[(2)]
      Let    $ r >    2L(m'+1)  )  \mathbf{d} +h>  4L(m'+1)   \mathbf{d} + a$, $a\ge 0$.
			The above composition, interpreting the $\mathfrak{F}_j$'s as maps in the $(\varrho , R)$
				variables  as in Lemma \ref{lem:ODEdomains},
				yields also  $\mathfrak{F}    \in C^{m'}(  \U _{-a}
  , \Ph ^{-h }
 )
$ for  $\U_  {-a} $ a sufficiently small neighborhood of the origin in $  \Ph ^{-a} $.

\item[(3)] For   $\U_  {-a} \subset  \Ph ^{-a}  $ like above and   for functions
$\mathcal{R}^{i, j}  _{a,m'}\in C^{m'} (\U_  {-a}  ,\R ) $ and $\mathbf{S}^{i, j}  _{a,m'}\in C^{m'} (\U_  {-a} ,\Sigma _{ a}) $, the following formulas hold:
\begin{equation} \label{eq:ODE2}\begin{aligned} & \Pi (R'):= \Pi (R)\circ  \mathfrak{F}  = \Pi (R)
 +\mathcal{R}^{i_1, m_1+1}  _{a,m'}(\Pi (R), R)  , \\& p':=p\circ  \mathfrak{F}    = p +\mathcal{R}^{i_1, m_1+1}  _{a,m'}(\Pi (R), R),\\& \Phi  _{p'} =   \Phi  _{p}  +\mathbf{S}^{i_1, m_1+1}  _{a,m'}(\Pi (R), R) .
\end{aligned}   \end{equation}
\item[(4)]  For a   function $F$ such that $F(e^{J\tau \cdot \Diamond }U)\equiv F(U)$    we have

\end{itemize}
\begin{equation*} \label{eq:ODE4}\begin{aligned} &
 F\circ \mathfrak{F}(U)= F\left ( \Phi  _{p}+P(p)  (R+  \mathbf{S}^{i_1,m_1 }_{k'   ,m '}) + \textbf{S}^{i_1,m_1 +1 }_{k' ,m'} \right ) \, , \,  k'=r-  7L (m'+1) \mathbf{d} .
\end{aligned}   \end{equation*}
\end{lemma}

\proof
Recall that by \eqref{eq:reg1} we have $\mathfrak{F}_j\in C^{m'}(\U _{\varepsilon _j ',h}^{s'_j}, \U _{\varepsilon  _j ,h}^{s_j}  )$ for $r-(m'+1)\textbf{d}> s'_j\ge s_j+m'\textbf{d}$ and appropriate choice  of the  $0<\varepsilon _j  '<\varepsilon _j   $ and for $h\in \Z \cap [0,r-(m'+1)\textbf{d}]$.
So for the composition we  have   $\mathfrak{F} \in C^{m'}(\U _{\varepsilon _L ' ,a  }^{\kappa}, \U _{\varepsilon  _1,h }^{s_1}  )$ for $a\le h$.
 The inequalities   $ r >  2L(m'+1)  \mathbf{d} + s'_L > 4L(m'+1) \mathbf{d} + s_1$, $s_1\ge\textbf{ d} $ can be accommodated
 since $r$ is assumed sufficiently large. This yields claim (1).

\noindent By Lemma \ref{lem:ODEdomains} we have    $\mathfrak{F}_j\in C^{m'}(\U_  {- h +j m'\textbf{d} } , \Ph ^{- h +(j-1)m'\textbf{d} }   )$  with $\U_  {- h +j m'\textbf{d} }\subset \Ph ^{-h+j m'\textbf{d} } $ a neighborhood of the origin. So for the composition  we have
 $\mathfrak{F}    \in C^{m'}(  \U _{-a  }
  , \Ph ^{-h }
 )
$ for $a\le h-Lm'\textbf{d}$. The conditions $ r >  2L(m'+1)  \mathbf{d} +h$,  $h>   4L(m'+1)    \mathbf{d} + a$ and $a\ge 0$,  can be accommodated
 since $r$ is assumed sufficiently large. This yields claim (2).

\noindent  We now prove \eqref{eq:ODE2}.
Let first  $L=1$.
By \eqref{eq:ODE1} we have $R':=( \mathfrak{F}_1)_R(\Pi (R),R)=e^{Jq _1\cdot \Diamond } ( R+ \textbf{S}^{i_1,m_1  }_{r-(m'+1)\mathbf{d},m'}  )$, where  we use  $M> m'$. Here we will omit the variables $(\Pi (R), R)$ in the $\mathbf{S}$'s
and $\resto$'s.
Then we have  for $a'=r-(m'+1)\mathbf{d}$
\begin{equation}\label{eq:pi1}
 \Pi (R') = \Pi ( R+\textbf{S}^{i_1,m_1  }_{a',m'}   ) = \Pi (R)
+\mathcal{R}^{i_1,m_1+1}_{a'-\mathbf{d},m'}.
\end{equation}
Here we have used
\begin{equation*}
 \begin{aligned} |\langle R, \Diamond \textbf{S}^{i_1,m_1  }_{a' ,m'}\rangle | \le  \| R \| _{\Sigma _{-a'+\mathbf{d}}}  \| \textbf{S}^{i_1,m_1  }_{a' ,m'} \| _{\Sigma _{a'}} . \end{aligned}
\end{equation*}
By $ p_j=\Pi _j-\Pi _j(R)+ \mathcal{R}^{0, 2} (\Pi (R) , R)
   $  we get
   \begin{equation}\label{eq:pi2}\begin{aligned} &
   p_j'=\Pi _j-\Pi _j(R')+ \mathcal{R}^{0, 2} (\Pi (R') , R')\\& = \Pi _j-\Pi _j(R )+ \mathcal{R}^{0, 2} (\Pi (R ) , R ) + \mathcal{R}^{i_1,m_1+1}_{a'-\mathbf{d} ,m'}= {p_j}+ \mathcal{R}^{i_1,m_1+1}_{a'-\mathbf{d} ,m'}  .
\end{aligned}
\end{equation}
   This yields  \eqref{eq:ODE2} for $L=1$ since $a\le r-4(m'+1)\textbf{d}<a'-\mathbf{d}$.
	We extend the proof    to the case $L>1$.
	We write here and below $\mathfrak{F} ':= \mathfrak{F}_1 \circ \cdots \circ \mathfrak{F} _{L-1} $.
	We suppose  that $ \mathfrak{F}'_R (\Pi (R),R)=e^{Jq  \cdot \Diamond } ( R+ \textbf{S}^{i_1,m_1  }_{a'_{L-1} ,m'}  )$  for $a'_{L-1}\le r-2(L-1) m' \textbf{d}$, which is true for $L-1=1$.   Then
	\begin{equation*}\begin{aligned} &
 R'=e^{J(q \circ\mathfrak{F} _{L} )  \cdot \Diamond } \left ( e^{Jq_L  \cdot \Diamond } ( R+ \textbf{S}^{i_L,m_L  }_{r-(m'+1)\mathbf{d},m'}  )+ \textbf{S}^{i_1,m_1  }_{a'_{L-1} ,m'}\circ\mathfrak{F} _{L}  \right ) \\&  = e^{J ({q\circ\mathfrak{F} _{L}+q_L} )  \cdot \Diamond }  \left  ( R+ \textbf{S}^{i_L,m_L  }_{r-(m'+1)\mathbf{d},m'}  )+e^{-Jq_L  \cdot \Diamond }\textbf{S}^{i_1,m_1  }_{a'_{L-1}-  m' \mathbf{d},m'}  \right ),
\end{aligned}
\end{equation*}
where $q_L=\resto ^{0,m_L+1  }_{r-(m'+1)\mathbf{d},m'}$ and where we used the last claim in Lemma \ref{lem:ODEdomains}. Since   $e^{-Jq_L  \cdot \Diamond }\textbf{S}^{i_1,m_1  }_{a'_{L-1}- m' \mathbf{d},m'} =\textbf{S}^{i_1,m_1  }_{a'_{L-1}-2m'\mathbf{d},m'}  $  we conclude that there is an expansion
 $ R'=e^{Jq  \cdot \Diamond } ( R+ \textbf{S}^{i_1,m_1  }_{a'_{L } ,m'}  )$
for $a'_{L }\le a'_{L-1}-2m' \mathbf{d}$. Then
\begin{equation}\label{eq:tranR}\begin{aligned} & \mathfrak{F}_R (\Pi (R),R)=e^{Jq  \cdot \Diamond } ( R+ \textbf{S}^{i_1,m_1  }_{a'_{L } ,m'}  )    \, , \quad  a'_{L }:=r- 2L m' \textbf{d}.\end{aligned}
\end{equation}
For $a'=a'_{L }$ formulas \eqref{eq:pi1}--\eqref{eq:pi2} continue to hold.
 By
$a< a'_{L }-\mathbf{d}  $ this yields \eqref{eq:ODE2}.

   \noindent  We consider the last statement of Lemma \ref{eq:ODE1}. For $a'=r-(m'+1)\mathbf{d}$
	we have
\begin{equation*}  \begin{aligned} &
 F(\mathfrak{F}_1(U))= F(\Phi  _{p'}+P(p') e^{Jq _1\cdot \Diamond } ( R+ \textbf{S}^{i_1,m_1  }_{a' ,m'}  ))=\\&   F(\Phi  _{p}+P(p)e^{Jq _1\cdot \Diamond } ( R+ \textbf{S}^{i_1,m_1  }_{a' ,m'}  )+\mathbf{S}^{i_1,m_1+1 } _{a'+\mathbf{d} ,m'}) =\\&   F\left (e^{Jq _1 \cdot \Diamond } \left (\Phi  _{p}+P(p)  (R+ \mathbf{S}^{i_1,m_1 }_{a' ,m'}) +Y\right ) \right )
\end{aligned}   \end{equation*}
with
 \begin{equation*}   \begin{aligned} &Y =(e^{Jq _1 \cdot \Diamond }-1)  \Phi  _{p}+  [P(p),e^{Jq _1 \cdot \Diamond }]  (R+ \textbf{S}^{i_1,m_1  }_{a' ,m'})  +e^{-Jq _1 \cdot \Diamond } \textbf{S}^{i_1,m_1 +1 }_{a'-\mathbf{d},m'}   .
\end{aligned}   \end{equation*}
We claim
\begin{equation}\label{eq:Ypull}  \begin{aligned} &Y = \textbf{S}^{i_1,m_1 +1 }_{a'-2m'\mathbf{d},m'} .
\end{aligned}   \end{equation}
To prove \eqref{eq:Ypull}
 we use  $(e^{Jq _1 \cdot \Diamond }-1)  \Phi  _{p}  = \mathbf{S}^{i_1,m_1 +1} _{r-(m'+1)\mathbf{d},m'}  =\mathbf{S}^{i_1,m_1 +1} _{a' ,m'} $. This  follows from  $\Phi  _{p}\in C^\infty (\mathcal{O},\mathcal{S})$ and
\begin{equation} \label{eq:reg} \begin{aligned} & \left | (e^{Jq _1 \cdot \Diamond }-1)  \Phi  _{p} \right | _{\Sigma _l}
\le  |q_{1j}| \int _0^1  \left |  e^{t  Jq _1 \cdot \Diamond } \Diamond  _j  \Phi  _{p} \right | _{\Sigma _l} dt\le C_l  |q_{1j}|    \left |  \Diamond  _j  \Phi  _{p} \right | _{\Sigma _l}.
\end{aligned}   \end{equation}

\noindent Schematically  we have, summing over repeated indexes and for $ \mathbf{e}_j,\mathbf{e}_j^*
\in \mathcal{S}$,
\begin{equation*}  \begin{aligned} &[P(p),e^{Jq  _1\cdot \Diamond }] = [e^{Jq _1 \cdot \Diamond }, P_{N_{g}}(p) ] = e^{Jq _1 \cdot \Diamond }\mathbf{e}_j
\langle \mathbf{e}_j ^*, \, \rangle -\mathbf{e}_j
\langle e^{-Jq  _1\cdot \Diamond }\mathbf{e}_j ^*, \, \rangle \\&  =(e^{Jq _1 \cdot \Diamond }-1) \mathbf{e}_j
\langle \mathbf{e}_j ^*, \, \rangle -\mathbf{e}_j
\langle (e^{-Jq _1 \cdot \Diamond }-1)\mathbf{e}_j ^*, \, \rangle \\& =\mathbf{S}^{0,m_1+1 }_{r-(m'+1)\mathbf{d},m'} \langle \mathbf{e}_j ^*, \, \rangle + \mathbf{e}_j\langle \mathbf{S}^{0,m_1+1 }_{r-(m'+1)\mathbf{d},m'}, \, \rangle .
\end{aligned}   \end{equation*}
This yields  for any $a^{\prime\prime} \le a'=r-(m'+1)\mathbf{d}$
\begin{equation*}    [P(p),e^{Jq_1\cdot \Diamond }] (R+ \textbf{S}^{i_1,m_1  }_{a^{\prime\prime},m'})=\textbf{S}^{i_1,m_1 +2 }_{a^{\prime\prime },m'}.  \end{equation*}
We have   $e^{-Jq _1 \cdot \Diamond } \textbf{S}^{i_1,m_1+1  }_{a'-\mathbf{d},m'}  =\textbf{S}^{i_1,m_1+1  }_{a'-(m' +1)\mathbf{d},
m'}$. Then  \eqref{eq:Ypull} is proved.
Then
 \begin{equation}\label{eq:pi3}  \begin{aligned} & F(\mathfrak{F}_1(U))=
 F\left ( \Phi  _{p}+P(p)  (R+ \mathbf{S}^{i_1,m_1 }_{a'-2m'\mathbf{d},m'}) + \textbf{S}^{i_1,m_1 +1 }_{a'-2m'\mathbf{d},m'}\right )
\end{aligned}   \end{equation}
for $a'=r-(m'+1)\mathbf{d}$. This proves the last sentence of our lemma for $L=1$.
For $L>1$ set once more $\mathfrak{F} ':= \mathfrak{F}_1 \circ \cdots \circ \mathfrak{F} _{L-1} $.
We assume by induction that $ F(\mathfrak{F}'(U))$ equals the rhs of
\eqref{eq:pi3}  for $a'=a'_{L-1}:=r-2(L-1) m' \textbf{d}$.
 Then using  $\mathbf{S}^{i_1,m_1 }_{l,m '}\circ \mathfrak{F}_L= \mathbf{S}^{i_1,m_1 }_{l- m' \mathbf{d},m '}$ from Lemma
 \ref{lem:ODEdomains}, by \eqref{eq:ODE2} for $\mathfrak{F}=\mathfrak{F}_L$ and by \eqref{eq:Ypull} with the index 1 replaced by index $L$, we get

\begin{equation*}  \begin{aligned} &
 F(\mathfrak{F}(U))= F\big (\Phi  _{p'}+\\& P(p') e^{Jq _L\cdot \Diamond } ( R+ \textbf{S}^{i_L,m_L  }_{r-(m'+1)\mathbf{d},m'}  )  +P(p')    \mathbf{S}^{i_1,m_1 }_{a'_{L-1}- m' \mathbf{d},m '}+ \textbf{S}^{i_1,m_1 +1 }_{a'_{L-1}- m' \mathbf{d},m'} \, \big ) \\& =F\left (e^{Jq_L \cdot \Diamond } \left [\Phi  _{p}+P(p)  (R+  \mathbf{S}^{i_1,m_1 }_{a'_{L-1}- m' \mathbf{d},m '}) + \textbf{S}^{i_1,m_1 +1 }_{a'_{L-1}-2m' \mathbf{d},m'} \right ] \right ) .
\end{aligned}   \end{equation*}

\noindent We conclude  that $ F(\mathfrak{F}(U))$ equals the rhs of
\eqref{eq:pi3}  for $a'_{L}=r-2L m' \textbf{d}$. In particular this proves the last sentence of our lemma for any $L $.

\qed

\begin{lemma}
  \label{lem:Taylor ex}
For fixed vectors  $\mathbf{u}$ and $\mathbf{v} $ and for  $B$ sufficiently regular with   $B(0)=0$,
we have
\begin{equation}\label{ExpEP0}\begin{aligned}
 &   B (|\mathbf{u} +\mathbf{v}  |^2_1)  =B \left (   |\mathbf{u} |_1^2
   \right )+ B (| \mathbf{v}  |_1^2)\\& + \sum _{ j =0}^{3}\int _{[0,1]^2} \frac{t^j}{j!}   (\partial _t ^{j+1})_{|t=0} \partial _s   [B (|s\mathbf{u} +t\mathbf{v} |^2_1) ]\  dt ds  \\&  +\int _{[0,1]^2}  dt ds \int _0^t  \partial _\tau ^5 \partial _s [B (|s\mathbf{u} +\tau \mathbf{v}  |_1^2)] \frac{(t-\tau )^3}{3!} \ d\tau  .
\end{aligned}\end{equation}
\end{lemma} \proof Follows by Taylor expansion in $t$ of

\begin{equation*}\begin{aligned}
 &   B (|\mathbf{u} +\mathbf{v}  |_1^2)  =B \left (   |\mathbf{u} |_1^2
   \right ) +\int _0^1 \partial _t[  B (| \mathbf{u} +t\mathbf{v} |^2_1)]   dt =\\&
	B \left (   |\mathbf{u} |^2_1
   \right ) + B (| \mathbf{v}  |^2_1)+\int _{[0,1]^2}  dt ds \ \partial _s \partial _t[  B (|s \mathbf{u} +t\mathbf{v} |_1^2)]      .
\end{aligned}\end{equation*} \qed

\begin{lemma}
  \label{lem:back}  Consider a  transformation
  $\mathfrak{F} =\mathfrak{F}_1 \circ \cdots \circ \mathfrak{F}_L$
  like in Lemma \ref{lem:ODE1}  and with $m_1=1$,   with same notations, hypotheses and conclusions.
	In particular we suppose $r$ and $M$ sufficiently large that the conclusions of  Lemma \ref{lem:ODE1}
	hold for preassigned sufficiently large $s=s'_L$, $k' $ and $m' $.
	Let  $k \le k' - \max\{\textbf{d}, \text{ord}(\mathcal{D} )\} $ and $m\le m'$.
Then there are a  $ \underline{{{\psi}}} (\varrho) \in C^\infty $  with $\underline{{\psi}} (\varrho) =O(| \varrho |^2)$  near 0
		and a   small  $\varepsilon >0$
 such that in
			$\U^{s}_{\varepsilon ,k}$   we have the expansion

\begin{align}  \label{eq:back1}   &
K \circ \mathfrak{F}=  \underline{{{\psi}}} (\Pi  (R)) +   \frac  {1}{2}\Omega (  \mathcal{H}_pP(p)R,P(p) R )   +  \resto ^{1,2} _{k,m}   +E_P(P(p)R)+\textbf{R} ^{\prime \prime}
\\&   \textbf{R} ^{\prime \prime }:=     \sum _{d=2}^4
\langle B_{d } (   R,\Pi (R) ),  (P(p)R)  ^{   d} \rangle
      +\int _{\mathbb{R}^3}
B_5 (x,  R, R(x),\Pi (R) )   (P(p)R)^{   5}(x) dx \nonumber
\end{align}
with:

  \begin{itemize}
\item $ \resto ^{1,2} _{k,m}  =  \resto ^{1,2} _{k,m}   (\Pi (R), R) $;

\item   $B_2(0,0  )=0$;

\item $(P(p)R)^d(x)$  represent $d-$products of components of $P(p)R$;

  \item
$B_{d
}( \cdot ,  R,\varrho  ) \in C^{m } ( \U _{-k},
\Sigma _k (\mathbb{R}^3, B   (
 (\mathbb{R}^{2N  })^{\otimes d},\mathbb{R} ))) $  for $2\le d \le 4$ with $\U _{-k}\subset \Ph ^{-k}$ a neighborhood of the origin;
 \item  for
$ \zeta \in \mathbb{R}^{2N  }$ with $|\zeta  |\le \varepsilon$
 and $( \varrho ,R) \in \U _{-k}$
 we have  for $i\le m$
\begin{equation} \label{eq:B5}\begin{aligned} &  \| \nabla _{ R,\zeta, \varrho  }
^iB_5(  R,\zeta  ,\varrho  ) \| _{\Sigma _k(\mathbb{R}^3,   B   (
 (\mathbb{R}^{2N  })^{\otimes 5},\mathbb{R} )} \le C_i .
 \end{aligned}  \end{equation}
\end{itemize}

\end{lemma}
\proof Here we will omit the variables $(\Pi (R), R)$ in the $\mathbf{S}$'s
and $\resto$'s.

\noindent By Lemma \ref{lem:ODE1}  for  $m\le m' \le M$,  $k+\max\{\textbf{d}, \text{ord}(\mathcal{D} )\}\le k' \le r-L(m'+2)\textbf{d}$,
  we have
 \begin{equation}    \begin{aligned} &
 {K}(\mathfrak{F}(U))= E (\Phi _p+P(p)R+P(p)\mathbf{S}^ {1,1}_{k',m'}+\mathbf{S}^{1,2}_{k',m'})- E\left (   \Phi _{p _0}\right ) \\& -(\lambda _j(p) + \resto ^{1,2} _{k,m})   \left ( \Pi _j (\Phi _p+P(p)R ) + \mathcal{R}^{1,2}_{k,m}  -\Pi _j\left (   \Phi _{p _0}\right ) \right ) ,
\end{aligned}
\end{equation}
where, by  \eqref{eq:ODE2}, we have used  $p':=p\circ \mathfrak{F}= p +\mathcal{R}^{1,2}_{k,m}$ and where by
$k\le   k'  -  \textbf{d} $
\begin{equation*}\begin{aligned} &
\Pi _j (\Phi _p+P(p)R+P(p)\mathbf{S}^{1, 1}_{k',m'}   	 +\mathbf{S}^{1, 2}_{k',m'} )=\Pi _j (\Phi _p+P(p)R ) +\mathcal{R}^{1,2 } _{k,m} .\end{aligned}
\end{equation*}
Set now $\Psi = \Phi _p +P(p)\mathbf{S}^{1,1}  _{k',m'}  	 +\mathbf{S}^{1,2}_{k',m'}  	$.  By \eqref{ExpEP0} for $\textbf{u}=\Psi$ and $\textbf{v}=
P(p)R$

 \begin{equation}   \label{eq:back2}     \begin{aligned} & E_P( \Psi +P(p)R )= E_P( \Psi   ) + E_P( P(p)R )  \\& +\sum _{ j =0}^{1}\int _{\R ^3}dx \int _{[0,1]^2}
 \frac{t^j}{j!}   (\partial _t^{j+1})_{|t=0} \partial _s   [B (|s \Psi+tP(p)R  |^2_1) ]  dt ds  \\&+ \sum _{ j =2}^{3}\int _{\R ^3}dx\int _{[0,1]^2} \frac{t^j}{j!}   (\partial _t^{j+1})_{|t=0} \partial _s [  B (|s \Psi  +tP(p)R  |_1^2) ]  dt ds \\& + \int _{\R ^3}dx \int _{[0,1]^2}  dt ds \int _0^t  \partial _\tau ^5 \partial _s [B (|s \Psi   +\tau P(p)R   |_1^2)] \frac{(t-\tau )^3}{3!} d\tau  .
\end{aligned}
\end{equation}
The last two lines can be incorporated in $ \textbf{R}^{\prime \prime}  	$.  For example, schematically
we have
\begin{equation*}
    \partial _\tau ^5   \partial _s B (|s\Phi _{p } +\tau P(p)R  |_1^2)
     \sim   \widetilde{B} ( s\Phi _{p } +\tau P(p)R) \  \Phi _{p } \ (P(p)R  )^5,
\end{equation*}
  for some $\widetilde{B}(Y)\in C^\infty (\R ^{2N}, B ^6(\R ^{2N}, \R )).$
  This produces a term which can be absorbed in the $B_5$ term of $\mathbf{R} ^{\prime  \prime} .$ In particular, \eqref{eq:B5}
   follows from \eqref{eq:growthB}. The terms in the third line of \eqref{eq:back2} can be treated similarly yielding terms
   which end in the  $B_d$  term of $\mathbf{R} ^{\prime  \prime}  $ with $d=j+1$.

  \noindent The second line of  \eqref{eq:back2}  equals

 \begin{equation}   \label{eq:back3}     \begin{aligned} &  \int _{\R ^3}dx \int _{[0,1]^2} dt ds
\sum _{ j =0}^{1}\frac{t^j}{j!}(\partial _t^{j+1})_{|t=0} \partial _s   \ \big \{ \ \  B (|s \Phi _p+tP(p)R  |^2_1)    +  \\&  +  \int _0^1 d\tau   \partial _\tau  [B (|s (\Phi _p
+\tau (P(p)\mathbf{S}^{1,1}_{k',m'} 	+\mathbf{S}^{1,2}_{k',m'}) +tP(p)R  |^2_1) ]  \  \  \big \}   .
\end{aligned}
\end{equation}
The contribution from the  last line of   \eqref{eq:back3}  can be incorporated in $\textbf{R}^{\prime \prime}+\resto  ^{1,2} _{k,m}  	$.

 \noindent  By  $k \le k' -\text{ord}(\mathcal{D} )$ we have
\begin{equation*}     \label{eq:back31}   \begin{aligned} & E_K( \Psi +P(p)R )
= E_K(  \Psi )  +\langle   {\mathcal D} \Phi _p,      P(p)R \rangle  \\& + \overbrace{\langle   {\mathcal D} (P(p)\mathbf{S}^{1,1}_{k',m'}+\mathbf{S}^{1,2}_{k'   ,m'}),      P(p)R \rangle}^{\resto  ^{1,2}_{k,m}}  	+E_K(  P(p)R ).\end{aligned}
\end{equation*}
Notice that from the
$j=0$ term  in the first line of   \eqref{eq:back3} we get
\begin{equation*} \begin{aligned} & 2\int _{\R^3}dx \int _{0}^1   ds
    \partial _s   [B' (|s \Phi _p  |^2_1) s \Phi _p\cdot _{1} P(p)R  ] =2\int _{\R^3}dx B' (|  \Phi _p  |^2_1)   \Phi _p\cdot _{1} P(p)R \\& =\langle \nabla E_P(\Phi _p), P(p) R\rangle .
\end{aligned}
\end{equation*}
By \eqref{eq:lagr mult} and \eqref{eq:begspectdec2}, that is $\nabla E(\Phi _p  )= \lambda   (p) \cdot \Diamond \Phi _p \in N_g (\mathcal{H}_p^{\ast})$,
  and by $P(p)R\in N_g^\perp (\mathcal{H}_p)$, we    have
\begin{equation*}
\langle   {\mathcal D} \Phi _p,      P(p)R \rangle + \langle \nabla E_P(\Phi _p), P(p) R\rangle =\langle \nabla E(\Phi _p),   P(p)R \rangle  =0.
\end{equation*}
The
$j=1$ term in the first line of   \eqref{eq:back3}  is $\frac{1}{2}
\langle \nabla ^2 E_P(\Phi _p)P(p) R, P(p) R\rangle$ which summed
to the $E_K(  P(p)R )$ in  \eqref{eq:back31} yields the $\frac  {1}{2}\Omega (  \mathcal{H}_pP(p)R,P(p) R )$ in \eqref{eq:back1}.

\noindent We have $  E_K ( \Psi   )+ E_P ( \Psi   )= E ( \Psi   )$ and
\begin{equation*}       \begin{aligned} & E ( \Psi   )
=   E (  \Phi _p )  +\overbrace{\langle \nabla E(\Phi _p),   P(p)\mathbf{S}^{1,1}_{k',m'}\rangle}^{0}+
 \overbrace{\langle \nabla E(\Phi _p),    \mathbf{S}^{1,2 } _{k',m'} \rangle }^{ \resto^{1,2 } _{k,m}}+  \resto^{1,2 } _{k,m}.\end{aligned}
\end{equation*}
The last term we need to analyze, for     $d(p):=E(\Phi _{p  }) -   \lambda  (p) \cdot      \Pi   (\Phi _{p  })$, is

 \begin{equation*}  \label{eq:back4}   \begin{aligned} &
E(\Phi _{p  }) - E(\Phi _{p _0 })-\sum _j\lambda _j(p) (\Pi _j (\Phi _{p  }) - \Pi _j (\Phi _{p _0 })) \\& =
d(p)- d(p_0)- \sum _j (\lambda _j(p_0)-\lambda _j(p))  p  _{0j}  =:  \widetilde{\psi} (p,p_0),
\end{aligned}
\end{equation*}
 where $ \widetilde{\psi} (p,p_0)=O ((p-p_0)^2)$  by   $\partial _{p_j}d(p)=-p\cdot \partial _{p_j} \lambda (p)$.
Notice that $\widetilde{\psi}\in  C^\infty (\mathcal{O}^2, \R) $.
Now recall that in the initial system
of coordinates  we have $p '=\Pi -\Pi (R') +\resto ^{0,2}(\Pi (R') ,R')$.   Substituting
$p'$ and $\Pi (R')$ by means of \eqref{eq:ODE2}, and $R'$ by means of  \eqref{eq:tranR}
  we conclude that $p=p_0 -\Pi (R)+\resto ^{0,2}_{k',m'} 	.$
Then  $ \widetilde{\psi} (p,p_0)=\underline{{{\psi}}} (\Pi  (R))+\resto ^{1,2}_{k,m} $
with
$\underline{{{\psi}}}(\varrho ):= \widetilde{\psi} (p_0-\varrho,p_0)$   a $C^\infty$ function  with
$\underline{{{\psi}}}(\varrho  )=O(|\varrho |^2)$ for $\varrho$ near $0$.

  \qed

\begin{lemma}
  \label{lem:back11}  Under the hypotheses and notation of Lemma \ref{lem:back},
	 for an $\textbf{R}'$ like $\textbf{R}^{\prime \prime}$, for a $\psi \in C^\infty $
	with $\psi (\varrho )=O(|\varrho|^2)$ near 0,   we have
\begin{align}  \label{eq:back11}   &
K \circ \mathfrak{F}=  {{\psi}} (\Pi  (R)) +   \frac  {1}{2}\Omega (  \mathcal{H}_{p _0} R,  R )
 +  \resto ^{1,2} _{k,m}  (\Pi (R),R) 	  +E_P( R)+\textbf{R}',
\\&   \textbf{R} ':=     \sum _{d=2}^4
\langle B_{d } (   R,\Pi (R) ),  R  ^{   d} \rangle
      +\int _{\mathbb{R}^3}
B_5 (x,  R, R(x),\Pi (R) )R^5   ( x) dx, \nonumber
\end{align}
 the $B_d$ for $d=2,...,5$ with similar properties of the functions in Lemma \ref{lem:back}.

 \end{lemma} \proof
We have
\begin{equation*}
 P(p)    R=R+(P(p)-P(p_0)) R=R    +\textbf{S}^{1,1}(p-p_0,R)=R+\textbf{S}^{1,1}(\Pi (R),R).
\end{equation*}
 Substituting $P(p)    R= R+\textbf{S}^{1,1}(\Pi (R),R)$
 in \eqref{eq:back1} we obtain that $ \resto ^{1,2}_{k,m}+\textbf{R}^{\prime \prime} $ is absorbed in
$ \resto ^{1,2}_{k,m}(\Pi (R),R) 	+\textbf{R}^{\prime } $.
This is elementary to see for the terms with $d\le 4$. We consider
the case $d=5$.

 \begin{equation*} \begin{aligned}&
  B_5 (x,  R, R(x),\Pi (R) )    R ^{   i}(x) (\textbf{S}^{1,1}) ^{5-i} \\& =
  \sum _{j=0}^{5-i} \frac{1}{j!}(\partial _t^j)_{|t=0} [B_5 (x,  R, tR(x),\Pi (R) )]  R ^{   i}(x) (\textbf{S}^{1,1}) ^{5-i}\\& + \int _0^1 \frac{(1-t)^{4-i}}{(4-i)!}\partial _t ^{5-i} [B_5 (x,  R, tR(x),\Pi (R) )]  R ^{   i}(x) (\textbf{S}^{1,1}) ^{5-i} \end{aligned}
\end{equation*}
The last term can be absorbed in the $d=5$ term of $\mathbf{R}'$. Similarly, all the other terms either are absorbed  in $\mathbf{R}'$
or, like for instance the $i=j=0$ term, they are   $ \resto ^{1,2}    $.

\noindent We write $E_P(P(p)R)=E_P(R-P_{N_g(p)} R) $ and use \eqref{ExpEP0}
for $\mathbf{u}=R$ and $\mathbf{v}=-P_{N_g(p)} R$.
We get  the sum of
$E_P( R)$ with a term which can be absorbed in $ \resto ^{1,2}_{k,m}(\Pi (R),R) 	 +\textbf{R}^{\prime } $.
We finally focus on
\begin{equation}   \label{eq:back12}   \begin{aligned}& \frac 12    \langle J^{-1}   \mathcal{H}_p P(p)R, P(p)R \rangle  =  \frac 12  \langle    \mathcal{D} P(p)R, P(p)R \rangle \\&- \lambda _j(p) \Pi _j(P(p)R) +  \frac 12    \langle  \nabla ^2 E _P (\Phi _p)P(p)R, P(p)R \rangle     .\end{aligned}
\end{equation}
We have
  \begin{equation}   \label{eq:back13}   \begin{aligned}    \langle    \mathcal{D} P(p)R, P(p)R \rangle
	&= \langle \mathcal{D}  R,  R \rangle  +\resto ^{1,2}_{k,m}(\Pi (R),R) 	 \\
	   \langle  \nabla ^2 E _P (\Phi _p)P(p)R, P(p)R \rangle & =
		\langle  \nabla ^2 E _P (\Phi _{p_0}) R,  R \rangle + \resto ^{1,2}_{k,m}(\Pi (R),R) 	 \\&+ \langle  ( \nabla ^2 E _P (\Phi _p)-\nabla ^2 E _P (\Phi _{p_0}))  R,  R \rangle  \\
	\lambda _j(p) &=\lambda _j(p_0)+  \resto ^{1,0} (\Pi (R)) + \resto ^{1,2}_{k,m}(\Pi (R),R) 	\\  \Pi _j(P(p)R) &= \Pi _j( R ) +\resto ^{1,2}_{k,m}(\Pi (R),R) 		 .\end{aligned} \nonumber
\end{equation}
Then we conclude that the right hand side of   \eqref{eq:back12}  is
\begin{equation}   \label{eq:back120}   \begin{aligned}    &
 \overbrace{\frac 12   \langle (\mathcal{D}  - \lambda  (p_0) \cdot \Diamond  +\nabla ^2 E _P (\Phi _{p_0})) R,  R \rangle }^{\frac 12\langle J^{-1}   \mathcal{H} _{p_0} R,  R \rangle }+\resto ^{2, 0} (\Pi (R) ) +\resto ^{1,2}_{k,m}(\Pi (R),R)
\\& + \frac 12\langle  ( \nabla ^2 E _P (\Phi _p)-\nabla ^2 E _P (\Phi _{p_0}))  R,  R \rangle  \end{aligned}
\end{equation}
where the last  term   can be absorbed in  the $d=2$ term of $\textbf{R}'$.
Setting  $\psi (\varrho ) =\underline{\psi }(\varrho  ) +\resto ^{2, 0} (\varrho  )$ with the $\resto ^{2, 0}$
in \eqref{eq:back120},   we get the desired result.
\qed

\bigskip
We have completed the part of this paper devoted to the Darboux Theorem.  The next step consists in the decomposition of $R$ into discrete and
continuous modes, and the search of a new coordinate system by an appropriate Birkhoff normal forms argument.

\section{Spectral coordinates   associated to $\mathcal{H} _{p_0}$}
\label{sec:speccoo}
We will consider the operator $\mathcal{H} _{p_0}$, which will be central
in our analysis henceforth.  We will list now  various hypotheses, starting
with the spectrum of $\mathcal{H} _{p_0}$ thought as an operator in  the natural complexification $L^2(\R ^3, \C ^{2N}) $ of $L^2(\R ^3, \R ^{2N}) $.

\begin{itemize}
\item[(L1)]    $\sigma _e(\mathcal{H} _{p_0}) $ is a
union of intervals  in  $\im \R $  with $0\not \in \sigma _e(\mathcal{H} _{p_0})$ and is symmetric with respect to 0.

\item[(L2)]    $\sigma _p(\mathcal{H} _{p_0})  $  is finite.

\item[(L3)] For any eigenvalue $\mathbf{e} \in \sigma _p(\mathcal{H} _{p_0})  \backslash \{ 0 \}$  the algebraic and geometric dimensions
    coincide and are finite.

\item [(L4)]  There is  a number   $\mathbf{n}\ge 1$ and positive numbers  $0<\mathbf{e}
_1' \le \mathbf{e} _2'  \le ...\le\mathbf{e} _\mathbf{n} ' $
such that $\sigma _p(\mathcal{H} _{p_0})$ consists exactly of the numbers
$\pm   \im \mathbf{e} _j'$ and 0. We assume that
there are fixed integers $\mathbf{n}_0=0< \mathbf{n}_1<...<\mathbf{n}_{l_0}=\textbf{n}$ such that
$\mathbf{e} _j' = \mathbf{e}  _i '$ exactly for $i$ and $j$
both in $(\mathbf{n}_l, \mathbf{n}_{l+1}]$ for some $l\le l_0$. In this case $\dim
\ker (\mathcal{H}_{p_0}  -\mathbf{e} _j'  )=\mathbf{n}_{l+1}-\mathbf{n}_l$.
 We assume there exist $N_j\in \mathbb{N}$ such that $ N_j+1 =\inf  \{ n\in \N : n \mathbf{e} _j'\in \sigma _e(\mathcal{H} _{p_0}) \}  $. We set $\mathbf{N}=\sup _j N_j$.  We assume that  $\mathbf{e} _j'\not \in \sigma _p(\mathcal{H} _{p_0})$ for all $j$.

\item[(L5)] If $\mathbf{e} _{j_1}'<...<\mathbf{e} _{j_i}'$ are i distinct
  $\lambda$'s, and $\mu\in \Z^k$ satisfies
  $|\mu| \leq 2N +3$, then we have
$$
\mu _1\mathbf{e} _{j_1}'+\dots +\mu _k\mathbf{e} _{j_i}'=0 \iff \mu=0\ .
$$

\end{itemize}
The following hypothesis holds quite generally.
\begin{itemize}
\item[(L6)] If $\varphi \in \ker (\mathcal{H}_{p_0} -\im  \mathbf{e}  )  $ for  $\im \mathbf{e} \in
\sigma _p (\mathcal{H}_{p_0}   )  $ then
$\varphi \in   \mathcal{S}(\R ^3, \C ^{2N})$.
\end{itemize}

\bigskip
\noindent
By \eqref{eq:begspectdec1},     $
\mathcal{H}_{p_0}  \xi   = \mathbf{e}   \xi
 $  implies  $
\mathcal{H}_{p_0} ^* J^{-1} \xi   =-  \mathbf{e}   J^{-1}\xi
 $.   Then $\sigma _p  (\mathcal{H}_{p_0} )=\sigma _p  (\mathcal{H}_{p_0} ^*)$. We denote it
by $\sigma _p$.

\noindent  By general argument we have:
\begin{lemma}
  \label{lem:Specdec} The following spectral decomposition remains determined:
\begin{align}  \label{eq:spectraldecomp} &
N_g^\perp (\mathcal{H}_{p_0} ^*)\otimes _\R  \C =  \big (\oplus _{\mathbf{e} \in
\sigma _p\backslash \{ 0\}}   \ker (\mathcal{H}_{p_0}  - \mathbf{e}
 ) \big) \oplus X_c ( {p_0} )\\& \nonumber  X_c ( {p_0} ):=
\left\{N_g(\mathcal{H}_{p_0} ^\ast)\oplus \big (\oplus _{\mathbf{e} \in
\sigma _p\backslash \{ 0\}}   \ker (\mathcal{H}_{p_0} ^*- \mathbf{e}
 ) \big)\right\} ^\perp .
\end{align}
\end{lemma}

We denote by $P_c$ the projection on $X_c ( {p_0} )$ associated to \eqref{eq:spectraldecomp}.  Set $\mathcal{H}:=\mathcal{H}_{p_0}P_c.$

\noindent The following  hypothesis is    important to solve the homological
equations in the Birkhoff normal forms argument.
\begin{itemize}
\item[(L7)] We have  $R_{\mathcal{H}}\circ \Diamond _j^i  \in C^\omega  ( \rho (\mathcal{H} ), B(\Sigma _{n },\Sigma _{n } )) $ for any $ n\in \N $, any $j=1,...,n_0$ and for any  $i =0,1$,  where $\rho (\mathcal{H} ) =\C \backslash  \sigma _e (\mathcal{H}_{p_0}   )$.

\end{itemize}

\noindent For the examples in Sect. \ref{sec:examples},
(L7) can be checked with standard arguments.

\noindent We   discuss now the choice  of a good frame of   eigenfunctions.

\begin{lemma}
  \label{lem:basis}  It is possible to choose   eigenfunctions   $\xi '\in  \ker (\mathcal{H}_{p_0} -\im  \mathbf{e} _j' )$
so that  $ \Omega ( \xi _j ',\overline{\xi} _k')=0$ for $j\neq k$ and  $ \Omega ( \xi _j' ,\overline{\xi}' _j)=-\im s_j$ with $s_j \in \{ 1 ,-1   \} $ .  We have $ \Omega ( \xi _j ', {\xi} _k')=0$ for all $j$ and $k$.
 We have $ \Omega (\xi , f )=    0$ for any eigenfunction
$\xi  $ and any $f\in X_c ( {p_0} )$.

\end{lemma}
\proof   First of all,  if  $\lambda ,\mu \in \sigma _p(\mathcal{H}_{p_0} )  $ are two eigenvalues with $\lambda
\neq 0$ and given two associated eigenfunctions $\xi _\mu $ and $\xi _\lambda$

 \begin{equation} \label{eq:basis1}
\begin{aligned}    \langle J^{-1}\xi _\lambda , \overline{\xi} _\mu \rangle &=
\frac 1 \lambda \langle J^{-1}\mathcal{H}_{p_0}\xi _\lambda , \overline{\xi} _\mu \rangle
=-\frac 1 \lambda \langle \mathcal{H}_{p_0}^* J^{-1}\xi _\lambda , \overline{\xi} _\mu \rangle
\\& =-\frac 1 \lambda \langle J^{-1}\xi _\lambda , \mathcal{H}_{p_0}  \overline{\xi} _\mu \rangle
=-\frac {\overline{\mu}}  \lambda  \langle J^{-1}\xi _\lambda , \overline{\xi} _\mu \rangle ,
\end{aligned}
\end{equation}
where for the second equality we used  \eqref{eq:begspectdec1}  and for the last one the fact that
 $
\mathcal{H}_{p_0}  \xi   = \mu   \xi
 $  implies    $
\mathcal{H}_{p_0}   \overline{\xi }  = \overline{\mu }   \overline{\xi }   .
 $  Then, for $\mathbf{e} _j\neq \mathbf{e} _k$ and associated eigenfunctions  $\xi _j $ and ${\xi} _k$
we get
$ \Omega ( \xi _j ,\overline{\xi} _k)= 0 $.  Notice that by a similar argument we have $ \Omega ( \xi _\lambda ,  {\xi} _\mu )=-\frac \mu \lambda   \Omega ( \xi _\lambda ,  {\xi} _\mu )$ and so  $ \Omega ( \xi _j ', {\xi} _k')\equiv 0$ .

\noindent   Since $\mathcal{H}_{p_0}  \xi   = \mathbf{e}   \xi
 $  implies  $
\mathcal{H}_{p_0} ^* J^{-1} \xi   =-  \mathbf{e}   J^{-1}\xi
 $, for  any eigenfunction
$\xi  $   of $\mathcal{H}_{p_0}$ then $J^{-1}\xi
 $  is  an  eigenfunction
    of $\mathcal{H}_{p_0}^*$. By the definition of  $X_c ( {p_0} )$ in \eqref{eq:spectraldecomp},
we conclude $ \Omega (\xi , f )=   \langle J^{-1}\xi , f\rangle = 0$
for any $f\in X_c ( {p_0} )$.

\noindent Let  $\im \textbf{e} \in \im \R\backslash \{ 0\}$ be an eigenvalue.  By  the above discussion,   the Hermitian form $\langle \im J^{-1}\xi , \overline{ \eta }\rangle
$ is non degenerate in  $ \ker (\mathcal{H}_{p_0}  - \im \mathbf{e}
 ) $. Then we can find a basis such that  $\langle \im J^{-1}\eta _j , \overline{ \eta }_k\rangle  =- |a_j|\text{sign}(a_j )
\delta _{jk}$, for appropriate non zero numbers $a_j\in \R$. Then set $\xi '=\sqrt{|a _j|}\eta _j.$
 \qed

\bigskip
\noindent We set $\xi _j=\xi '_j$ and  $\textbf{e} _j=\textbf{e} '_j$
if $s_j=1$.

\noindent We set $\overline{\xi} _j=\xi '_j$ and  $\textbf{e} _j=-\textbf{e} '_j$
if $s_j=-1$.

\noindent  Notice that if $f\in  X_c ( {p_0} )$ then also $\overline{f}\in  X_c ( {p_0} )$. This implies that for    $R \in  N_g^\perp (\mathcal{H}_{p_0} ^*)\otimes _\R \C$ with real entries, that is if $R=\overline{R}$, then we have

\begin{equation}
  \label{eq:decomp2}
  R (x) =\sum _{j=1}^{\mathbf{n}}z_j \xi _j (x) +
\sum _{j=1}^{\mathbf{n}}\overline{z}_j\overline{\xi  }_j( x )
+ f (x), \quad   f \in X_c (p_0).
\end{equation}
with $f=\overline{f}$.

\noindent By Lemma \ref{lem:basis}  we have, for the $s_j$ of Lemma \ref{lem:basis},
\begin{equation}
  \label{eq:H2}
 \frac{1}{2} \Omega (\mathcal{H}_{p_0}  R, R ) =  \sum _{j=1}^\mathbf{n}  \mathbf{e} _j |z_j|^2
+\frac{1}{2} \Omega (\mathcal{H}_{p_0}  f, f )=:H_2.
\end{equation}

\noindent Consider the map  $R\to (z,f)$  obtained from \eqref{eq:decomp2}.  In terms of
the pair $(z,f)$, the Fr\'echet derivative $R'$ can be expressed as

\begin{equation*}
    R'=\sum _{j=1}^{\mathbf{n}}(dz_j \xi _j+d\overline{z}_j \overline{\xi} _j)  +f'.
\end{equation*}
We have

\begin{equation}
  \label{eq:OmegaCoo}
   \Omega  (   R' , R ')   =-   \im   \sum _{j=1}^{\mathbf{n}}  dz_j\wedge d \overline{z}_j
+  \Omega (   f' , f ')  .
\end{equation}
For a function $F$ independent of $\tau $ and $\Pi$ let   us decompose $X_F$ as of spectral decomposition \eqref{eq:decomp2}:
\begin{equation*}
    X_F =\sum _{j=1}^{\mathbf{n}}(X_F)_{z_j} \xi _j (x) +
\sum _{j=1}^{\mathbf{n}}(X_F)_{\overline{z}_j}\overline{\xi  }_j( x )
+ (X_F)_{f}, \quad   (X_F)_{f} \in X_c (p_0).
\end{equation*}
By $i_{X_F}\Omega =dF$ and by
\begin{equation*} \begin{aligned} &
 dF= \partial _{z_j}F dz_j+\partial _{\overline{z}_j}F d\overline{z}_j+ \langle \nabla _fF, f'\ \rangle  \\&
i_{X_F}\Omega =- \im  (X_F)_{z_j}  d\overline{z}_j+ \im  (X_F)_{\overline{z}_j}  dz_j+ \langle  J^{-1} (X_F)_{f}, f'\ \rangle ,
\end{aligned}\end{equation*}
we get
\begin{equation*} \begin{aligned} &
(X_F)_{z_j}=\im  \partial _{\overline{z}_j}F\ , \quad (X_F)_{\overline{z}_j}=-\im   \partial _{z_j}F\ , \quad (X_F)_{f}=J\nabla _fF.
\end{aligned}\end{equation*}
This implies
\begin{equation} \label{eq:poiss}\begin{aligned} & \{  F,G  \} :=dF(X_G) = \im \partial _{z_j}F\partial _{\overline{z}_j}G
-\im \partial _{\overline{z}_j}F\partial _{z_j}G + \langle \nabla _fF, J\nabla _fG \rangle.
\end{aligned}\end{equation}
Hence, for $H_2$ defined in \eqref{eq:H2}, for   $z=(z_1,....,z_\mathbf{n})$, using standard multi index notation and by \eqref{eq:begspectdec1},
 we have:
\begin{equation} \label{eq:poiss1}\begin{aligned} & \{  H_2,z^\mu \overline{z} ^\nu  \} =-\im \mathbf{e} \cdot (\mu - \nu ) z^\mu \overline{z}^\nu  \ ; \quad  \{  H_2, \langle J ^{-1}\varphi ,f\rangle   \} =  \langle J ^{-1}\mathcal{H}\varphi ,f\rangle .
\end{aligned}\end{equation}

\subsection{Flows in spectral coordinates}
\label{subsec:flowsspec}

We restate Lemma \ref{lem:ODE} for a special class of transformations.
\begin{lemma}
  \label{lem:chi}
	Consider
\begin{equation}
\label{eq:chi1}\chi   =\sum _{|\mu +\nu |=M_0 +1} b_{\mu\nu} (\Pi (f))  z^{\mu} \overline{z}^{\nu} + \sum _{|\mu +\nu |=M_0 } z^{\mu} \overline{z}^{\nu}
 \langle  J ^{-1}  B_{\mu   \nu
}(\Pi (f))
  , f \rangle
\end{equation}
with $ b_{\mu\nu}(\varrho)= \resto   ^{i,0}  _{r,M}(\varrho)$ and $ B_{\mu\nu}(\varrho)=  \textbf{S}      ^{i,0}  _{r,M}(\varrho)$ with $i\in \{  0,1\}$ fixed and $r,M\in \N$ sufficiently large   and with
\begin{equation}
\label{eq:symm}  \overline{b}_{\mu\nu}  = {b}_{\nu\mu}   \   , \quad   \overline{B}_{\mu   \nu
}  =B_{\nu\mu} ,
\end{equation}
(so that $\chi$ is real valued for $f=\overline{f}$).
Then we have what follows.

\begin{itemize}
\item[(1)]
 Consider the vectorfield $X_\chi  $ defined
with respect to $\Omega _0$.
 Then, summing on repeated indexes (with the equalities defining the   field $X_\chi ^{st}$), we have:
 \begin{equation*}   \begin{aligned} &(X_\chi  ) _{z_j}= \im \partial _{\overline{z}_j}
 \chi =: (X_\chi ^{st} ) _{z_j} \,  , \quad (X_\chi  ) _{\overline{z}_j}= -\im \partial _{z_j}
 \chi  =: (X_\chi ^{st} ) _{\overline{z}_j} \, , \\&  (X_\chi  ) _{f} =\partial  _{\Pi _j(f)}\chi \, P_c^*(p_0)
     J \Diamond _j f +
     (X_\chi ^{st} ) _{f}  \text{ where }  (X_\chi ^{st} ) _{f}:=z^{\mu} \overline{z}^{\nu}B_{\mu   \nu
}(\Pi (f)).
 \end{aligned}
  \end{equation*}

\item[(2)]
Denote by $\phi ^t$ the flow of  $X_\chi$  provided by Lemma \ref{lem:ODE} and
set      $(z^t,f^t)=  (z,f)\circ \phi ^t$.   Then we have
\begin{equation} \label{eq:quasilin51}
\begin{aligned} &    z^t  =   z  +
 \mathcal{Z}(t)  \, \, \quad
  & f^t  =e^{Jq(t )\cdot \Diamond } ( f+ \textbf{S}(t ) )
\end{aligned}
\end{equation}
where, for   $(k,m)$ with   $k\in \Z\cap [0,r-(m+1)\textbf{d}]$
and  $1\le m \le M$,
 for $ B_{\Sigma  _{-k}}$
a sufficiently small neighborhood of 0 in $ \Sigma  _{-k}\cap X_c(p_0) $  and for $B_{\C ^{\mathbf{n} }} $ (resp.$B_{\R ^{n_0}}$)
 a  neighborhood of 0 in $\C ^{\mathbf{n} }$ (resp.$ {\R ^{n_0}}$)

\begin{equation} \label{eq:ODEpr21}\begin{aligned} &
     \textbf{S} \in C^m((-2,2)\times  B_{\C ^{\mathbf{n} }}   \times    B_{\Sigma  _{-k}}
\times B_{\R ^{n_0}} , \Sigma _{k}
 )  \\&  {q}   \in C^m((-2,2)\times  B_{\C ^{\mathbf{n} }}   \times    B_{\Sigma  _{-k}}
\times B_{\R ^{n_0}} , \R ^{ n_0}
 ) \\ &  \mathcal{Z }  \in C^m((-2,2)\times  B_{\C ^{\mathbf{n} }}   \times    B_{\Sigma  _{-k}}
\times B_{\R ^{n_0}} , \C ^{  \mathbf{n}}
 ),
\end{aligned}   \end{equation}
with for fixed $C$
 \begin{equation} \label{eq:symbol11}   \begin{aligned}   &   | q (t,z,f,\varrho ) |\le C (|z|+\| f\| _{\Sigma _{-k}}) ^{M_0+1} \\&
  |\mathcal{Z} (t,z,f,\varrho ) |+ \| \textbf{S} (t,z,f,\varrho ) \| _{\Sigma _{k}} \le C (|z|+\| f\| _{\Sigma _{-k}}) ^{M_0 } .
\end{aligned}    \end{equation}
We have $ \textbf{S} (t,z,f,\varrho )= \textbf{S}_1 (t,z,f,\varrho )+ \textbf{S}_2 (t,z,f,\varrho )$ with
\begin{equation} \label{eq:symbol12}   \begin{aligned}   &    \textbf{S}_1 (t,z,f,\varrho ) =\int _0^t
 (X_\chi ^{st} ) _{f}\circ \phi ^{t'} dt' \\&
  \| \textbf{S}_2 (t,z,f,\varrho ) \| _{\Sigma _{k}} \le C (|z|+\| f\| _{\Sigma _{-k}}) ^{2M_0+1 }(|z|+\| f\| _{\Sigma _{-k}}+|\varrho |)^i .
\end{aligned}    \end{equation}

\item[(3)] The flow $\phi ^t$ is canonical: for   $s,s',k$ as in Lemma \ref{lem:ODE},    the map $\phi  ^t \in C^l(  \U ^{s'}_{\varepsilon _1,k}
  , \widetilde{\Ph}  ^{s })$
 satisfies
   $\phi ^{t*}\Omega _0=\Omega _0$
  { in } $C^{\infty}(   \U ^{s'}_{\varepsilon _2,k}
  , B^2 (\widetilde{ \Ph}^{s'}, \R )
 )  $ for $\varepsilon  _2>0$ sufficiently small.
\end{itemize}

  \end{lemma}
	\proof   First of all notice that  $\chi$ does not depend on $\tau$ and $\Pi$ so that the only nonzero component
	of  $X_\chi$ is $(X_\chi  )_R=J\nabla _R\chi $.  The latter is  of the form indicated in claim (1) by a direct
	computation. Claim (2) follows now by Lemma \ref{lem:ODE}.
	
\noindent To prove Claim (3) we need to make rigorous  the following  formal computation
\begin{equation*}
 \frac{d}{dt}\phi ^{t*}\Omega _0=\phi ^{t*}L_{X_\chi}\Omega _0=\phi ^{t*}di_{X_\chi}\Omega _0=\phi ^{t*}d^2\chi =0.
\end{equation*}
  To make sense of this we can proceed as  in Corollary \ref{cor:darboux}. We skip the proof.
	\qed

\begin{lemma}
  \label{lem:ExpH11}  Consider a  transformation
  $\mathfrak{F} =\mathfrak{F}_1 \circ \cdots \circ \mathfrak{F}_L$
  like in Lemma \ref{lem:ODE1}  and with $m_1=2$ and for fixed $r$ and $M$ sufficiently large.
	Denote by $(k',m')$ the pair  $(k,m)$ of  Lemma \ref{lem:back11} and
	consider a pair  $(k,m)$ with $k\le k'$ and $m\le m' -(2\textbf{N}+5)$.
Set $H':=K\circ \mathfrak{F} $.  Consider decomposition
  \eqref{eq:decomp2}.  Then on a domain $\U^s_{\varepsilon , k} $ like \eqref{eq:domain0} we have
\begin{equation}  \label{eq:ExpH11} \begin{aligned} & H '=  { \psi} (\Pi (f)) +H_2'   +\textbf{R}\ ,
\end{aligned}
\end{equation}      for a $\psi \in C^\infty $
	with $\psi (\varrho )=O(|\varrho|^2)$ near 0
and with   what follows.\begin{itemize}
\item[(1)]
 We have

\begin{equation}  \label{eq:ExpH2} H_2 '=
\sum _{\substack{ |\mu +\nu |=2\\
\mathbf{e}  \cdot (\mu -\nu )=0}}
 a_{\mu \nu} ( \Pi  (f) )  z^\mu
\overline{z}^\nu + \frac{1}{2} \langle J ^{-1} \mathcal{H}_{p
_0} f, f\rangle .
\end{equation}

\item[(2)] We have $\textbf{R}  =  {\textbf{R} _{ -1 }} +  {\textbf{R} _{0}} + {\textbf{R} _1 } +
 {\textbf{R} _2} +\resto ^{1,2}_{k,m+2}(\Pi (f), f)+
 {\textbf{R} _3}  +
 {\textbf{R} _4} $, with:

\begin{equation*}   \begin{aligned} &{\textbf{R} _{ -1 }}=
 \sum _{\substack{ |\mu +\nu |=2\\
\mathbf{e} \cdot (\mu -\nu )\neq 0  }} a_{\mu \nu } (\Pi (f)
)z^\mu
\overline{z}^\nu  +\sum _{|\mu +\nu |  = 1} z^\mu
\overline{z}^\nu  \langle  J ^{-1}G_{\mu \nu }(\Pi  (f)
),f\rangle  ;\end{aligned}
\end{equation*}
For $\mathbf{N}$ as  in (L4) of this section,
    \begin{equation*}    \begin{aligned} &   {\textbf{R} _0}=   \sum _{|\mu +\nu |= 3}^{2\mathbf{N}+1} z^\mu
\overline{z}^\nu   a_{\mu \nu
}( \Pi (f) ) ;
 \end{aligned}
\end{equation*}
\begin{equation*}    \begin{aligned} &   {\textbf{R} _1}=   \sum _{|\mu +\nu
|= 2}^{2\mathbf{N} } z^\mu
\overline{z}^\nu \langle J ^{-1}  G_{\mu \nu }(  \Pi (f)
), f\rangle ;
 \end{aligned}
\end{equation*}
\begin{equation*}    \begin{aligned} &   {\textbf{R} _2}=
\langle \mathbf{B}_{2 } (   \Pi (f) ),   f  ^{   2} \rangle
     \text{ with $\mathbf{B}_{2 } (0)=0$}
\end{aligned}\nonumber
\end{equation*}
  where   $f^d(x)$    represents schematically $d-$products of components of $f$;
\begin{equation*}    \begin{aligned} &   {\textbf{R} _3}=   \sum _{ \substack{ |\mu +\nu |=\\=  2N+2}} z^\mu
\overline{z}^\nu   a_{\mu \nu
}( z,f,\Pi (f) )     +\sum _{ \substack{ |\mu +\nu
|=\\= 2N+1}} z^\mu
\overline{z}^\nu  \langle J ^{-1}  G_{\mu \nu }(z,f,  \Pi (f)
), f\rangle ;
 \end{aligned}
\end{equation*}
\begin{equation*}    \begin{aligned} &   {\textbf{R} _4}=  \sum _{d=2}^4
\langle B_{d } (   z ,f,\Pi (f) ),   f  ^{   d} \rangle
      +\int _{\mathbb{R}^3}
B_5 (x,  z ,f, f(x),\Pi (f) )  f^{   5}(x) dx\\& +   \widehat{\textbf{R}} _2(   z ,f,\Pi (f))+   E_P (  f) \text{ with $B_2(0,0,\varrho )=0$.}
\end{aligned}\nonumber
\end{equation*}

\item[(3)] For $\delta _j:=( \delta _{1j}, ..., \delta _{mj}),$
\begin{equation}  \label{eq:ExpHcoeff1} \begin{aligned} &
a_{\mu \nu } ( 0 ) =0 \text{ for $|\mu +\nu | = 2$  with $(\mu
, \nu )\neq (\delta _j, \delta _j)$ for all $j$,} \\& a_{\delta _j
\delta _j } ( 0 ) =\lambda _j (\omega _0)  ,
\\& G_{\mu \nu }(  0 ) =0 \text{ for $|\mu +\nu | = 1$ }.
\end{aligned}
\end{equation}
These $a_{\mu \nu } ( \varrho )$ and $G_{\mu \nu }( x,\varrho
)$ are $C ^{m }$ in all variables with $G_{\mu \nu }( \cdot ,\varrho )
\in C^m  ( \mathrm{U},\Sigma _k (\mathbb{R}^3,\mathbb{C}^{2N  }))$,
  for   a
small neighborhood $\mathrm{U}$ of $( 0,0,0)$ in
$ \mathbb{C}^{ \mathbf{n}}\times (\Sigma _{-k}\cap X_c(p_0))\times \R ^{n_0}$
(the space of the $(z,f,\varrho )$),
 and they satisfy
symmetries analogous to \eqref{eq:symm}.

\item[(4)] We have
$a_{\mu \nu
}(   z,  \varrho   ) \in C^{ m }( \mathrm{U},
 \mathbb{C} ) $ .

\item[(5)] $G_{\mu \nu
}( \cdot ,  z,  \varrho ) \in C^{ m } ( \mathrm{U},
\Sigma _k (\mathbb{R}^3,\mathbb{C}^{2N  }))) $.

\item[(6)] $B_{d
}( \cdot ,  z,f,\varrho  ) \in C^{m } ( \mathrm{U},
\Sigma _k (\mathbb{R}^3, B   (
 (\mathbb{C}^{2N  })^{\otimes d},\mathbb{R} ))) $, for $2\le d \le 4$.
$\mathbf{B}_{2
}( \cdot ,   \varrho  )$ satisfies the same property.

\item[(7)] Let
$ \zeta \in \mathbb{C}^{2N  }$. Then for
  $B_5(\cdot ,  z ,f, \zeta  ,\varrho )$   we have  (the derivatives are not in the holomorphic sense)
\begin{equation} \label{H5power2}\begin{aligned} &\text{for  $|l|\le  m $ ,
}  \| \nabla _{  z ,f,\zeta, \varrho  }
^lB_5(  z,f,\zeta  ,\varrho  ) \| _{\Sigma _k(\mathbb{R}^3,   B   (
 (\mathbb{R}^{2N  })^{\otimes 5},\mathbb{R} )} \le C_l .
 \end{aligned}\nonumber \end{equation}
 \item[(8)]
\begin{equation}\label{eq:Rhat0}\begin{aligned} &
\widehat{\textbf{R}} _2
 \in C^{m} ( \mathrm{U} ,\C
 ),  \\&    | \widehat{\textbf{R}} _2 (z ,f, \varrho
)|  \le C (|z|+  \| f \| _{\Sigma _{-k}}) \| f \| _{\Sigma _{-k}}^2;
\end{aligned}\end{equation}

 \end{itemize}
\end{lemma}
\proof     We need to express $R$ in terms of $(z,f)$
 using \eqref{eq:decomp2} inside \eqref{eq:back11}.

\noindent   We have $\Pi (R)= \Pi (f) +\resto ^{0,2}(R).$ Then, succinctly,
\begin{equation*}\begin{aligned}
&
\resto ^{1,2}_{k',m' }(\Pi (R),R)=\sum _{a+b=2}	^{2\textbf{N}+1} \frac{1}{a!b!} \langle
 \nabla ^a_{ \varrho } \nabla ^b_{ R }\resto ^{1,2}_{k',m' }(\Pi (f),0)  , (\resto ^{0,2}(R) )^{a }  R ^{b\otimes} \rangle  + \\&   \sum _{\substack{a+b\\ =2\textbf{N}+2}} \int _0^1
\frac{(1-t)^{2\textbf{N}+1}}{a!b!} \langle \nabla ^a_{ \varrho } \nabla ^b_{ R }
\resto ^{1,2}_{k',m' }(\Pi (f)+t\resto ^{0,2}(R),tR) ,  (\resto ^{0,2}(R) )^{a }  R ^{b\otimes}    \rangle  dt,
   \end{aligned}
\end{equation*}
with   $(k',m')$ the pair    $(k,m)$  of Lemma \ref{lem:back11}.
We substitute \eqref{eq:decomp2}, that is $R=z\cdot \xi +\overline{z}\cdot \overline{\xi} +f$.
 For   $m\le m' -(2\textbf{N}+2)$ and $k\le k'$,
 the terms from the $R^{b\otimes}$
  of degree in $f$  at most  1,   go into $\textbf{R}_i$ with $i=-1,0,1,3$ and $H'_2$.
	For   $m\le m' -(2\textbf{N}+4)$,
 the remaining  terms are absorbed in $\resto ^{1,2}_{k ',m+2 }(\Pi (f),f)+\widehat{\textbf{R}}_2(z,f,\Pi (f)) $.

\noindent   We focus now on the $d=5$ term in \eqref{eq:back11}.  We substitute $R=z\cdot \xi +\overline{z}\cdot \overline{\xi} +f$.
This schematically   yields, for a $\widetilde{B}_5$ satisfying claim (7) with the pair $(m',k')$,

\begin{equation}\label{eq.b5}\begin{aligned}
& \sum _ {j=0}^5\int _{\mathbb{R}^3}
\widetilde{B}_5 (x, z, f, f(x),\Pi (f) )  (z\cdot \xi +\overline{z}\cdot \overline{\xi} ) ^{5-j}  f^j   ( x) dx
  . \end{aligned}
\end{equation}
 For $j=5$ we get a term that can be absorbed  in the $B_5$ term in $\textbf{R} _4$.  Expand  the $j<5$ terms
in \eqref{eq.b5} as

\begin{equation*}\label{eq.b51}\begin{aligned}
& \sum _ {i=0}^{4-j}\int _{\mathbb{R}^3} \frac{1}{i!} (\partial _t^i ) _{|t=0}
\widetilde{B}_5 (x, z, f,t f(x),\Pi (f) )  (z\cdot \xi +\overline{z}\cdot \overline{\xi} ) ^{5-j}  f^{i+j }  ( x) dx+\\&
 \int _{\mathbb{R}^3} \frac{1}{(4-j)!} \int _0^1   \partial _t^{5-j} [\widetilde{B}_5 (x, z, f,t f(x),\Pi (f) )]  (z\cdot \xi +\overline{z}\cdot \overline{\xi} ) ^{5-j}  f^5   ( x) dx
  . \end{aligned}
\end{equation*} go into the $B_d$ term in $\textbf{R} _4$
The last term fits   in the $B_5$ term in $\textbf{R} _4$ by $m\le m'-5$.  The terms in the first line
go into the $B_d$ of $\textbf{R}_4$
for $d=i+j\ge 2$ .  The terms with $ i+j< 2$ can be treated like the $\resto ^{1,2}_{k',m' }(\Pi (R),R)$ for
 $m\le m' -(2\textbf{N}+5)$ and $k\le k'$.

\noindent We  focus on $E_P(R)=E_P(z\cdot \xi +\overline{z}\cdot \overline{\xi} +f)$. We use Lemma \ref{lem:Taylor ex} for $\textbf{v}=f$ and  $ \textbf{u}= z\cdot \xi +\overline{z}\cdot \overline{\xi}.$ Then
\begin{equation*}       \begin{aligned} & E_P(  R )= E_P( f   ) + E_P( z\cdot \xi +\overline{z}\cdot \overline{\xi} ) +\\& \int _{\R ^3}dx
\sum _{ j =0}^{3}\int _{[0,1]^2} \frac{t^j}{j!}   (\partial _t^{j+1})_{|t=0} \partial _s   [B (|s ( z\cdot \xi +\overline{z}\cdot \overline{\xi} )+tf  |^2_1) ]  dt ds   \\& + \int _{\R ^3}dx \int _{[0,1]^2}  dt ds \int _0^t  \partial _\tau ^5 \partial _s [B (|s ( z\cdot \xi +\overline{z}\cdot \overline{\xi} )   +\tau f   |_1^2)] \frac{(t-\tau )^3}{3!} d\tau  .
\end{aligned}
\end{equation*}
By   $B(0)=B'(0)=0$,   we have $E_P( z\cdot \xi +\overline{z}\cdot \overline{\xi} )=\resto ^{0,4} (R)$.
It is easy to conclude that this term easily fits into $\textbf{R}_0 +\textbf{R}_3$. Similarly, the $j=0$ term
fits in $ \textbf{R}_1+\textbf{R}_3$. The  $j\ge 1$ terms
fit  in   the $B_{j+1}$ term in $  \textbf{R}_4$.   The last line fits in the $B_ {5}$ term in $  \textbf{R}_4$.

\noindent  The symmetries \eqref{eq:symm} for the coefficients in $H'_2+\textbf{R}_{-1}+\textbf{R}_0+ \textbf{R}_1 $ are an elementary consequence of the fact that $H'$ is real valued.

\qed

\begin{remark}
\label{rem:differences}
 Given a Hamiltonian $H'$ expanded as in Lemma \ref{lem:ExpH11} and given
a transformation $\mathfrak{F} $, we cannot obtain the expansion of
Lemma \ref{lem:ExpH11}  for $H'\circ\mathfrak{F}$ analysing   one by one the terms of the expansion of
$H'$. This works in the set up of \cite{Cu1,Cu2}  but not here  (see in particular the discussion on the exponential under formula \eqref{eq:key-201}
later).
\end{remark}

\section{ Birkhoff normal forms}
\label{sec:Normal form}
In this section we arrive at the main result of the paper.

\subsection{Homological equations}
\label{subsec:homological}

       We consider   $a_{\mu \nu}^{(\ell )} ( \varrho   )\in C ^{ \widehat{m}} (U,C)$ for $k_0\in \N$ a fixed number and $U$ a neighborhood of 0 in $\R ^{n_0}$.
Then we set

\begin{equation}  \label{eq:H2birk} H_2^{(\ell )}  (\varrho ):=
\sum _{\substack{ |\mu +\nu |=2\\
\mathbf{e}    \cdot (\mu -\nu )=0}}
 a_{\mu \nu}^{(\ell )} ( \varrho   )  z^\mu
\overline{z}^\nu + \frac{1}{2} \langle J^{-1} \mathcal{H} f, f\rangle .
\end{equation}

\begin{equation} \label{eq:lambda}
\textbf{e} _j  ( \varrho  ) :=   a
_{\delta _j\delta _j} ^{(\ell )}(\varrho   ), \quad \textbf{e} (\varrho )=
(\lambda _1(\varrho ), \cdots, \lambda _m(\varrho )).\end{equation}
We assume  $\textbf{e} _j (0 ) =  \textbf{e}  _j     $ and   $a_{\mu \nu}^{(\ell )}  ( 0 ) =0 $ if $(\mu ,\nu)\neq (\delta _j,\delta _j)$ for
all $j$, with $\delta _j$ defined in \eqref{eq:ExpHcoeff1}.

\begin{definition}

\label{def:normal form} A function $Z( z,f, \varrho )$ is in normal form if  $
  Z=Z_0+Z_1
$
where $Z_0$ and $Z_1$ are finite sums of the following type:
\begin{equation}
\label{e.12a}Z_1= \sum _{\mathbf{e} (0)   \cdot(\nu-\mu)\in \sigma _e(\mathcal{H} _{p_0})}
z^\mu \overline{z}^\nu \langle  J^{-1} G_{\mu \nu}( \varrho   ),f\rangle
\end{equation}
with $G_{\mu \nu}( x ,\varrho )\in  C^{m} (  U,\Sigma _{k }(\R ^3, \C ^{2N}))$ for  fixed $k,m\in \N$ and $U\subseteq \R ^{n_0}$ a neighborhood of 0;
\begin{equation}
\label{e.12c}Z_0= \sum _{  \mathbf{e}(0)  \cdot(\mu-\nu)=0} g_{\mu   \nu}
( \varrho  )z^\mu \overline{z}^\nu
\end{equation}
and $g_{\mu   \nu} ( \varrho  )\in  C^{m} (   U,
\mathbb{C})$.
We assume furthermore that the above coefficients
satisfy the symmetries in \eqref{eq:symm}: that is $\overline{g}_{\mu \nu}=g_{\nu\mu  }$ and $\overline{G}_{\mu \nu}=G_{\nu \mu }$.
 \end{definition}

\begin{lemma}
\label{lem:NLhom1} We consider $\chi =\chi (b,B)$ with
\begin{equation}
\label{eq:chi}\chi (b,B) =\sum _{|\mu +\nu |=M_0+1} b_{\mu\nu}   z^{\mu} \overline{z}^{\nu} + \sum _{|\mu +\nu |=M_0} z^{\mu}
\overline{z}^{\nu}
 \langle  J^{-1}  B_{\mu   \nu
}
  , f \rangle
\end{equation}
for $b_{\mu\nu} \in \C $ and $B_{\mu   \nu
} \in  \Sigma  _{\widehat{k} }(\R ^3, \C ^{2N})\cap X_c(p_0)$ with $\widehat{k} \in \N$, satisfying the symmetries in \eqref{eq:symm}. Here we interpret
the polynomial  $\chi$ as a function with parameters $b=(b_{\mu\nu})$ and $B=(B_{\mu\nu})$. Denote by $X_{\widehat{k} }$ the space of the pairs $(b,B)$.
Let us also consider given polynomials with  $K=K( \varrho  ) $  and $\widetilde{K}=\widetilde{K}( \varrho ,b,B  ) $ where:

\begin{equation}\label{eq:Krho}
K(\varrho ) :=\sum _{|\mu +\nu |=M_0+1} k_{\mu\nu} ( \varrho   ) z^{\mu} \overline{z}^{\nu} + \sum _{|\mu +\nu |=M_0} z^{\mu} \overline{z}^{\nu}
 \langle  J^{-1}  K_{\mu   \nu
}(\varrho  )
  , f \rangle ,
\end{equation}
with $k_{\mu\nu} ( \varrho   ) \in C^ {\widehat{m} } ( U, \C )$ and $K_{\mu\nu} ( \varrho   ) \in C^{\widehat{m} } ( U,    \Sigma  _{\widehat{k}}(\R ^3, \C ^{2N})  \cap X_c(p_0))$ for $U$ a neighborhood of 0 in $\R ^{n_0}$,
satisfying the symmetries in \eqref{eq:symm};

\begin{equation}\label{eq:tildeKrho} \begin{aligned} &\widetilde{K} ( \varrho ,b,B  ) :=
 \sum _{|\mu +\nu |=M_0+1} \widetilde{k}_{\mu\nu}  ( \varrho , b,B  ) z^{\mu} \overline{z}^{\nu} \\&+
\sum _{i=0}^1\sum _{j=1}^{n_0}
 \sum _{|\mu +\nu |=M_0} z^{\mu}
\overline{z}^{\nu}
 \langle   J ^{-1} \Diamond  ^i_j  K_{j\mu   \nu
}^i(\varrho , b,B   )
  , f \rangle ,\end{aligned}
\end{equation}
with $\widetilde{k}_{\mu\nu}   \in C^{\widehat{m}} ( U \times X _{ \widehat{k} }, \R )$ and $\widetilde{K}_{j\mu\nu}^{i}   \in C^{\widehat{m}} ( U\times X_{\widehat{ k} },   \Sigma _{ \widehat{k} }(\R ^3, \C ^{2N})  \cap X_c(p_0))$,
satisfying the symmetries in \eqref{eq:symm}.
 Suppose also that the sums \eqref{eq:Krho} and \eqref{eq:tildeKrho} do not
 contain terms in normal form and that $\widetilde{K}( 0 ,b,B  )=0$. Then there exists
 a  neighborhood  $V\subseteq U$ of  0 in $\R ^{n_0}$  and a unique choice of functions
 $(b(\varrho  ), B(\varrho   )) \in C ^{\widehat{m}} (V,X _{\widehat{k}} )$ such that for $\chi (\varrho  )
  =\chi (b(\varrho  ), B(\varrho  ))$,    $\widetilde{K}(\varrho  ) =\widetilde{K} (\varrho ,b(\varrho  ),B(\varrho  ))$  we have

\begin{equation}
\label{eq:homologicalEq} \left\{ \chi (\varrho  ) ,H_2 ^{(\ell )}(\varrho  )  \right \} ^{st}= K   (\varrho)+\widetilde{K}(\varrho  ) +Z(\varrho    )
\end{equation}
where  $\{
\cdots   \} ^{st} $ is the  bracket  \eqref{eq:poiss} for $\varrho$
fixed and
where $Z(\varrho  )$ is in normal form and homogeneous of degree $M_0+1$ in $(z,f)$.

\end{lemma}
\proof    Summing on repeated indexes,  by \eqref{eq:poiss1} we get \begin{equation} \label{eq:homologicalEq1}  \begin{aligned}&    \{
H_2^{(\ell )}, \chi  \} ^{st} =    -\im  \mathbf{e}  ( \varrho )\cdot (\mu -
\nu) z^{\mu} \overline{z}^{\nu} b_{\mu \nu }(\varrho  ) \\&- z^{\mu} \overline{z}^{\nu}\langle f , J ^{-1}( \im \mathbf{e}  ( \varrho )\cdot (\mu -
\nu)- \mathcal{H} )B _{\mu \nu }(\varrho  )
\rangle   + \widehat{ K}  (\varrho   ,b(\varrho  ),B(\varrho  ))   ,
\end{aligned}
\end{equation}
  \begin{equation} \label{eq:homologicalEq2}  \begin{aligned}  \widehat{ K}  (\varrho ,  b,B) :&=   \sum _{\substack{|\mu +\nu |=2\\ (\mu , \nu )\neq (\delta _j , \delta _j ) \ \forall \ j}  }  a_{\mu \nu}^{(\ell   )}(\varrho  )
\big [\sum _{|\mu '+\nu '|=M_0  +1}
\{z^\mu
\overline{z}^\nu ,z ^{\mu'}
\overline{z}^{ \nu '} \}  b_{\mu ' \nu '}   \\& +
\sum _{|\mu '+\nu '|=M_0   }
\{z^\mu
\overline{z}^\nu ,z ^{\mu'}
\overline{z}^{ \nu '} \}  \langle  J^{-1}  B_{\mu ' \nu '}
  , f \rangle
 \big ]  .
\end{aligned}
\end{equation}
 $\widehat{ K} $    is a homogeneous polynomial of the same type of the above ones and we have $\widehat{ K}  (0  , b,B) =0$.   In particular,  $\widehat{ K} $  satisfies the symmetries \eqref{eq:symm}  by (for $f=\overline{f}$)

\begin{equation*}    \begin{aligned}  & ( a_{\mu \nu}^{(\ell   )}   b_{\mu ' \nu '}
\{z^\mu
\overline{z}^\nu ,z ^{\mu'}
\overline{z}^{ \nu '} \}  )^*   = a_{\nu \mu }^{(\ell   )}   b_{\nu '\mu '  }
\{z^\nu
\overline{z}^\mu ,z ^{\nu'}
\overline{z}^{ \mu '} \}    \,
 \\&  ( a_{\mu \nu}^{(\ell   )}   \langle  J^{-1}  B_{\mu ' \nu '}
  , f \rangle
\{z^\mu
\overline{z}^\nu ,z ^{\mu'}
\overline{z}^{ \nu '} \}  )^*   = a_{\nu \mu }^{(\ell   )}   \langle  J^{-1}  B_{\nu ' \mu '}
  , f \rangle
\{z^\nu
\overline{z}^\mu ,z ^{\nu'}
\overline{z}^{ \mu '} \}
\end{aligned}
\end{equation*}
which follow by  $(\im \partial _{z_j}F\partial _{\overline{z}_j}G
-\im \partial _{\overline{z}_j}F\partial _{z_j}G)^*=\im \partial _{z_j}F^*\partial _{\overline{z}_j}G^*
-\im \partial _{\overline{z}_j}F^*\partial _{z_j}G^*$, where in these formulas $a^*=\overline{a}$,  and by the symmetries \eqref{eq:symm} for $\chi$ and for $H_2^{(\ell )}$   .

 \noindent Denote by $\widehat{ Z}  (\varrho ,  b,B) $ the sum of
monomials in normal form of $\widetilde{{K}}$ and set $\textbf{K } :=\widetilde{{K}} +\widehat{ K}  -\widehat{ Z}   $.
 We look at   \begin{equation} \label{eq:homologicalEq10}  \begin{aligned} &   - \im \mathbf{e}  (\varrho )\cdot (\mu -
\nu) z^{\mu} \overline{z}^{\nu} b_{\mu \nu }    -z^{\mu} \overline{z}^{\nu}    \langle f , J ^{-1}(\im \mathbf{e}   (\varrho ) \cdot (\mu -
\nu)- \mathcal{H} )B _{\mu \nu }
\rangle    \\& +  \textbf{K }(\varrho   ,b ,B )  +K(\varrho  ) =0
\end{aligned}
\end{equation}
 that is  at

\begin{equation} \label{eq:NLhom11}\begin{aligned}&
k_{\mu \nu} (\varrho )  + \textbf{{k}}_{\mu \nu}(\varrho   ,b,B)  - b_{\mu \nu}(\varrho  )  \im  \mathbf{e}  (\varrho )
\cdot (\mu - \nu)=0   \\&  B_{\mu   \nu } (\varrho  )  =- R_{\mathcal{H}}(\im \mathbf{e}   (\varrho ) \cdot (\mu - \nu) )\left [ {K}_{\mu \nu }(\varrho ) +  \textbf{{K}} _{\mu \nu}(\varrho   ,b,B)\right ]
   ,
\end{aligned}
\end{equation}
with $\textbf{{k}}_{\mu \nu}$ and  $\textbf{{K}} _{\mu \nu}$ the coefficients of $\textbf{K}$.
Notice that when    $ \textbf{k}_{\mu \nu}(0  ,b,B) =0$ and $\textbf{K}_{\mu \nu}(0 ,b,B) =0$,  for $\varrho = 0$ there is a unique solution $(b,B)\in X_{ \widehat{k} }$ given by
\begin{equation} \label{eq:NLhom12}\begin{aligned}&
  b_{\mu \nu}(0)
 =\frac{   k_{\mu\nu}(0)}{\im
    \mathbf{e}     \cdot(\mu-\nu)}   \, , \quad  B_{\mu   \nu } (0 )  =-  R _{\mathcal{H}}      (\im \mathbf{e}   \cdot
 (\mu -\nu )  )  K_{\mu   \nu
}(0)  .
\end{aligned}
\end{equation}
  Lemma \ref{lem:NLhom1} is then  a consequence of   the Implicit Function Theorem by   Hypothesis (L7) in Sect. \ref{sec:speccoo}.
\qed

In the particular case $M_0=1$ we need a slight variation of   Lemma \ref{lem:NLhom1}.

\begin{lemma}
\label{lem:NLhom2}  Suppose now $M_0=1$ and assume  the notation  of  Lemma \ref{lem:NLhom1},
  assuming  $K(0 )=0$,  $\widetilde{K}( 0 ,0,0  )=0$  and $\nabla _{b,B}\widetilde{K}( 0 ,0,0  )=0$. We furthermore
consider function $a_{\mu \nu}^{\mu ' \nu '} \in C^{\widehat{m}} (U\times X_{\widehat{k}}, \C)$
with  $| a_{\mu \nu}^{\mu ' \nu '} (\varrho , b, B)|\le C \| (b, B)\| _{X_{\widehat{k}}}$ and we set

\begin{equation}
\label{eq:modop} \begin{aligned} & \left\{ \chi (\varrho  ) ,H_2 ^{(\ell )}(\varrho  )  \right \} ^{\widetilde{st}}=
\left\{ \chi (\varrho  ) ,H_2 ^{(\ell )}(\varrho  )  \right \} ^{st} \\&  +
 \sum _{\substack{|\mu +\nu |=1\\  |\mu '+\nu '|=1}}
	a_{\mu \nu}^{\mu ' \nu '} ( \varrho ,b(\varrho ),B(\varrho )) z^{\mu} \overline{z}^{\nu}\langle
	\mathcal{H} B_{\mu ' \nu '}(\varrho ), f
\rangle .\end{aligned}
\end{equation}
Then, the same conclusions of Lemma  \ref{lem:NLhom1} hold for
\begin{equation}
\label{eq:modhomologicalEq} \left\{ \chi (\varrho  ) ,H_2 ^{(\ell )}(\varrho  )  \right \} ^{\widetilde{st}}= K   (\varrho)+\widetilde{K}(\varrho  ) +Z(\varrho    ) .
 \end{equation}

\end{lemma}
\proof     Like above we get  to
\begin{equation*} \label{eq:NLhom111}\begin{aligned}&
k_{\mu \nu} (\varrho )  + \textbf{{k}}_{\mu \nu}(\varrho   ,b,B)  - b_{\mu \nu}  \im  \mathbf{e}  (\varrho )
\cdot (\mu - \nu)=0   \\&  B_{\mu   \nu }   =- R_{\mathcal{H}}(\im \mathbf{e}   (\varrho ) \cdot (\mu - \nu) )  [ {K}_{\mu \nu }(\varrho ) +  \textbf{{K}} _{\mu \nu}(\varrho   ,b,B) + \sum _{\mu ' \nu '}a_{\mu \nu}^{\mu ' \nu '}(\varrho , b ,B )  \mathcal{H} B_{\mu ' \nu '}   ]
  .
\end{aligned}
\end{equation*}
For $(\varrho , b, B) =(0 ,0,0) $ both sides are 0. Then
  Lemma \ref{lem:NLhom2}  follows by Implicit Function Theorem.
\qed

\subsection{The Birkhoff normal forms}
\label{subsec:mainbirkhoff}

 Our   goal in this section is to prove the
following result where $N$ is as of (L4) in Sect. \ref{sec:speccoo}.

\begin{theorem}
\label{th:main}   For any integer $2\le \ell \le  2\mathbf{N}+1$
we have  transformations $\mathfrak{F} ^{(\ell )} = \mathfrak{F}_1 \circ \phi _2\circ ...\circ \phi _\ell   $, with   $\mathfrak{F}_1$  the transformation in
Corollary \ref{cor:darboux}
the $\phi _j$'s
  like in Lemma \ref{lem:chi}, such that
  the conclusions of  Lemma \ref{lem:ExpH11} hold,
    that is such that we have the following expansion
\begin{equation*}
    H ^{(\ell )}:=K\circ  \mathfrak{F} ^{(\ell )}   ={ \psi} (\Pi (f))+ H_2^{(\ell )}   + \resto^{1,2}_{k,m+2}(\Pi (f) , f)+ \sum _{j=-1}^4 \textbf{R}_{j}^{(\ell )},
\end{equation*}
with $H_2^{(\ell )}$ of the form \eqref{eq:ExpH2} and with the following additional properties:
\begin{itemize}
\item[(i)]  $\textbf{R} _{-1} ^{(\ell )} =0$;

\item[(ii)]  all the nonzero terms in $\textbf{R} _0 ^{(\ell )} $ with $|\mu +\nu |\le \ell $  are in normal form,
that is $\lambda \cdot (\mu -\nu )=0$;

\item[(iii)] all the nonzero terms in $\textbf{R} _1^{(\ell )}  $ with $|\mu +\nu |\le \ell -1$  are in normal form,
that is $\lambda \cdot (\mu -\nu )\in \sigma _e(\mathcal{H} _{p_0})$.
\end{itemize}
\end{theorem}

\proof   The proof of  Theorem \ref{th:main} is   by induction.   There are two distinct
parts in the proof,   \cite{Cu2,Cu1,bambusi}.  Here we  follow  the ordering of
\cite{bambusi}.  In the first part
we
assume that for some $\ell \ge 2$ the statement of the theorem is true, and we show
that it continues to be true for $\ell+1$.
The proof of case  $\ell =2$, which presents some additional complications, is dealt   in the
second part.

In the proof we will get polynomials \eqref{eq:chi1}    with $M_0=1,..., 2\textbf{N}  $ with decreasing
$(r,M)$ as  $M_0$ increases.  Nonetheless, in view of the fact that in Lemma \ref{lem:fred12}  the $n$ is arbitrarily large and  that   $(r,M)$ decreases by a fixed amount at each step, these  $(r,M)$ are
arbitrarily large. This is exploited in Theorem \ref{theorem-1.1} later.

\bigskip
\noindent  {\it The step $\ell +1>2$.} We can assume that
 $H^{(\ell)}$ have the desired properties for  indexes $(k',m')$ (instead of $(k,m)$) arbitrarily large.
We consider the representation \eqref{eq:ExpH11} for $H^{(\ell)}$ and we set $\mathbf{h}=H^{(\ell)}(z,f,\varrho )$    replacing
  $\Pi (f)$ with $\varrho$ in \eqref{eq:ExpH11}.  Then  $\mathbf{h}=H^{(\ell)}(z,f,\varrho )$ is
	$C^{2\mathbf{N}+2}$
	near 0 in
$\Ph ^{s _0}=\{  (\varrho , R ) \}$ for $m '\ge    2\mathbf{N}+2 $  for $s_0>\max \{
 \text{ord}(\mathcal{H}_{p_0}), 3/2 \} $
by
   Lemma \ref{lem:ExpH11}.
 So we have equalities
\begin{align}
   \label{eq:derivmain1}  & a_{\mu \nu}(\varrho )
=\frac{1}{\mu !\nu ! }  \partial _z ^\mu \partial _{\overline{z}}^\nu
 \mathbf{h}_{|(z,f,\varrho )= (0,0,\varrho )}   \  ,  \quad |\mu +\nu |\le 2\mathbf{N}+1 ,\\&  \label{eq:derivmain2}J^{-1}G_{\mu \nu}(\varrho )
=\frac{1}{\mu !\nu ! } \partial _z ^\mu \partial _{\overline{z}}^\nu \nabla _f
 \mathbf{h} _{|(z,f,\varrho )= (0,0,\varrho )} \  ,  \quad |\mu +\nu |\le 2\mathbf{N}  .
\end{align}
We consider now  a yet unknown $\chi$ as in  \eqref{eq:chi1}   with $M_0=\ell$,
 $i=0$,  $M=m'$ and $r=k'$.
Set $\phi :=\phi ^1$, where  $\phi ^t$ is the flow
of Lemma \ref{lem:chi}. We are seeking $\chi$  such that $H^{(\ell )}
\circ\phi $ satisfies the conclusions of Theorem \ref{th:main} for $\ell +1$.

 \noindent   We know  that $H^{(\ell )}
\circ\phi $ satisfies the conclusions of  Lemma \ref{lem:ExpH11}. Therefore, to prove the induction step,
    all we need to do is to check that    the  expansion of    $H^{(\ell )}
\circ\phi $    satisfies  $\mathbf{R}_{-1}=0$ and that the only terms in
   $\mathbf{R}_{0} $ and $\mathbf{R}_{ 1} $ of degree $\le \ell +1$
   are in normal form. We have
\begin{equation}  \label{eq:backH}    \begin{aligned} &   H_2 ^{(\ell )}\circ
\phi  = H_2^{(\ell )} +  \int _0^1 \{  H_2^{(\ell )}  , \chi \}  ^{st}\circ \phi ^t dt  + \int _0^1  (\partial _{\varrho _j} a_{\mu \nu }z^\mu {\overline{z}}^\nu \{  \Pi _j(f)    , \chi \}  )  \circ \phi ^t dt.
\end{aligned}
\end{equation}
By  \eqref{eq:homologicalEq1}--\eqref{eq:homologicalEq2} we have
for $\varrho =\Pi (f)$
\begin{equation}  \label{eq:key0}   \begin{aligned} &
 \{  H_2 ^{(\ell )}   , \chi \} ^{st}=    -\im  \sum _{|\mu +\nu |=\ell +1}\mathbf{e}  ^{(\ell)}  ( \varrho )\cdot (\mu -
\nu) z^{\mu} {\overline{z}}^{\nu} b_{\mu \nu }  ( \varrho ) \\& -  \sum _{|\mu +\nu |=\ell  }  z^{\mu} {\overline{z}}^{\nu} \langle  J^{-1} ( \im \mathbf{e}  ^{(\ell)}  (   \varrho   )\cdot (\mu -
\nu)- \mathcal{H} )B _{\mu \nu }( \varrho   ),f
\rangle  + \\&  \sum _{\substack{|\mu +\nu |=2\\ (\mu , \nu )\neq (\delta _j , \delta _j ) \ \forall \ j}  }  a_{\mu \nu}^{(\ell   )}( \varrho   )
\big [\sum _{|\mu '+\nu '|=\ell +1}
\{ z^\mu
{\overline{z}}^\nu ,  z^{\mu'}
{\overline{z}}^{ \nu '} \}  b_{\mu ' \nu '} ( \varrho  )  \\& +
\sum _{|\mu '+\nu '|=\ell  }
\{  z^\mu
{\overline{z}}^\nu ,z ^{\mu'}
{\overline{z}}^{ \nu '} \}  \langle  J^{-1}  B_{\mu ' \nu '}( \varrho   )
  , f \rangle
 \big ] .
\end{aligned}
\end{equation}
By Lemma \ref{lem:chi}  for $M_0=\ell$,
 $i=0$,  $M=m'$ and $r=k'$   for first and last formula and by the proof
of Lemma \ref{lem:ODE}, in particular by  \eqref{eq:gronwall1}, we have
\begin{align}
 &    z \circ \phi ^t  =   z  +
 \resto ^{0, \ell }_{k^{\prime \prime }, m'} (t,\Pi (f) , R) \,  , \,  \Pi  (f)   \circ \phi ^t   = \Pi  (f)
+\resto ^ {0, \ell +1} _{k^{\prime \prime },m'}(t, \Pi (f),R)  \, ,  \nonumber \\&  \label{eq:backcoo}
    f \circ \phi ^t =e^{J\resto ^{0, \ell +1 }_{k^{\prime \prime }, m'} (t,\Pi (f) , R) \cdot \Diamond } ( f+ \textbf{S} ^{0, \ell }_{k^{\prime \prime }, m'}(t,\Pi (f) , R) )
\end{align}
for     $ k^{\prime \prime } \le k'-(m'+1)\textbf{d} $.  Then, substituting \eqref{eq:backcoo} in \eqref{eq:key0} we get, if $
k\le k^{\prime \prime }-\text{ord}(\mathcal{H}_{p_0})$, where $\text{ord}(\mathcal{H}_{p_0})\le \max \{\text{ord}(\mathcal{D}), \textbf{d}\}  $,
   for $1\le m \le m'$ and exploiting  that an $\resto ^ {0,2\ell  } _{k ,m }$ is also an $\resto ^ {0,\ell +2} _{k ,m }$ for $\ell \ge 2$,

\begin{equation}  \label{eq:key-11}   \begin{aligned} & \int _0^1\{  H_2 ^{(\ell )}   , \chi \} ^{st}   \circ \phi ^t dt = \{  H_2 ^{(\ell )}   , \chi \} ^{st}    +\resto ^ {0,\ell +2} _{k ,m }(\Pi (f),R) .
\end{aligned}
\end{equation}

\noindent We have
\begin{equation*}      \begin{aligned} & \{  \Pi _j(f)    , \chi \}    =\sum _{k=1}^{n_0}\{  \Pi _j(f)    , \Pi _k(f) \} \partial _{\Pi _k(f)} \chi + \sum _{|\mu '+\nu '|=\ell  }z ^{\mu'}
{\overline{z}}^{ \nu '}\langle P_c^*(p_0) \Diamond _j f,   B_{\mu ' \nu '} \rangle.
\end{aligned}
\end{equation*}
We have, for $P_d(p_0)=1-P_c(p_0)$ the projection on  the direct sum of
$N_g(\mathcal{H}_{p_0})$  and the complement of $X_c(p_0)$ in
\eqref{eq:spectraldecomp}, and using $JP_c^*(p_0)=P_c (p_0)J$  which follows from  \eqref{eq:begspectdec1},
\begin{equation}  \label{eq:keyPoiss1}   \begin{aligned} &
\{ \Pi _i(f), \Pi _j(f) \}= \langle P_c^*(p_0) \Diamond _i f, J P_c^*(p_0) \Diamond _j f\rangle \\& = \langle   \Diamond _i f, P_d (p_0)J  \Diamond _j f\rangle = \resto ^ {0,2} (f).
\end{aligned}
\end{equation}
Notice also that,  for $B_{\mu \nu} \in \Sigma _{k'} $ independent of $\Pi (f)$ and for $|\mu +\nu |=\ell$,  we have
\begin{equation}  \label{eq:keyPoiss2}   \begin{aligned} &
\{ \Pi _i(f),  z ^{\mu }
{\overline{z}}^{ \nu  }    \langle  J^{-1}  B_{\mu   \nu }
  , f \rangle \}= z ^{\mu }
{\overline{z}}^{ \nu  }   \langle P_c^*(p_0) \Diamond _i f,  B_{\mu   \nu } \rangle = \\&
      z ^{\mu }
{\overline{z}}^{ \nu  }   \langle    f,  \Diamond _i B_{\mu   \nu } \rangle - z ^{\mu }
{\overline{z}}^{ \nu  }   \langle   P_d^*(p_0) \Diamond _i   f,   B_{\mu   \nu } \rangle
= \resto ^ {0,\ell +1} _{k'-\mathbf{d} ,\infty} (R) +\resto ^ {0,\ell +1}(R)   .
\end{aligned}
\end{equation}
By \eqref{eq:keyPoiss1}--\eqref{eq:keyPoiss2}  we conclude that $\{  \Pi _j(f)    , \chi \}    =\resto ^ {0,\ell +1} _{k'-\textbf{d} ,m '}(\Pi (f),R) $.
By   \eqref{eq:backcoo} we get for $m\le m'$
\begin{equation*}     \begin{aligned} &  \{  \Pi _j(f)    , \chi \} \circ \phi ^t  =
 \resto ^ {0,\ell +1} _{k'-\textbf{d} ,m'}\left (\Pi  (f)+ \resto ^ {0, \ell +1} _{k^{\prime \prime },m'}   (t, \Pi (f),R) ,\textsl{S}      \right ) ,\\& \text{for }  \textsl{S}:=e^{J\resto  ^{0, \ell +1 }_{k^{\prime \prime }, m'}
(t,\Pi (f) , R) \cdot \Diamond } \left ( R+ \textbf{S} ^{0, \ell }_{k^{\prime \prime }, m'}(t,\Pi (f) , R) \right ).
\end{aligned}
\end{equation*}
Then
\begin{equation}  \label{eq:key-21}   \begin{aligned} &  \{  \Pi _j(f)    , \chi \} \circ \phi ^t  =   \resto ^ {0, \ell +1} _{k^{\prime \prime } -m'\textbf{d},m'}(t, \Pi (f),R) .
\end{aligned}
\end{equation}
By   \eqref{eq:backcoo}  and \eqref{eq:key-21}      the last term in \eqref{eq:backH}  is $\resto ^ {0,\ell +2} _{k ,m }(\Pi (f),R)$  for $k\le  k^{\prime \prime } -m'\textbf{d}$. This and  \eqref{eq:key-11} yield
for $ {k}=\min \{   k^{\prime   } -(2m'+1)\textbf{d} ,k^{\prime   } -( m'+1)\textbf{d} -\text{ord} (\mathcal {H}_{p_0})  \} $
\begin{equation} \label{eq:key001}\begin{aligned} &   H_2^{(\ell )}  \circ
\phi  = H_2^{(\ell )}  +    \{  H_2^{(\ell )}  , \chi \}  ^{st}+ \resto ^ {0,\ell  +2}_{\widetilde{k},m}(\Pi (f),R).
\end{aligned}
\end{equation}
 A  second  observation is that     $\mathbf{h}=(H^{(\ell)}\circ \phi ) (z,f,\varrho )$ is  $C^{2\textbf{N}+2}$
  in $\Ph ^{s_0} =\{  (\varrho , R ) \}$  for $m\ge  2\textbf{N}+2$.
We can compute again the corresponding coefficients in \eqref{eq:derivmain1}--\eqref{eq:derivmain2}.
Because of  \eqref{eq:symbol11},   for $|\mu +\nu |
\le \ell $ in \eqref{eq:derivmain1} and   for $|\mu +\nu |
\le \ell -1$ in \eqref{eq:derivmain2} these coefficients are the same of   $\mathbf{h}= H^{(\ell)}  (z,f,\varrho )$.

\noindent  A third observation is that  for $j= 3,4$ we have for $\mathbf{k}=\textbf{R}^{(\ell )} _j\circ \phi$
\begin{equation}  \label{eq:derivatives} \begin{aligned} &      \partial ^\mu _{z} \partial ^\nu _{{\overline{z}}}   \mathbf{k} _{| (0,0, \varrho ) } =0 \text{ for } |\mu |+|\nu  |\le \ell +1 \\&    \partial ^\mu _{z} \partial ^\nu _{{\overline{z}}} \nabla _f \mathbf{k} _{| (0,0, \varrho ) }=0 \text{ for } |\mu |+|\nu  |\le \ell  .
\end{aligned}
\end{equation}
By Lemma \ref{lem:ODEbis}  for $l=m$, $s=k$ and $r=k'$,
we have for $k\le k'-(2m+1) \textbf{d}$
\begin{equation}   \label{eq:key31}    \begin{aligned} &  \Pi _j(f)   \circ \phi  =  \Pi _j(f)   \circ \phi _0
+\resto ^ {0,2\ell +1} _{k,m}(\Pi (f),R),
\end{aligned}
\end{equation}
with $\phi _0=\phi _0^1$  and $ \phi _0^t$ the flow defined  as in Lemma 
\ref{lem:ODEbis} using the field $X_\chi ^{st}$.
Then we have
\begin{equation}   \label{eq:key32back}     \begin{aligned} &  \Pi _j(f)   \circ \phi  _0= \Pi _j(f)    +\int _0^1
  \left  \langle \Diamond _j
 (X_\chi ^{st} ) _{f} (\Pi (f), R \circ \phi _0^t) ,  f \circ \phi _0^t \right  \rangle dt
.
\end{aligned}
\end{equation}
By the definition of  $X_\chi ^{st}$ and    by   formulas \eqref{eq:backcoo} for  $\phi _0^t$,
which are simpler because there are no   phase
factors,   by $|\mu +\nu | =\ell $  the integrand in \eqref{eq:key32back}  is
\begin{equation*}        \begin{aligned} &  \left  (z  +
 \resto ^{0, \ell }_{k^{\prime \prime }, m} (t,\Pi (f) , R)\right ) ^\mu  \left (\overline{z}  +
 \resto ^{0, \ell }_{k^{\prime \prime }, m} (t,\Pi (f) , R)\right ) ^\nu    \\& \times  \left  \langle \Diamond _j
 B_{\mu , \nu } (\Pi (f)),  f+ \textbf{S} ^{0, \ell }_{k^{\prime \prime }, m}(t,\Pi (f) , R)   \right  \rangle \\& =  z   ^\mu  \overline{z}    ^\nu       \langle \Diamond _j
 B_{\mu , \nu } (\Pi (f)),  f     \rangle +\resto ^{0, 2\ell }_{k^{\prime \prime } , m}(t,\Pi (f) , R).
\end{aligned}
\end{equation*}
 Then for  $k\le  k^{\prime \prime } $ we have

\begin{equation}   \label{eq:key32}     \begin{aligned} &  \Pi _j(f)   \circ \phi  _0=  \Pi _j(f)    +
   \langle \Diamond _j
 (X_\chi ^{st} ) _{f}  ,  f \rangle
+\resto ^ {0, 2\ell  } _{k,m}(\Pi (f),R).
\end{aligned}
\end{equation}
 By $\ell \ge 2$ we have $2\ell  \ge \ell +2$ and so   $\resto ^ {0, 2\ell  } _{k,m} $ is an $\resto ^ {0, \ell +2} _{k,m}$.

 \noindent  By
$ {\psi} (\varrho)=O(|\varrho|^2)$ near 0,  we conclude that

\begin{equation} \label{eq:key-101}  \begin{aligned} &   {\psi} ( \Pi (f)  )  \circ \phi  =   {\psi} ( \Pi (f)  )+\widetilde{K} '+\resto  ^{1,\ell +2} _{k,m}(\Pi (f),R),
\end{aligned}
\end{equation}
with $\widetilde{K} '$ a polynomial as in \eqref{eq:tildeKrho} with $M_0=\ell ,$ with $ \widetilde{K} '(0,b,B)=0$
and
  $( \widehat{k} ,  \widehat{m}) =(k',m')$ satisfying. Notice that  it was to get the last equality, which follows from \eqref{eq:key32},  that
	we introduced the flow
$\phi _0^t$.

\noindent We now focus on $\mathbf{R}_2$.  We have by \eqref{eq:backcoo}

\begin{align} \label{eq:key-201}    &  \mathbf{R}_2  \circ \phi  = \langle \mathbf{B}_{2 } (   \Pi (f') ),   (f')  ^{   2} \rangle = \\&  \langle \mathbf{B}_{2 } \left (    \Pi (f)    +\resto ^ {0,\ell +1} _{k,m} (\Pi (f),R)\right ),   \left (e^{J\resto ^{0,\ell +1 } _{k^{\prime \prime },m'} (\Pi (f),R) \cdot \Diamond } ( f+ \textbf{S}^{0,\ell } _{k^{\prime \prime },m'}(\Pi (f),R) )\right )  ^{   2} \rangle
 .\nonumber
\end{align}
In our present set up the exponential $e^{J\resto ^{0,\ell +1 } _{k^{\prime \prime },m'}  \cdot \Diamond }$ cannot be moved to the
$\mathbf{B}_{2 }$ by a change of variables in the integral
as  in \cite{Cu1}.  Fortunately we know already that $H^{(\ell )}\circ \phi$  has the expansion of
Lemma \ref{lem:ExpH11} and that all we need to do is to compute some derivatives of $\mathbf{R}_2  \circ \phi $.

\noindent
Using  the expansion  in \eqref{eq:key-201}  and formula
\eqref{eq:symbol12},   for $i=0$   now, we set
\begin{equation}\label{eq:key-301} \begin{aligned} &
 \mathfrak{R}_2:= \langle \mathbf{B}_{2 } (    \Pi (f)    ),   (    f+ \textbf{S}^{i,\ell } _{k^{\prime \prime },m'}( \Pi (f)  ,R)  )  ^{   2} \rangle = \\&   \left \langle \mathbf{B}_{2 } (      \Pi (f)       ),   \left [    f+ \int _0^1 (X^{st}_\chi )_f \circ \phi ^t dt + \textbf{S}^{i,2\ell +1 } _{k^{\prime \prime },m'}(  \Pi (f)   ,R)  \right ]  ^{   2} \right \rangle  =   \\& \langle \mathbf{B}_{2 } (      \Pi (f)       ),       f ^{   2} \rangle+ 2  \int _0^1 \langle \mathbf{B}_{2 } (     \Pi (f)     ),   (X^{st}_\chi )_f \circ \phi ^t   \  f   \rangle  dt +  \resto^{i,2\ell  } _{k^{\prime \prime } ,m' }(  \Pi (f)   ,R)  . \end{aligned}
\end{equation}
We have that      $\mathbf{k}=
\mathbf{R}_2 \circ \phi -   \mathfrak{R}  _2$   is $C^{\ell +1}$ and satisfies \eqref{eq:derivatives}.
Hence the analysis of  $\mathbf{R}_2   \circ \phi$ reduces to that of  $\mathfrak{R}_2$.
By \eqref{eq:backcoo},  for $k\le k^{\prime\prime}$,  $m\le m'-1$ and $\ell >1$    we have
\begin{equation}\label{eq:key-401} \begin{aligned} & \int _0^1    X^{st}_\chi  \circ \phi ^t    dt  = X^{st}_\chi
+ \textbf{S}_{k^{\prime\prime},m'-1}^{0,2\ell -1}(  \Pi (f)   ,R)=  X^{st}_\chi
+ \textbf{S}_{k,m}^{0, \ell  +1}(  \Pi (f)   ,R)  . \end{aligned}
\end{equation}
This implies
\begin{equation}\label{eq:key-501} \begin{aligned} &
 \mathfrak{R}_2 =  \langle \mathbf{B}_{2 } (      \Pi (f)       ),       f ^{   2} \rangle+\widetilde{K} ^{\prime \prime}
  +
\resto _{k,m}^{0, \ell  +2}(\Pi (f),R) \quad , \\&  \widetilde{K}^{\prime \prime } :=2 \langle \mathbf{B}_{2 } (      \Pi (f)       ),       f (X^{st}_\chi  )_f  \rangle
 . \end{aligned}
\end{equation}
Then  $  \widetilde{K}^{\prime \prime } $ is  a polynomial like in \eqref{eq:tildeKrho} for the pair $(\widehat{k},\widehat{m})=(k',m')$ satisfying $ \widetilde{K} ^{\prime \prime }(0,b,B)=0$ by $B_2(\varrho )=0$ for $\varrho =0$.

\noindent By \eqref{eq:backcoo}  and for the pullback of the term
  $\resto _{k',m'+2}^{1, 2 }(      \Pi (f)  ,    f )$ in Lemma
\ref{lem:ExpH11}  we have   for $\varrho =\Pi (f )  $
  \begin{equation}  \label{eq:key-50110}
\begin{aligned} &\resto _{k',m'+2}^{1, 2 }(      \Pi (f')  ,    f ')= \resto _{k',m'+2}^{1, 2 }(     \varrho   ,    f ') \\&+
 \int _0^1 (\nabla _{\varrho }\resto _{k',m'+2}^{1, 2 })(    \varrho  +t
 \resto _{k^{\prime\prime},m'+2}^{0, \ell +1 }(      \varrho      ,    f  )  ,    f ')  \cdot \resto _{k^{\prime\prime},m'}^{0, \ell +1 }(       \varrho      ,    f  )   dt  \\& =   \resto _{k',m'+2}^{1, 2 }(     \varrho  ,    f ') +      \resto _{k,m}^{0, \ell +3 }  (      \varrho    ,   R )
   \end{aligned}
\end{equation}
for $k\le  k^{\prime\prime} - m \textbf{d}$ and $m\le m' $, by elementary analysis of the
second line.

\noindent  Applying again  \eqref{eq:backcoo}  we have
 \begin{equation}  \label{eq:key-50111}
\begin{aligned}       \resto _{k',m'+2}^{1, 2 }(       \varrho    ,    f ') &=  \resto _{k',m'+2}^{1, 2 }\left (       \varrho    ,    e^{J\resto ^{0, \ell +1 }_{k^{\prime \prime }, m' } ( \varrho    , R) \cdot \Diamond } \left ( f+ \textbf{S} ^{0, \ell }_{k^{\prime \prime }, m'}( \varrho    , R)\right  ) \right  ) \\&=  \resto _{k',m'+2}^{1, 2 }\left (       \varrho    ,   f+ \textbf{S} ^{0, \ell }_{k^{\prime \prime }, m'}( \varrho    , R) \right  )+    \resto _{k,m}^{1, \ell +2 }  (        \varrho       ,   R )
   \end{aligned}
\end{equation}
for $k\le  k^{\prime\prime} - m \textbf{d}$ and $m\le m'-1$. Next,  by Lemma \ref{lem:chi}, \eqref{eq:symbol12}
and by \eqref{eq:key-401}, we have
$  \textbf{S} ^{0, \ell }_{k^{\prime \prime }, m'}( \varrho    , R) =( X^{st}_\chi )_f   +\textbf{S}_{k^{\prime\prime},m'-1}^{0, \ell  +1}(  \varrho   ,R) $
and
\begin{equation*}  \begin{aligned} &  \resto _{k',m'+2}^{1, 2 }\left (     \varrho   ,       f+( X^{st}_\chi )_f   +\textbf{S}_{k^{\prime\prime},m }^{0, \ell  +1}(  \varrho   ,R) \right )   =\resto _{k',m'+2}^{1, 2 }(     \varrho   ,       f ) + \\&  \int _0^1 \langle \nabla _R \resto _{k',m'+2}^{1, 2 }\left (
   \varrho   ,       f+t( X^{st}_\chi )_f   +t\textbf{S}_{k^{\prime\prime},m }^{0, \ell  +1}(  \varrho   ,R)
	\right ) , ( X^{st}_\chi )_f   + \textbf{S}_{k^{\prime\prime},m }^{0, \ell  +1}
(  \varrho   ,R) \rangle   dt \\& =\resto _{k',m'+2}^{1, 2 }(     \varrho   ,       f ) +\langle \nabla _f \resto _{k',m'+2}^{1, 2 }(     \varrho   ,       f   ) , ( X^{st}_\chi )_f    \rangle   +\resto _{k ,m }^{1, \ell +2 }(     \varrho   ,       R )
   \end{aligned}
\end{equation*}
where we have used $\ell \ge 2$, $k\le k^{\prime\prime} \le k^{\prime }$ and $m\le m' -1$.
Notice that   we have  that $\resto _{k',m' +2}^{1, 2 }(     \varrho   ,       f )$ is an $\resto _{k ,m+2}^{1, 2 }(     \varrho   ,       f )$.
Finally we have
\begin{equation}\label{eq:key-5013} \begin{aligned} &
   \langle \nabla _f \resto _{k',m'+2}^{1, 2 }(     \varrho   ,       f   ) , ( X^{st}_\chi )_f    \rangle    =\widetilde{K} ^{\prime \prime \prime }+\overline{\textbf{R}}_2\quad ,
	\\& \widetilde{K} ^{\prime \prime \prime }:= \langle \nabla _f ^2\resto _{k',m'+2}^{1, 2 }(     \varrho   ,      0   ) f,    ( X^{st}_\chi )_f    \rangle   ,
   \end{aligned}
\end{equation}
with $\overline{\textbf{R}}_2 $ a term we can absorb in $\widehat{\textbf{R}}_2 $ and with $ \widetilde{K} ^{\prime \prime \prime }$ like in \eqref{eq:tildeKrho} for the pair $(\widehat{k},\widehat{m})=(k',m')$ satisfying $ \widetilde{K} ^{\prime \prime \prime }(0,b,B)=0$.

\noindent We  set
\begin{equation} \label{eq:key-601}  \begin{aligned} &      \textbf{R}^{(\ell )} _0 +\textbf{R}^{(\ell )} _1   =Z  '+ K +  \textbf{R}  _{01}  \  ,
\end{aligned}
\end{equation}
where:  $Z  '$ is the sum of  the monomials in normal form of degree  $\le \ell+1$;
$K$, which is like in \eqref{eq:Krho},  is the sum of the the monomials     of degree
equal to $ \ell+1$  not in normal form;
   $\textbf{R}  _{01}$  is the sum  of the monomials of  degree  $> \ell+1$. By induction there are no monomials not in normal form of degree $\le \ell $  so that each of  the monomials of the lhs of  \eqref{eq:key-601}
	go into exactly one of the three terms of the rhs.

 \noindent  We define  ${Z}^{\prime  \prime  } $ and $\widetilde{K}$ by setting

\begin{equation}  \label{eq:key-701} \begin{aligned} &     \widetilde{K}^{\prime  }+\widetilde{K}^{\prime \prime }+\widetilde{K} ^{\prime \prime \prime }=  {Z}^{\prime  \prime  } +\widetilde{K},
\end{aligned}
\end{equation}
collecting in   ${Z}^{\prime  \prime  } $ all monomials of the lhs in normal form (all of degree  $  \ell+1$) and      in   $\widetilde{K} $ all monomials of the lhs not in normal form.
Here  $\widetilde{K} $  is  like in \eqref{eq:tildeKrho} for $(\widehat{k},\widehat{m})=(k',m')$ with  $ \widetilde{K}  (0,b,B)=0$.

      \noindent Applying Lemma \ref{lem:NLhom1} for $(\widehat{k},\widehat{m})=(k',m')$ we
can choose $\chi $ such that for $Z=Z  ' +{Z}^{\prime  \prime  }  $ we have
\begin{equation}\label{eq:key1}
 \{  H_2^{(\ell )}  , \chi \}  ^{st}+Z+ K+\widetilde{K}=0.
\end{equation}
Then $H^{(\ell +1)}:=H^{(\ell )}\circ \phi $ satisfies the conclusions of
  Theorem \ref{th:main} for $\ell+1$.

\bigskip

{\it The step $\ell +1=2$.}
Set $H^{(1 )}=K\circ \mathfrak{F}_1$.  We are seeking a transformation
$\phi$ as in the previous  part such that  $H^{(2)}:=H^{(1 )}\circ \phi $ has term  $\mathbf{R} _{-1}^{(2)}=0 $ in its expansion
in Lemma \ref{lem:ExpH11}.
The argument is similar to the previous one, but this time $\chi $  has degree $\ell +1$ with
  $\ell =1$.
So  the  steps in the previous argument where we   exploited  $\ell \ge 2 $   need to
be reframed.
\noindent  We know  that $H^{(1 )}$ satisfies Lemma \ref{lem:ExpH11} for $L=1 $
for some pair that we denote by $(k',m')$  rather than $(k,m)$.

\noindent The proof  of \eqref{eq:key-11} is different from the previous one. By \eqref{eq:ODE1bis}
 we have for some $(k,m)$ appropriately smaller than  $(k',m')$

\begin{equation}   \label{eq:keyfirst}    \begin{aligned} &
 \{  H_2 ^{(1 )}   , \chi \} ^{st}\circ \phi ^t=  \{  H_2 ^{(1 )}   , \chi \} ^{st}\circ \phi ^t _0+  \resto ^{0,4}_{k,m}(\Pi (f),R).
\end{aligned}
\end{equation}
 The following  linear transformation

\begin{equation*} \label{eq:linop}   \begin{aligned} &   (Z,\overline{Z},F) \to  \begin{pmatrix}    \im \nu _j b_{\mu \nu}(\Pi (f)) \frac{Z^\mu \overline{Z}^\nu  }{\overline{Z}_j} +  \im \nu _j  \frac{Z^\mu \overline{Z}^\nu  }{\overline{Z}_j} \langle J^{-1}B_{\mu \nu}(\Pi (f)),F\rangle
   \\  -\im \mu _j b_{\mu \nu}(\Pi (f)) \frac{ {Z}^\mu \overline{{Z}}^\nu  }{ {Z}_j} -  \im \mu _j  \frac{Z^\mu \overline{Z}^\nu  }{ {Z}_j} \langle J^{-1}B_{\mu \nu}(\Pi (f)),F\rangle
\\ B_{\mu \nu}(\Pi (f))  {Z}^\mu \overline{{Z}}^\nu
 \end{pmatrix}
\end{aligned}
\end{equation*}
  depends linearly on $(b (\varrho),B (\rho))$, for $\varrho =\Pi (f) $.
	Then
	\begin{equation} \label{eq:linop1}   \begin{aligned} &    z_j \circ \phi ^t _0=z_j+
 a  _j(t,b,B)\cdot z +   b _j(t,b,B)\cdot \overline{z} + \sum _{\mu \nu} c _{ j\mu  \nu} (t,b,B)
\langle J^{-1}B_{\mu \nu} ,f\rangle
\end{aligned}
\end{equation}
	for $a _j,b  _j\in C^\infty ([0,1]\times X_{k'}, \C^{\textbf{n}})$
	with $|a  _j| + |b    _j| \le C \| (b,B)\| _{  X_{k'}}$
	and
	$c _{ j\mu  \nu} \in C^\infty ([0,1]\times X_{k'}, \C )$.
	Similarly
	\begin{equation} \label{eq:linop2}   \begin{aligned} &   f \circ \phi ^t _0=f+
 \textbf{a}   (t,b,B)\cdot z +  \textbf{ b}  (t,b,B)\cdot \overline{z} + \sum _{\mu \nu} \textbf{c}  _{ \mu  \nu} (t,b,B)
\langle J^{-1}B_{\mu \nu} ,f\rangle
\end{aligned}
\end{equation}
	with  $\textbf{a}  ,\textbf{b}  \in C^\infty ([0,1]\times X_{k'}, \Sigma _{k'}^{\textbf{n}})$
	with $\| \textbf{a}    \| _{\Sigma _{k'}^{\textbf{n}}} + \|\textbf{b}     \| _{\Sigma _{k'}^{\textbf{n}}} \le C \| (b,B)\| _{  X_{k'}}$
	and
	$\textbf{c} _{  \mu  \nu} \in C^\infty ([0,1]\times X_{k'}, \Sigma _{k'} )$.  These coefficients
	satisfy appropriate symmetries that ensure $\overline{f \circ \phi ^t _0}=f \circ \phi ^t _0$.

  \noindent We  have
    \begin{equation}  \label{eq:key-1310}   \begin{aligned} & \{  H_2 ^{(1 )}   , \chi \} ^{st}
   \circ \phi ^t _0= \{  H_2 ^{(1 )}   , \chi \} ^{st}
  (\Pi (f),  R\circ \phi _0^t) + \resto ^{1,4}_{k,m} (t,\Pi (f), R).
\end{aligned}
\end{equation}
To compute $ \{  H_2 ^{(1 )}   , \chi \} ^{st}
  (\Pi (f),  R\circ \phi _0^t)$ we  replace the
  $R$ in \eqref{eq:key0}  with $R\circ \phi _0^t $. The coordinates  of the latter can be expressed in terms
	of $R$ by \eqref{eq:linop1}--\eqref{eq:linop2}.
  When we substitute $(z,f)$ in \eqref{eq:key0}
	using
	\eqref{eq:linop1}--\eqref{eq:linop2},  
	by an elementary
computation we obtain
\begin{equation*}     \begin{aligned} & \{  H_2 ^{(1 )}   , \chi \} ^{st}
  (\varrho,  R\circ \phi _0^t)=  \{  H_2 ^{(1)}    , \chi \} ^{st} (\varrho,  R )   \\& + \sum _{\substack{|\mu +\nu |=1\\  |\mu '+\nu '|=1}}
	a_{\mu \nu}^{\mu ' \nu '} (t,\varrho ,b(\varrho ),B(\varrho )) z^{\mu} \overline{z}^{\nu}\langle  \mathcal{H} B_{\mu \nu}(\varrho ), f
\rangle   +A^t + \underline{{\mathbf{R}}}^t.
\end{aligned}
\end{equation*}

Here:
\begin{itemize}

\item    $a_{\mu \nu}^{\mu ' \nu '} (t,\varrho ,b ,B ) \in C^{m'}$
with  $a_{\mu \nu}^{\mu ' \nu '} (t,0 ,0 ,0 ) =0$;

\item   we have \begin{equation*}     \begin{aligned} &
A^t=\sum _{|\mu +\nu |=2}
	\alpha _{\mu \nu}  (t,\varrho ,b(\varrho ),B(\varrho ))z^{\mu} \overline{z}^{\nu} \\& +\sum _{l=0}^1  \sum _{j=1}^{n_0}
	\sum _{ |\mu +\nu |=1 }
	  z^{\mu} \overline{z}^{\nu} \langle \Diamond _j^l  A^l_{\mu   \nu  }(t, \varrho , b(\varrho ),B(\varrho )), f
\rangle   ,
\end{aligned}
\end{equation*}
$\alpha _{\mu \nu}(t,\varrho ,b,B)$ and $A^l_{\mu   \nu  }(t,\varrho ,b,B)$ are  $C^{m'}$
with
 for $i=2$
\begin{equation}
	\label{eq:corr1}     |\alpha _{\mu \nu} (t,\varrho ,b,B) | + \| A^l_{\mu   \nu  }(t,\varrho ,b,B) \| _{\Sigma _{k'}}\le C \| ( b,B)\| _{X_{k'}}^i;
\end{equation}
\item  $\underline{\mathbf{R}} ^t( \varrho ,z,f)$ is $C^m$ in $(t,\varrho ,z,f)\in \R ^{n_0+1} \times \C ^{\mathbf{n}} \times \Sigma _{-k}$ with $(\varrho ,z,f)$     near $(0 ,0,0)$, with    for $i=2$
\begin{equation}
	\label{eq:corr2}     | \underline{{\mathbf{R}}}^t |  \le C  \| ( b,B)\| _{X_{k'}}^2 \| f\| _{\Sigma _{-k}}^2  .
\end{equation}
\end{itemize}
Then, in the notation of  Lemma \ref{lem:NLhom2}
    \begin{equation}  \label{eq:key-110}   \begin{aligned} & \int _0^1\{  H_2 ^{(1 )}   , \chi \} ^{st}
 \circ   \phi ^t _0 dt = \{  H_2 ^{(1)}    , \chi \} ^{\widetilde{st}}    + A
 +\underline{{\mathbf{R}}} +\resto ^{1,4}_{k,m}(\Pi (R) , R),
\end{aligned}
\end{equation}
 with $A=\int _0^1 A^tdt  $ and $\underline{{\mathbf{R}}}=\int _0^1 \underline{{\mathbf{R}}}^tdt  $  are like $A^1$ and $\underline{{\mathbf{R}}}^1$.
  Then, using also \eqref{eq:keyfirst},
	we get the following analogue of \eqref{eq:key001}: \begin{equation} \label{eq:key-001}\begin{aligned} &   H_2^{(1 )}  \circ
\phi  = H_2^{(1 )}  +    \{  H_2^{(1 )}  , \chi \}  ^{\widetilde{st}} + A+ \underline{{\mathbf{R}}}  + \resto ^ {0,4}_{k,m}(\Pi (f),R).
\end{aligned}
\end{equation}
\eqref{eq:key31} remains true also for $\ell =1$.  We consider  \eqref{eq:key32back} and expand
 \begin{equation*}  \begin{aligned} &
\langle \Diamond _j
 (X_\chi ^{st} ) _{f} (\Pi (f), R \circ \phi _0^t) ,  f \circ \phi _0^t    \rangle    =
\langle \Diamond _j
 (X_\chi ^{st} ) _{f} (\Pi (f), R)   ,  f     \rangle + A^t + \textbf{R} ^t,
\end{aligned}
\end{equation*}
with $A^t$     and
  $\textbf{R} ^t$    like the previous ones  but such that  \eqref{eq:corr1}--\eqref{eq:corr2}
	hold for $i=1$.  This yields
\begin{equation}   \label{eq:key320}     \begin{aligned} &  \Pi _j(f)   \circ \phi  _0= \Pi _j(f)
+A'+\underline{{\mathbf{R}}}'.
\end{aligned}
\end{equation}
Here        $ \underline{{\mathbf{R}}}'$  is like  $ \underline{{\mathbf{R}}} ^1$  such that  \eqref{eq:corr2}
	holds for $i=1$.
	$A'$  is like $A^1$  such that \eqref{eq:corr1}
	holds for $i=1$.
	
\noindent By $ {\psi} (\varrho)=O(|\varrho|^2)$ near 0 and \eqref{eq:key31}  we get the first equality in

\begin{equation} \label{eq:key101}  \begin{aligned} &   {\psi} ( \Pi (f)  )  \circ \phi  = {\psi} ( \Pi (f)  )  \circ \phi _0 +
 \resto ^ {1, 3} _{k,m}(\Pi (f),R)\\&=
  {\psi} ( \Pi (f)  )+\widetilde{K} '+\resto  ^{1,2} _{k',m'}(\Pi (f),f)+\resto ^ {1, 3} _{k,m}(\Pi (f),R),
\end{aligned}
\end{equation}
    where $\widetilde{K} '=\resto  ^{1,2} _{k',m'}(\Pi (f),R)$ is a polynomial in $R$ as in \eqref{eq:tildeKrho} with $\widetilde{K} ' (0,b,B)=0$. The second line in  \eqref{eq:key101} follows by  $ {\psi} (\varrho)=O(|\varrho|^2)$, by the fact that  $ {\psi} (\varrho)$ is smooth 
		and  by \eqref{eq:key320}. Notice that by choosing $m\le m' -2$ we have  $\resto  ^{1,2} _{k',m'}(\Pi (f),f)=\resto  ^{1,2} _{k ,m+2}(\Pi (f),f).$
		
 \noindent The discussion of $\textbf{R}\circ \phi $ is similar to the previous one   after \eqref{eq:key-201} . This time, though,
   by  \eqref{eq:ODE1bis}   we write 

\begin{equation}\label{eq:key401} \begin{aligned} & \int _0^1    X^{st}_\chi  \circ \phi ^t    dt  =
\int _0^1    X^{st}_\chi  \circ \phi ^t_0    dt + \textbf{S}_{k,m}^{0,3} (\Pi (f),R)  . \end{aligned}
\end{equation}
By \eqref{eq:linop1}--\eqref{eq:linop2}  we get

\begin{equation}\label{eq:key402} \begin{aligned} & \int _0^1    X^{st}_\chi  \circ \phi ^t _0   dt
=
X^{st}_\chi    +  \mathbf{A}   \text{ in $\Ph ^{k'}$} , \end{aligned}
\end{equation}
with $(z,f)\to  \mathbf{A}(\varrho ,z,f)$ linear, with $C^{m'} $
dependence in $\varrho$ and with
\begin{equation}\label{eq:key403} \begin{aligned} &  \|  \mathbf{A} (\varrho , z, f)  \|_ {  \Ph ^{k'}}\le C  \| (b (\varrho ),B(\varrho )) \| _{X_{k'}}  (|z|+\| f \| _{\Sigma _{-k'}}). \end{aligned}
\end{equation}
This yields, for $\mathfrak{R}_2$ defined as in \eqref{eq:key-301},
\begin{equation*}  \begin{aligned} &
 \mathfrak{R}_2 =    \left \langle \mathbf{B}_{2 } (      \Pi (f)       ),   \left [    f+ \int _0^1 (X^{st}_\chi )_f \circ \phi ^t_0 dt   \right ]  ^{   2} \right \rangle +\resto ^{1,3 } _{k ,m }(  \Pi (f)   ,R) =\\&  \langle \mathbf{B}_{2 } (      \Pi (f)       ),       f ^{   2} \rangle+ 2\langle \mathbf{B}_{2 } (      \Pi (f)       ), f  \mathbf{A}\rangle +  \langle \mathbf{B}_{2 } (      \Pi (f)       ),    \mathbf{A}^2\rangle  +\resto ^{1,3 } _{k ,m }(  \Pi (f)   ,R)  , \end{aligned}
\end{equation*}
where we have used $ \mathbf{B}_{2 } (      0       )=0$ for the reminder.

\noindent We have
\begin{equation*}  \begin{aligned} &
   2\langle \mathbf{B}_{2 } (      \Pi (f)       ), f  \mathbf{A}\rangle +  \langle \mathbf{B}_{2 } (      \Pi (f)       ),    \mathbf{A}^2\rangle  = \widetilde{K} ^{\prime \prime}
   +  \underline{\mathbf{R}}^{\prime \prime}  , \end{aligned}
\end{equation*}
with  $\underline{{\mathbf{R}}}^{\prime \prime}$ like $\underline{{\mathbf{R}}}$    and  with
$\widetilde{K} ^{\prime \prime}$ like \eqref{eq:tildeKrho}
with $\widetilde{K} ^{\prime \prime}(0,b,B)=0$, by $\mathbf{B}_{2 } (      0       )=0$, and with $(\widehat{k},\widehat{m})=(k',m')$.
Summing up, we have
\begin{equation}\label{eq:key501} \begin{aligned} &
 \mathfrak{R}_2 =  \langle \mathbf{B}_{2 } (      \Pi (f)       ),       f ^{   2} \rangle+\widetilde{K} ^{\prime \prime}
  +   \underline{\mathbf{R}}^{\prime \prime}+\resto^{1,3  } _{k ,m }(  \Pi (f)   ,R).
   \end{aligned}
\end{equation}
Notice that the reduction of  $\mathbf{R}_2\circ \phi $ to $ \mathfrak{R}_2$  continues to hold also for $\ell= 1$.

\noindent We consider $\resto _{k',m'+2}^{1, 2 }\circ \phi$ from  the $ \resto _{k',m'+2}^{1, 2 }$ term in the expansion of $\mathbf{R}$ in Lemma
\ref{lem:ExpH11}.
Then, by  \eqref{eq:key-50110}   and by \eqref{eq:key401}--\eqref{eq:key402}, for $\varrho =\Pi (f )$ we have
 \begin{equation*}  \begin{aligned} &
\resto _{k',m'+2}^{1, 2 }(      \Pi (f')  ,    f ')= \resto _{k',m'+2}^{1, 2 }(      \varrho  ,    f   +( X^{st}_\chi )_f +\mathbf{A} + \textbf{S}_{k,m}^{0,3} ) + \resto _{k,m}^{0, 4 }(      \varrho  ,    R )
 . \end{aligned}
\end{equation*}
 The first term in the rhs can be expanded for $\varrho =\Pi (f )$  as
  \begin{equation*}  \begin{aligned} &
     \resto _{k',m'+2}^{1, 2 }(     \varrho   ,       f   +( X^{st}_\chi )_f +\mathbf{A} )     +\resto _{k ,m }^{1, 4 }(     \varrho   ,       R ).
   \end{aligned}
\end{equation*}
We have for $\varrho =\Pi (f )$
\begin{equation*}  \begin{aligned} &
     \resto _{k',m'+2}^{1, 2 }(     \varrho   ,       f   +( X^{st}_\chi )_f +\mathbf{A} )   = \mathfrak{B}_2(\varrho ) ( f   +( X^{st}_\chi )_f +\mathbf{A} )^2+ \resto _{k ,m }^{1, 3 }(     \varrho   ,       R ),
   \end{aligned}
\end{equation*}
with $\mathfrak{B}_2(\varrho )$ a $C^{m'}$ function with values
in $B ^2(\Sigma_{-k'},\Sigma_{ k'}) $ with $\mathfrak{B}_2(0 )=0$.
 Considering the binomial expansion we get for $\varrho =\Pi (f )$

\begin{equation*}  \begin{aligned} &
\resto _{k',m'+2}^{1, 2 }(      \Pi (f')  ,    f ') = \mathfrak{B}_2(\varrho )   f      ^2 +\widetilde{K}^{\prime\prime\prime} +\underline{\mathbf{R}} ^{\prime\prime\prime}+\resto _{k ,m }^{0, 3 }(     \varrho   ,       R ) ,
   \end{aligned}
\end{equation*}
  with $\underline{\mathbf{R}} ^{\prime\prime\prime}$  like $\underline{\mathbf{R}}  $ and with $\widetilde{K}^{\prime\prime\prime}$ like
  \eqref{eq:tildeKrho} with $\widetilde{K}^{\prime\prime\prime}(0,b,B)=0$ and $(\widehat{k},\widehat{m})=(k',m')$.

\noindent We now set  $K= \textbf{R}^{(1 )} _{-1}$  and  with the $A$ of  \eqref{eq:key-110} we write

\begin{equation}  \label{eq:key701} \begin{aligned} &      \widetilde{K}^{\prime  }+\widetilde{K}^{\prime \prime }+\widetilde{K} ^{\prime \prime \prime }+A  ={Z}^{\prime  \prime  } +\widetilde{K}  ,
\end{aligned}
\end{equation}
where in ${Z}^{\prime  \prime  }$ we collect the null terms of the lhs
and in $\widetilde{K}$ the other terms.
Now we have   $K(0 )=0$,    $\widetilde{K}( 0 ,0,0  )=0$  and $\nabla _{b,B}\widetilde{K}( 0 ,0,0  )=0$. By   Lemma  \ref{lem:NLhom2} for $(\widehat{k},\widehat{m})=(k',m')$ we
can choose $\chi $ such that for   we have
\begin{equation} \label{eq:key702}
 \{  H_2^{(\ell )}  , \chi \}  ^{\widetilde{st}}+Z^{\prime \prime }+ K+\widetilde{K}=0.
\end{equation}
Then $H^{(2)}:=H^{(1 )}\circ \phi $ satisfies the conclusions of
  Theorem \ref{th:main} for $\ell =2$.

\qed

Summing up, we have proved the following result, whose proof we sketch now.

\begin{theorem}\label{theorem-1.1}  For fixed $p_0\in \mathcal{O}$   and for sufficiently large $l\in \N$,
there are a fixed $k\in \N$, an  $\epsilon >0$,  an $1\ll s' \ll l$  and a   $1\ll k  \ll k'$ such that for solutions $\widehat{U}(t)$ to
\eqref{eq:NLSvectorial}  with $\Pi (U) =p_0$  with $|\Pi (\widehat{R} (t)) |+\| \widehat{R} (t) \| _{\Sigma _{-k}}<\epsilon $ and $\widehat{R} (t) \in \Sigma _{l}$, there
exists a $C^0$ map $\Phi :\U ^l_{\epsilon  , k}  \to \U ^{s'}_{\epsilon  ' , k'} $ such that
 \begin{equation}  \label{eq:1.11}\begin{aligned} &
R :=\Phi    _{R} (  \Pi (\widehat{R} ), \widehat{R} ) = e^{Jq(  \Pi (\widehat{R} ), \widehat{R} )\cdot \Diamond } ( \widehat{R} + \textbf{S} ( \Pi (\widehat{R} ), \widehat{R}     ))
,
\end{aligned}   \end{equation}
\begin{equation} \label{eq:1.12} \begin{aligned} \text{with   }  &
     \textbf{S}  \in C^2((-2,2)
\times B_{\R ^{ n_0}}\times   B_{\Sigma _{-k}} , \Sigma _{s'}
 )    \\&  {q}  \in C^2((-2,2)
\times B_{\R ^{ n_0}} \times   B_{\Sigma _{-k}}, \R ^{ n_0}
 )
\end{aligned}   \end{equation}
such that $    \| \textbf{S} (  \Pi (\widehat{R} )  ,\widehat{R}  )
\| _{\Sigma _{ s'} }\le C   \epsilon\| \widehat{R} \| _{\Sigma _{ -k}} $
and such that  splitting $R(t) $ in spectral coordinates $(z(t), f(t))$ the latter satisfy

\begin{equation} \label{eq:HamSyst}\begin{aligned} & \dot z_j =\im  \partial _{\overline{z}_j}H\ ,   \quad \dot f=J\nabla _fH
\end{aligned}\end{equation}
where $H$ is a given function satisfying the properties of $H^{(2\textbf{N}+1)}$ in Theorem \ref{th:main}.
\end{theorem}
  \proof     Since
	 in Lemma \ref{lem:fred12} we can pick arbitrary $n$,
	we see by   the proof of  Theorem \ref{th:main}
 that we can suppose  that the $2\textbf{N}+1$ transformations $\phi _\ell$ are defined by flows
\eqref{eq:ODE}  with pair $(r,M)$ with $r$ and $M$ as large as needed.

\noindent 	Starting with an appropriate  $\U ^{s}_{\varepsilon _0, \kappa _0}$, we know that there  is a map
	$ \mathfrak{F}  : \U ^{s'}_{\varepsilon _1, \kappa '} \to \U ^{s }_{\varepsilon _0, \kappa _0}$  as regular as needed
	which satisfies the conclusions of  Theorem \ref{th:main}.
	In particular here
	we have $s'\gg s$ and $ 1\ll  \kappa ' \ll  \kappa _0$ and in   $\U ^{s'}_{\varepsilon _1, \kappa '}$
	 we get the system  \eqref{eq:HamSyst} by pulling back  the system  which exists in  $\U ^{s }_{\varepsilon _0, \kappa _0}$.

\noindent We choose now  $l\gg s'$, $1 \ll k\ll  \kappa '$ and
	sufficiently small  $\epsilon  $ and $ \delta $ with
	$\U ^l_{\delta  , k}  \subset
	\U ^{s }_{\varepsilon _0, \kappa _0}$ and $\U ^l_{\epsilon  , k}  \subset
	\U ^{s'}_{\varepsilon _1, \kappa '} $.  Here $l$ and $ \kappa '$ can be as large as we want, thanks to
	our freedom to choose  $(r,M)$.

	\noindent  By choosing $\delta $ small we can assume
	   $\U ^l_{\delta  , k} \subset \mathfrak{F}  ( \U ^{s'}_{\varepsilon _1, \kappa '})
	 $.  This follows from \eqref{eq:main1} which implies $\mathfrak{F} ^{-1}(\U ^l_{\delta   , k})
	\subset \U ^l_{\epsilon  , k} $. Finally   we set $\Phi =  \mathfrak{F} ^{-1}$ where
	$\mathfrak{F} ^{-1} : \U ^l_{\delta   , k}  \to   \U ^{s'}_{\varepsilon _1, \kappa '}$.
	
	\noindent Formula \eqref{eq:1.11}  and the information on $\textbf{S}$ has been   proved in the course of the proof  of  Lemma \ref{lem:ODE1}.  The  information on the phase function $q$ can be proved by a similar
	induction argument, which we skip here.
	
	\qed

\begin{remark}
\label{rem:bamb1}   The paper \cite{bambusi}
   highlights  in the Introduction  and states in  Theorem 2.2,
	that it is able  to  treat   all      solutions of the NLS  near ground states    in $H^1 $. But
   in fact, in    \cite{bambusi} there is no explicit proof of this.  While
	\cite{bambusi} does not state  the regularity properties
	of the maps in   Theorem 3.21 and   Theorem 5.2   \cite{bambusi},
		from the context they appear to be just continuous. Even if we assume that they are   \textit{almost smooth} transformations
		(but see Remark \ref{rem:formal1} above), nonetheless an  explanation
	is required	    on   why they preserve the structure needed to make sense of the NLS.
		But while pullbacks of the Hamiltonian
are    analyzed, the pullbacks of differential forms and the making sense of them,
are not discussed    in  \cite{bambusi}.
  For example, there is no explicit discussion on why  $\mathfrak{F} ^{t*} \Omega _t$   makes sense
	 in formula (3.42)  \cite{bambusi},  i.e. \eqref{eq:fdarboux} here.
\end{remark}

\begin{remark}
\label{rem:formal2} In the 2nd version of  \cite{bambusi} there is an incorrect effective
Hamiltonian.  If we  use  the correct definition of  the  symbols $ \mathbf{{S}}^{i, j} $ which we give above,
the functions $\Phi _{\mu \nu }$ used in the normal form expansion   in \cite{bambusi} are in $\mathcal{W}^j$ for some large $j$, rather than in  $ \cap _{j\ge 0}\mathcal{W}^j$.   In  pp. 25--27 in the 2nd version  of \cite{bambusi}, the
  $\mathcal{W}^j$'s  are defined using the  classical
pair  of operators $L_\pm $, see \cite{W2}, and are   closed subspaces of  $ H^{j-1}(\R^3)$ of finite codimension.
	This last fact seems to be unnoticed in   \cite{bambusi} and leads to the breakdown of the proof in the 2nd version  of \cite{bambusi}, as we explain below.
The space $\mathcal{W}^2$, for example,
 is defined   by  first considering  $ \langle L_+u, u\rangle  $
for
    $u\in \ker ^\perp L_-\cap \ker ^\perp L_+\subset L^2$. Notice that  $ \langle L_+u, u\rangle  \ge 0 ,$
				see  Prop. 2.7  \cite{W2} or Lemma 11.12 \cite{RSS}. Proceeding like in  Lemma 11.13 \cite{RSS}
		it can be shown that for  $  u\in \ker ^\perp L_-\cap \ker ^\perp L_+\subset L^2$ with $u\neq 0$ we have
		$\| u \| _L^2 :=\langle L_+u, u\rangle   >0.$ Then consider  	
  the completion of $\ker ^\perp L_- \cap \ker ^\perp L_+\cap C_0^\infty$ by the norm $\| u \| _L $.
	This  completion is  exactly
  $\ker ^\perp L_-\cap \ker ^\perp L_+\cap H^1(\R^3)$. Then  $\mathcal{W}^2$ is
   a closed subspace of finite codimension of   the latter space. Specifically,
   $\mathcal{W}^2$ is in the continuous spectrum part
   in the spectral decomposition of  the operator $L_-L_+$, which is selfadjoint for
	$\langle  u, v \rangle _L := \langle  L^{-1}_- u, v \rangle $ in $\ker ^\perp L_-$.
	Notice that, under hypotheses analogous to (L1)--(L6) in Sect. \ref{sec:speccoo},  $L_-L_+$ has finitely many eigenvalues and its eigenfunctions  are Schwartz functions.
	Likewise,  also the  other $\mathcal{W}^j $'s  are closed subspaces of  $ H^{j-1}(\R^3)$ of finite codimension.  Later in the 2nd version of \cite{bambusi}, at
		p.41,  the Strichartz estimates hinge on  the false  inclusion
of $  \mathcal{W}^j $, or of $ \mathcal{W}^\infty $, in $L^{ \frac{6}5 }(\R ^3, \C )$. 
Additional mistakes appear in the justification of  the     Fermi Golden rule. While formulas 
$R^\pm _{L_0}(\rho )\Phi $ in (St.2)--(St.3)  on p. 38 of the 2nd version   make sense because  $ \Phi \in H^{k ,s }$
for $s>0$ appropriate,  analogous formulas  $R^\pm _{B}(\rho )\Phi $  in  (6.50) and elsewhere in Sect. 6.2,
are undefined when we know only that $\Phi  \in \mathcal{W}^\infty $. In fact even  $R^\pm _{-\Delta }(\rho )\Phi $
is   undefined for $\rho \ge 0$  for such $\Phi$'s. So in particular, in the    2nd version of \cite{bambusi},
the discussion of the    Fermi Golden rule is purely formal.

   The above ones are not   simple  oversights.  Rather, they stem from the fact that,
   in the 2nd version of \cite{bambusi},  the
   homological equations  are solved    only in these   $\mathcal{W}^j$'s, while it is unclear if they can be solved  in spaces
    with spacial weights like the $H^{k,n}$ or the $\Sigma _n$ for $n>0$, as we remarked in an early version of \cite{Cu1}. 
   The 3rd version of   \cite{bambusi}
	credits  our remark for having stimulated
	changes  in this part of the paper.   These  changes are   classified   in  the 3rd version of \cite{bambusi}
	as   plain   simplifications.
	This might leave the   wrong  impression that the proof  in  the 2nd version of \cite{bambusi},
	although more complicated than in the 3rd version, is still correct.    

\end{remark}

\section{The NLS and the Nonlinear Dirac Equation}
\label{sec:examples}
We give a sketchy discussion of few examples.

\textbf{The Nonlinear Schr\"odinger equation.}  We consider the equation

\begin{equation*}
 \im U_{t }=-\Delta U +  2B'(|U|^2) U  \ .
\end{equation*}
 Here $N=1$, ${\mathcal D}=-\Delta $, $|\ | _1=|\ |$, $J=\begin{pmatrix}  0 &
1  \\
-1 & 0
 \end{pmatrix}$.
There are four invariants:
 \begin{equation*}
 \begin{aligned} &  Q(U)=\Pi _4(U)=\frac 12\langle U , U \rangle   \text{ and }
\Pi _j(U)=\frac 12\langle U , J \frac{\partial }{\partial x_j} U \rangle  \text{ for $j\le 3.$}
\end{aligned}
\end{equation*}
For  fixed $v\in \R ^3$ we have
\begin{equation*}
 \begin{aligned} &  Q(e^{-\frac 12 J v\cdot x}U)=Q(U)  \ , \
\Pi _j(e^{-\frac 12 J v\cdot x}U)=\Pi _j( U)  -\frac{v_j}2 Q(U)   \text{ for $j\le 3 $ and} \\&
E(e^{-\frac 12 J v\cdot x}U) = E(U)-\sum _{j=1}^{3} v_j\Pi _j( U) + \frac{v^2}{2}Q(U).
\end{aligned}
\end{equation*}
There is    well established theory guaranteeing under appropriate hypotheses existence  of   open sets $\mathcal{O}\subseteq \R ^+$
and
  $(\phi _ {\omega },0)\in C^\infty (\mathcal{O},{ \mathcal S}(\R ^3, \R ^2))$  such that
\begin{equation*}
  \Delta \phi _ {\omega } -\omega \phi _ {\omega }+2B'( \phi _ {\omega }  ^2)\phi _ {\omega }=0\quad\text{for $x\in \R^3$}.
\end{equation*}
More precisely it is possible to prove exponential decay to 0 of  $\phi _ {\omega }(x)$ as $x\to \infty$.

  \noindent For $v\in \R ^3$
arbitrary
we get
$\Phi _p (x) =e^{-\frac 12 J v\cdot x} (\phi _ {\omega } (x),0)$ where $p_4= \Pi _4(\phi _ {\omega })$
and $p_j=-\frac 12 v_j  p_4$ for $j\le 3.$
We have $\lambda _4 (p)=-\omega - \frac {v^2}4$ and
$\lambda _j (p)=-v_j$ for $j\le 3.$   Notice that for $\frac{d}{d\omega } Q(\phi _ {\omega }) \neq 0$
this yields \eqref{eq:nondegen}. Notice that

\begin{equation*}
 \begin{aligned} &   \nabla ^2
E(e^{-\frac 12 J v\cdot x}U) = e^{-\frac 12 J v\cdot x}
 \left (  \nabla ^2 E(U)-  J  v\cdot \nabla _x + \frac{v^2}{4}  \right )   e^{\frac 12 J v\cdot x}
\end{aligned}
\end{equation*}
and that $v\cdot \nabla _x \circ  e^{-\frac 12 J v\cdot x} =e^{-\frac 12 J v\cdot x}\circ   (v\cdot \nabla _x -J \frac {v^2}2) $
and
\begin{equation*}
 \begin{aligned} &   \nabla ^2
E(\Phi _p (x) ) -\lambda (p) \cdot \Diamond = e^{-\frac 12 J v\cdot x}
 \left (  \nabla ^2 E((\phi _ {\omega },0))-  J  v\cdot \nabla _x + \frac{v^2}{4}  \right )   e^{\frac 12 J v\cdot x}  \\&
+ Jv\cdot \nabla _x  e^{-\frac 12 J v\cdot x} e^{\frac 12 J v\cdot x} + (\omega +\frac {v^2}4)e^{-\frac 12 J v\cdot x} e^{\frac 12 J v\cdot x}.
\end{aligned}
\end{equation*}
They imply
\begin{equation}\label{eq:conjNLS}
 \begin{aligned} &   {\mathcal H}_p =e^{-\frac 12 J v\cdot x}
 {\mathcal H}_\omega  e^{\frac 12 J v\cdot x} \, , \quad  {\mathcal H}_\omega
:=J
(\nabla ^2 E((\phi _ {\omega },0)) +\omega ).
\end{aligned}
\end{equation}
The multiplier operator $e^{-\frac 12 J v\cdot x}$ is an isomorphism in all spaces $\Sigma _n$ so all the information on the spectrum of
${\mathcal H}_p $ is obtained from the spectrum of ${\mathcal H}_\omega .$
We have   $ {\mathcal H}_\omega  = {\mathcal H}_{0\omega}     +V $ where $ {H}_{0\omega}:=    J     (-\Delta   +\omega )$
and
\begin{equation*}
 \begin{aligned} &    V:= 4 J \begin{pmatrix}   - B' (\phi _\omega ^2 )-2 B ^{\prime  \prime}   (\phi
_\omega ^2) \phi_\omega^2  &
0  \\
0 &   - B' (\phi _\omega ^2 )
 \end{pmatrix}  .
\end{aligned}
\end{equation*}
This yields   $\sigma  _e( {\mathcal H}_\omega) =\sigma   ( {H}_{0\omega}) =(-\infty , -\omega ]\cup [\omega , \infty )$  and that  $\sigma  _p( {\mathcal H}_\omega)$
is finite with finite multiplicities.  The fact that $\sigma  _p( {\mathcal H}_\omega)    $ is in the complement of
$ \sigma _e(\mathcal{H} _{\omega })$  is expected to be  true generically.
Set  ${\mathcal H} ={\mathcal H}_\omega P_c(\omega )$ for
  $P_c(\omega )$  the projection on $X_c({\mathcal H}_\omega)$.

\begin{lemma} \label{lem:a5}  The statement in  (A5) is true.
\end{lemma}
	\proof Notice that $\Sigma  _n$ is invariant by Fourier transform so that \eqref{eq:a71}
	is equivalent to the fact that for  the following multiplier operator
(that is an operator $\psi (x)$ which maps $u\to (\psi u)(x):=\psi (x)u(x)$)
we have
\begin{equation} \label{eq:a71bis}  \begin{aligned}
 &
\quad \text{ $\|     (1+\epsilon ^2+\epsilon ^2|x| ^2)^{ -2}  \| _{B(\Sigma  _n ,\Sigma  _ n ) } \le C_n<\infty$ $\forall$ $|\epsilon |\le 1 $ and $n\in \N$.}
 \end{aligned}
\end{equation}
Similarly  \eqref{eq:a72}  is equivalent to
\begin{align} \label{eq:a72bis}
 &    strong-\lim _{\epsilon\to 0} (1+\epsilon ^2+\epsilon ^2|x| ^2)^{ -2} =1 \text{  in $B(\Sigma  _n ,\Sigma  _ n )$ } \\&  \lim _{\epsilon\to 0}  \|  (1+\epsilon ^2+\epsilon ^2|x| ^2)^{ -2} -1 \|  _{B(\Sigma  _n ,\Sigma  _ {n'} ) }=0
\quad \text{  for any   $n'\in \N$ with $n'<n$. } \nonumber
 \end{align}
Both  \eqref{eq:a71bis}--\eqref{eq:a72bis}  are elementary to check
using  the first definition of $\Sigma _n$ in Sect.\ref{sec:setup},
computing commutators of the multiplier operators with $\partial ^\alpha _x$  and computing elementary bounds on the derivatives of the multipliers.\qed

\begin{lemma} \label{lem:a6}  The statement in  (A6) is true.
\end{lemma}
	\proof  Using the Fourier transformation  like in Lemma \ref{lem:a5},  (A6) is equivalent to
	 the statement that for any  $n\in \N$ and $c>0$ there a $C$ s.t.
	the following multiplier operator satisfies

	\begin{equation*}
		\|     e^{  (1+\epsilon ^2+\epsilon ^2|x| ^2)^{ -2}    J(   \tau _4 - \sum _{j=1}^3 x_j \tau _j  }  ) \| _{B(\Sigma _n,\Sigma _n)}\le C
	\end{equation*}
	  for any $|\tau |\le c$  and  any $|\epsilon |\le 1 $. This too is elementary to check.

\begin{lemma} \label{lem:l7}  The statement in  (L7) is true.
\end{lemma}
	\proof  From $ \sigma  ({\mathcal H} ) =\sigma _e(\mathcal{H} _{\omega })$ we have
$R_{\mathcal H}\in C^\omega (\rho  ({\mathcal H} ), B(L^2,L^2))$.

\noindent We have $R _{ {\mathcal H}_{0\omega}}$ and $R _{ {\mathcal H}_{0\omega}}\partial _{x_j }$ are in $ C^\omega (\rho  ({\mathcal H} ), B(\Sigma _n,\Sigma _n)) $
for any $n\in \N$.   By conjugation by Fourier transform this is equivalent to the statement that for $z\in  \rho  ({\mathcal H}_{0\omega})$ and $i=0,1$, we have

\begin{equation}  \begin{aligned} &     \xi _j ^i
\begin{pmatrix} ( |\xi |^2+\omega  -z)^{-1}   & 0 \\   0 &  - ( |\xi |^2+\omega  +z)^{-1} \end{pmatrix}
\in B(\Sigma _n,\Sigma _n).
\end{aligned}  \nonumber \end{equation}
This is  elementary, using   the first definition of $\Sigma _n$ in Sect.\ref{sec:setup}.

  \noindent  We have for $i=0,1$
 \begin{equation}\label{eq:resolv1}
  \begin{aligned} &    R _{ {\mathcal H } } (z) \partial _{x_j }^i  = R _{ {\mathcal H}_{0\omega}} (z) P_c(\omega )  \partial _{x_j }^i -
    R _{ {\mathcal H}_{0\omega}}(z) VR _{ {\mathcal H } } (z) \partial _{x_j }^i .
 \end{aligned} \end{equation}
From \eqref{eq:resolv1} we derive,   for $\| \ \| = \| \ \| _{B(L^2,L^2)}$.
\begin{equation}\label{eq:resolv11}
  \begin{aligned} &    \| R _{ {\mathcal H } } (z)\partial _{x_j }^i  \|   \le
	 \|  (1+ R _{ {\mathcal H}_{0\omega}}(z) V) ^{-1}\|
	\| R _{ {\mathcal H}_{0\omega}} (z) P_c(\omega )  \partial _{x_j }^i  \|   ,
 \end{aligned}
 \end{equation}
which yields the $n=0$ case.

\noindent From \eqref{eq:resolv1} we derive
 \begin{equation*}\label{eq:resolv2}
  \begin{aligned} &    \| R _{ {\mathcal H } } (z)\partial _{x_j }^i  \| _{B(\Sigma _n,\Sigma _n)}  \le C  \|  R _{ {\mathcal H}_{0\omega}} (z)\partial _{x_j }^i  \| _{B(\Sigma _n,\Sigma _n)}\\& +C
     \|  R _{ {\mathcal H}_{0\omega}} (z) \| _{B(\Sigma _n,\Sigma _n)}    \| \langle x \rangle ^nV \| _{W^{n,\infty}}
       \| R _{ {\mathcal H } } (z) \partial _{x_j }^i \| _{B(H^n,H^n)}.
 \end{aligned}
 \end{equation*}
\noindent The  last factor is bounded. Indeed for $\mathbf{v}=R _{ {\mathcal H } } (z)\partial _{x_j }^i\mathbf{u}$ we have
 \begin{equation*}\label{eq:resolv3}
  \begin{aligned} &    \partial ^\alpha _x \mathbf{v}= R _{ {\mathcal H } } (z)\partial ^\alpha _x\partial _{x_j }^i \mathbf{u}  +R _{ {\mathcal H } } (z) [V,\partial ^\alpha _x]\partial _{x_j }^i\mathbf{u}
 \end{aligned}
 \end{equation*}
 and induction in  $  n$   yields the desired bounds $\|    \mathbf{v} \| _{H^n}\le C \|   \mathbf{u} \| _{H^n} $.
\qed

\bigskip
\textbf{The Nonlinear Dirac Equation.}
Here the unknown $U$ is   $\C^4$-valued, $u^*$ its complex conjugate and   for $m>0$
\begin{equation}\label{Eq:NLDE}
\im U_t -D_m U -Vu+2B'(U\cdot \beta  {U}^* )\beta U=0
\end{equation}
where  we assume  for the moment $V=0$ and  where
$D_m=-\im  \sum _{j=1}^3\alpha_j\partial_{x_j}+m\beta$,
with for $j=1,2,3$
\begin{equation*}
\alpha _j=  \begin{pmatrix}  0 &
\sigma _j  \\
\sigma_j & 0
\end{pmatrix} , \  \beta =  \begin{pmatrix}  I _{\C^2} &
0 \\
0 & -I _{\C^2}
 \end{pmatrix}
 , \end{equation*} \begin{equation*}
\sigma _1=\begin{pmatrix}  0 &
1  \\
1 & 0
 \end{pmatrix} \, ,
\sigma _2= \begin{pmatrix}  0 &
\rmi   \\
-\rmi  & 0
 \end{pmatrix}  \, ,
\sigma _3=\begin{pmatrix}  1 &
0  \\
0 & -1
 \end{pmatrix}.
\end{equation*}
Notice that the symmetry group \eqref{Eq:NLDE} is not Abelian. In \cite {boussaidcuccagna}
there is a symmetry restriction  on the  solutions considered, by  looking only at functions
such that for any   $x\in \R ^3 $ we  have  { $U (-x) =\beta U (
x )$ and $U (-x_1,-x_2,x_3)=S_3U (x_1, x_2,x_3)$} with
$S_3:=\begin{pmatrix}\sigma_3&0\\0&\sigma_3\end{pmatrix}$.
We need to redefine the spaces $\Sigma _n$ in the proof, introducing these symmetries.  This does not affect the proof.

 \noindent There is a unique invariant $Q(U)=\frac 12 \| u \| _{L^2}.$ In this case $\Diamond _1U=U$ for any $u$. Hence  all the changes of variables are diffeomorphism within each space $\Ph ^K$ (or  $\widetilde{\Ph} ^K$).

\noindent   (A5)--(A6) in this case are  elementary. In fact (A5) is unnecessary, (A6) is necessary only for $\epsilon =0$, in which case is trivial. (L7) is necessary only for $i=0$   (given that the only $\Diamond _j$ is the identity)
and can be proved in a way similar to Lemma \ref{lem:l7}.

\textbf{Nonlinear Dirac   Equation with a Potential.}  Pick $V\in {\mathcal S}(\R  ^3, B(\C^4))$ with
$V(x)$ selfadjoint for the scalar product in $\C^4$ for any  $x\in \R ^3$. Then generically  $\sigma _p(D_m  +V) \subset (-m,m)$. Suppose    $\sigma _p(D_m  +V)  =\{  e_0,..., e_{\textbf{n}  } \}$   with $ e_0<...< e_{\textbf{n}  }$.
Then   bifurcation yields corresponding families of  small standing waves $e^{-\im \omega t} \phi _\omega (x)$
of \eqref{Eq:NLDE}.  For generic $V$ the $e_j$ have multiplicity 1. If we focus on
  $e_0$,   for generic smooth  $B'(r )$
there will be a smooth family $\omega \to \phi _\omega $ in $C^\infty (\mathcal{O}, \Sigma _n)$
for any $n$, with $\mathcal{O}$ an open interval one of whose endpoints is  $e_1$. Then it can be shown that for generic $V$
the hypotheses (L1)--(L6) in Sect. are true, as well as all the previous hypotheses.   Indeed in this case, taking
$\omega $ sufficiently close to $e_0$,
we have eigenvalues with $\textbf{e}_j'$ arbitrarily close to $e_j-e_0$. Generically this yields (L4)--(L5).
 The multiplicity of the $\im \textbf{e}_j'$ is 1. We have $\sigma _e(\mathcal{H}_\omega )= (-\infty , -m +|\omega |]
\cup [m-|\omega |, \infty )$.   An eigenvalue    $\lambda $ of $\mathcal{H}_\omega$
is either $\lambda =0$, or $\lambda =\pm \im \textbf{e}_j'$ for some $j$.
This in particular yields (L1)--(L3).

Department of Mathematics and Geosciences,  University of Trieste, via Valerio  12/1  Trieste, 34127  Italy

{\it E-mail Address}: {\tt scuccagna@units.it}

\end{document}